\documentclass{article}

\usepackage[utf8]{inputenc}
\usepackage[a4paper, top=2.5cm, left=2.5cm, right=2.5cm, bottom=3cm]{geometry}

\usepackage{graphicx}%
\usepackage{multirow}%
\usepackage{amsmath,amssymb,amsfonts}%
\usepackage{amsthm}%
\usepackage{mathrsfs}%
\usepackage[title]{appendix}%
\usepackage{xcolor}%
\usepackage{textcomp}%
\usepackage{manyfoot}%
\usepackage{booktabs}%
\usepackage{algorithm}%
\usepackage{algorithmicx}%
\usepackage{algpseudocode}%
\usepackage{listings}%
\usepackage{caption}%
\usepackage{hyperref}%

\DeclareMathOperator{\st}{\textnormal{s.t.}}
\DeclareMathOperator{\diag}{diag}
\DeclareMathOperator{\argmax}{argmax}
\DeclareMathOperator{\convex}{conv}

\newcommand{\norm}[1]{\left\lVert#1\right\rVert}
\newcommand{\sqnorm}[1]{\left\lVert#1\right\rVert_2^2}

\newcommand{\numc}{l}
\newcommand{\numd}{n}
\newcommand{\numu}{m}

\newcommand{\dvec}{\boldsymbol{x}}
\newcommand{\uvec}{\boldsymbol{\xi}}
\newcommand{\uscal}{\xi}
\newcommand{\cvec}{\boldsymbol{c}}
\newcommand{\dmat}{\boldsymbol{A}}
\newcommand{\umat}{\boldsymbol{D}}
\newcommand{\uconst}{\boldsymbol{d}}
\newcommand{\umap}{\boldsymbol{h}}
\newcommand{\U}{\mathcal{U}}
\newcommand{\pmat}{\boldsymbol{P}}
\newcommand{\qvec}{\boldsymbol{q}}
\newcommand{\betavec}{\boldsymbol{\beta}}

\newcommand{\uvecbase}{\boldsymbol{\zeta}}
\newcommand{\uscalbase}{\zeta}

\newcommand{\vvec}{\boldsymbol{v}}
\newcommand{\svec}{\boldsymbol{s}}
\newcommand{\unitvec}{\boldsymbol{e}}
\newcommand{\nullvec}{\boldsymbol{0}}
\newcommand{\unitytary}{\boldsymbol{1}}

\theoremstyle{definition}
\newtheorem{definition}{Definition}
\newtheorem{assumption}{Assumption}

\theoremstyle{plain}
\newtheorem{theorem}{Theorem}
\newtheorem{lemma}[theorem]{Lemma}
\newtheorem{proposition}[theorem]{Proposition}

\title{Designing Tractable Piecewise Affine Policies for Multi-Stage Adjustable Robust Optimization}
\author{%
	Simon Thomä%
	\footnote{Chair of Operations Management, RWTH Aachen, 52072 Aachen, Germany},
	Grit Walther$^*$,
	Maximilian Schiffer%
	\footnote{TUM School of Management \& Munich Data Science Institute, Technical University of Munich, 80333 Munich, Germany}%
	\\
	simon.thomae@om.rwth-aachen.de
}

\date{\today}

\begin{document}
	
	\maketitle
	
	\begin{abstract}
		\noindent
		We study piecewise affine policies for multi-stage adjustable robust optimization (ARO) problems with non-negative right-hand side uncertainty. 
		First, we construct new dominating uncertainty sets and show how a multi-stage ARO problem can be solved efficiently with a linear program when uncertainty is replaced by these new sets. 
		We then demonstrate how solutions for this alternative problem can be transformed into solutions for the original problem. 
		By carefully choosing the dominating sets, we prove strong approximation bounds for our policies and extend many previously best-known bounds for the two-staged problem variant to its multi-stage counterpart. 
		Moreover, the new bounds are - to the best of our knowledge - the first bounds shown for the general multi-stage ARO problem considered. 
		We extensively compare our policies to other policies from the literature and prove relative performance guarantees. 
		In two numerical experiments, we identify beneficial and disadvantageous properties for different policies and present effective adjustments to tackle the most critical disadvantages of our policies. 
		Overall, the experiments show that our piecewise affine policies can be computed by orders of magnitude faster than affine policies, while often yielding comparable or even better results.
	\end{abstract}
	
	\section{Introduction}
	\label{sec:introduction}
	
	In practice, most decision-making problems have to be solved in view of uncertain parameters. 
	In the operations research domain, two fundamental frameworks exist to inform such decisions. 
	The \textit{stochastic optimization} framework captures uncertainty by using probability distributions and aims to optimize an expected objective. 
	Initially introduced in the seminal work of Dantzig  \cite{Dantzig1955}, the framework has been intensively studied and scholars used it to solve a vast variety of problems including production planning \cite{Beebe1968, Schmidt1996}, relief networks \cite{Noyan2016}, expansion planning \cite{Xie2018, Zou2018}, and newsvendor problems \cite{Porteus1990}. 
	While stochastic optimization performs well on many problem classes, finding tractable formulations is oftentimes challenging. 
	Additionally, the data needed to approximate probability distributions might not always be available, and gathering data is often a costly and time-consuming process. 
	The \textit{robust optimization} (RO) framework overcomes many of these shortcomings by capturing uncertainty through distribution-free uncertainty sets instead of probability distributions. 
	By choosing well-representable uncertainty sets, RO often offers computationally tractable formulations that scale well on a variety of optimization problems. 
	In recent years, these favorable properties have led to a steep increase in research interest, see, e.g. \cite{BenTal2002, BenTal2009, Bertsimas2011, Gorissen2015, Yanikoglu2019}. 
	The flexibility of uncertainty sets and scaleability of solution methods also make RO very attractive for applicational purposes and it has been widely applied to many operations management problems \cite{Lu2021}.
	
	In traditional RO, all decisions must be made before the uncertainty realization is revealed. 
	However, oftentimes some decisions can be delayed until after (part of) the uncertainty realization is known in real-world situations. 
	As a consequence, RO may lead to excessively conservative solutions. 
	To remedy this drawback, Ben-Tal et al. \cite{BenTal2004} introduced the concept of \textit{adjustable robust optimization} (ARO) where some decisions can be delayed until the uncertainty realization is (partly) known. 
	In general ARO, uncertainty realizations are revealed over multiple stages and decisions can be made after each reveal. 
	A decision made in stage $t$ can thus be modeled as a function of all uncertainties associated with previous stages $t'\leq t$.
	
	While ARO improves decision-making in theory, ARO is equivalent to RO on some special problem instances, where static decision policies yield optimal adjustable solutions \cite{BenTal2004, Bertsimas2015a}. 
	Using similar arguments to the ones used for the optimality of static solutions, Marandi and den Hertog \cite{Marandi2018} identified conditions where optimal adjustable decisions are independent of some uncertain parameters. 
	In general, static policies do not yield optimal solutions in the ARO setting. 
	Elucidating this, Haddad-Sisakht and Ryan \cite{HaddadSisakht2018} identified a collection of sufficient conditions that imply the suboptimality of static policies and a strict improvement of ARO over RO. %
	
	In general, even the task of finding optimal adjustable solutions in the special case of two-stage ARO proves to be computationally intractable \cite{BenTal2004}. 
	Accordingly, recent works developed many approximation schemes for ARO that often yield good and sometimes even optimal results in practice, see, e.g. \cite{BenTal2004, Bertsimas2010, Iancu2013, SimchiLevi2019}. 
	In the special case with only two decision stages, the first stage decisions are fixed before any uncertain parameters are known and the second stage decisions use full knowledge of the uncertainty realization. 
	Two-stage ARO already has many applications in practice and has widely been studied in the literature, see, e.g. \cite{BenTal2020, Bertsimas2012, Bertsimas2015a, Bertsimas2016, Bertsimas2021, Hanasusanto2015, Housni2021, SimchiLevi2019a}. 
	
	Still, many real-world problems show inherent multi-stage characteristics and cannot be modeled by two-stage ARO. 
	Examples include variants of inventory management \cite{BenTal2009a}, humanitarian relief \cite{BenTal2011}, and facility location \cite{Baron2011}. 
	The transition from two-stage to multi-stage ARO introduces two main challenges. 
	First, multi-stage ARO problems entail \textit{nonanticipativity restrictions} that disallow decisions to utilize future information. 
	Second, many approaches that solve the problem by iteratively splitting the uncertainty space, like scenario trees \cite{Hoeyland2001}, and adaptive partitioning \cite{Bertsimas2016, Postek2016}, grow exponentially in the number of stages. 
	As a consequence, many results found for two-stage ARO do not readily generalize to multi-stage scenarios. 
	Against this background, we design piecewise affine policies for multi-stage ARO that overcome the previously mentioned challenges and extend, although it is not straightforward, many of the best know approximation bounds for two-stage problems to a multi-stage setting.
	
	In the remainder of this section, we formally introduce our problem (\mbox{Section \ref{sec:problem_description}}), discuss closely related work (Section \ref{sec:related_work}), and summarize our main contributions (Section \ref{sec:contributions}).
	
	\subsection{Problem Description}
	\label{sec:problem_description}
	
\begin{figure}[tp]
		\centering
		\includegraphics{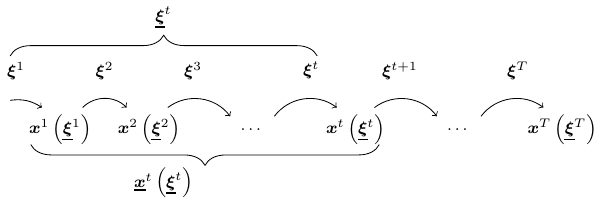}
		\caption{
			Illustration of multi-stage decision making over $T$ stages.
			In each stage $t$ a fraction $\uvec^t$ of the uncertainty is realized and decisions $\dvec^t$ are made.
			Here, decisions $\dvec^t$ may only depend on those uncertainties $\underline{\uvec}^t$ whose realization is known in stage $t$
		}
		\label{fig:multi_stage_decision_making}
	\end{figure}

	In this work we study multi-stage adjustable robust optimization with covering constraints and a positive affine uncertain right hand side.
	Specifically, we consider the following problem:
	\begin{equation}
		\label{eq:problem_formulation}
		\begin{aligned}
			Z_{AR}(\U) =& \min_{\dvec(\uvec)} \max_{\uvec\in\U}  &&\cvec^\intercal \dvec(\uvec)\\
			&\st && \dmat \dvec(\uvec) \geq \umat\uvec + \uconst & \forall \uvec\in\U
		\end{aligned}
	\end{equation}
	with $\dmat\in\mathbb{R}^{\numc\times \numd}$, $\cvec\in\mathbb{R}^{\numd}$, $\umat\in\mathbb{R}^{\numc\times \numu}_+$, $\uconst\in\mathbb{R}^\numc_+$, and compact $\U\subset \mathbb{R}^{\numu}_+$.
	Here, $\numu$ is the number of uncertain parameters, $\numd$ is the number of decisions, and $\numc$ is the number of constraints.
	To model the problem's $T$ stages, we split the uncertainty vector $\uvec$ into $T$ sub-vectors $\uvec = \left(\uvec^1, \dots, \uvec^T\right)$ with $\uvec^t$ being the uncertainty vector realized in stage $t$.
	In the following, we denote by $\underline{\uvec}^t:=\left(\uvec^1, \dots, \uvec^t\right)$ the vector of all uncertainties with known realization in stage $t$.
	Similarly, the adjustable decision vector $\dvec(\uvec)$ divides into $\dvec(\uvec) := \left(\dvec^1(\underline{\uvec}^1), \dots, \dvec^T(\underline{\uvec}^T)\right)$, where the decision $\dvec^t$ made in stage $t$ has to preserve \textit{nonanticipativity} and may only depend on those uncertainties $\underline{\uvec}^t$ whose realization is known in stage $t$.
	We explicitly allow $\uvec^1$ to be zero-dimensional making the initial decision $\dvec^1$ non-adjustable.
	Finally, we denote by $\underline{\dvec}^t(\uvec):=\left(\dvec^1(\underline{\uvec}^1), \dots, \dvec^t(\underline{\uvec}^t)\right)$ the vector of all decisions in the first $t$ stages.
	Figure \ref{fig:multi_stage_decision_making} visualizes the multi-stage decision process with nonanticipativity restrictions.
	
	Unless explicitly stated otherwise, we assume w.l.o.g. that the following assumption holds throughout the paper.
	\begin{assumption}
		\label{ass:uncertainty_assumption}
		$\U\subseteq [0,1]^{m}$ is convex, full-dimensional with $e_i\in\U$ for all $i\in\{1,\dots, \numu\}$, and \textit{down-monotone}, i.e., $\forall \uvec\in\U, \nullvec\leq\uvec'\leq\uvec \colon \uvec'\in\U$.
	\end{assumption}
	Down-monotonicity holds because $\umat,\uconst,\uvec$ are all non-negative, and thus constraints become less restrictive for smaller values of $\uvec$. 
	Convexity holds due to the linearity of the problem, and $\unitvec_i\in\U\subseteq[0,1]^\numu$ holds as $\U$ is compact and $\umat$ can be re-scaled appropriately. 
	We note that the non-negativity assumption of the right-hand side does restrict the problem space. 
	As Bertsimas and \mbox{Goyal \cite{Bertsimas2012}} point out, this assumption prevents the introduction of uncertain or constant upper bounds. 
	However, upper bounds in other decision variables are still possible as $\dmat$ is not restricted, and $\umat, \uconst$ can be zero. 
	Overall, \mbox{Problem (\ref{eq:problem_formulation})} covers a wide range of different problem classes including network design \cite{Ordonez2007, Wang2020}, capacity planning \cite{Laguna1998, Mulvey1995, Paraskevopoulos1991}, as well as versions of inventory management \cite{BenTal2009a, Solyali2016} where capacities are unbounded or subject to the decision makers choice.
	
	In the context of multi-stage decision making, some works require \textit{stagewise} uncertainty, i.e., that the uncertainty set $\U$ consists of uncertainty sets $\U_1,\dots,\U_T$ for each stage, see, e.g., \cite{Bertsimas2010, Bertsimas2011, Georghiou2019}. 
	Like other approaches based on decision rules \cite{Kuhn2011, Bertsimas2015b}, we do not need these restricting assumptions. 
	However, we show how to utilize the existence of such a structure in Section \ref{sec:bounds}.
	
	\subsection{Related Work}
	\label{sec:related_work}
	
	Feige et al. \cite{Feigea2007} show that already the two-stage version of Problem (\ref{eq:problem_formulation}) with $\umat = \boldsymbol{1}, \uconst=\boldsymbol{0}$ and $\dmat$ being a 0-1-matrix is hard to approximate with a factor better than $\Omega(\log \numu)$, even for budgeted uncertainty sets. 
	As it is thereby impossible to find general solutions for $\dvec$, a common technique to get tractable formulations is to restrict the function space.
	
	In this context, Ben-Tal et al. \cite{BenTal2004} consider $\dvec$ to be affine in $\uvec$.
	Specifically, they propose $\dvec^t$ to be of the form ${\dvec^t(\underline{\uvec}^t) = \boldsymbol{P}^t \underline{\uvec}^t + \boldsymbol{q}^t}$.
	Affine policies have been found to deliver good results in practice \cite{Adida2006, BenTal2005} and are even optimal for some special problems \cite{Bertsimas2010, Iancu2013, SimchiLevi2019}.
	Further popular decision rules include
	segregate affine \cite{Chen2008, Chen2009},
	piecewise constant \cite{Bertsimas2011},
	piecewise affine \cite{BenTal2020, Georghiou2015}, and
	polynomial \cite{Bampou2011, Bertsimas2011b} policies, as well as 
	combinations of these \cite{Rahal2021}.
	For surveys on adjustable policies we refer to  Delage and Iancu \cite{Delage2015} and Yan{\i}ko{\u{g}}lu et al. \cite{Yanikoglu2019}.
	
	A key question that arises when using policies to solve ARO problems is how good the solutions are compared to an optimal unrestricted solution. 
	To answer this, many approximation schemes for a priori and a posteriori bounds have been proposed. 
	In the context of a posteriori bounds the focus lies on finding tight upper and lower approximation problems. 
	Hadjiyiannis et al. \cite{Hadjiyiannis2011} estimate the suboptimality of affine decision rules using sample scenarios from the uncertainty set. 
	Similar sample lower bounds are used by Bertsimas and Georghiou \cite{Bertsimas2015b} to bound the performance of piecewise affine policies. 
	Kuhn et al. \cite{Kuhn2011} investigate the optimality of affine policies by using the gap between affine solutions on the primal and the dual of the problem. 
	Georghiou et al. \cite{Georghiou2015} generalize this primal-dual approach to affine policies on lifted uncertainty sets. 
	Building on both of the previous approaches, Georghiou et al. \cite{Georghiou2020} propose a convergent hierarchy of policies that combine affine policies with extreme point scenario samples. 
	Daryalal et al. \cite{Daryalal2022} construct lower bounds by relaxing nonanticipativity and stage-connecting constraints in multi-stage ARO. 
	They then use these lower bounds to construct primal solutions in a rolling horizon manner.
	
	In the context of a priori bounds, most approximation schemes have been proposed for the two-staged version of Problem (\ref{eq:problem_formulation}). 
	For general uncertainty sets on the two-stage version of (\ref{eq:problem_formulation}),  Bertsimas and Goyal \cite{Bertsimas2012} show that affine policies yield an $O(\sqrt{\numu})$ approximation if $\cvec$ and $\dvec$ are non-negative. 
	They further construct a set of instances where this bound is tight, showing that no better general bounds for affine policies exist. 
	Using geometric properties of the uncertainty sets, Bertsimas and Bidkhori \cite{Bertsimas2015} improve on these bounds for some commonly used sets including budgeted uncertainty, norm balls, and intersections of norm balls. 
	Ben-Tal et al. \cite{BenTal2020} propose new piecewise affine decision rules for the two-stage problem that on some sets improve these bounds even further. 
	In addition to strong theoretical bounds, this new approach also yields promising numerical results that can be found by orders of magnitude faster than solutions for affine adjustable policies. 
	For budgeted uncertainty sets and some generalizations thereof, Housni and Goyal \cite{Housni2021} show that affine policies even yield optimal approximations with an \textit{asymptotic bound}, i.e., asymptotic behavior of the approximation bound, see, e.g. \cite{Sedgewick2013}, of $O\left(\frac{\log \numu}{\log\log \numu}\right)$.
	This bound was shown to be tight by Feige et al. \cite{Feigea2007} for reasonable complexity assumptions, namely 3SAT cannot be solved in $2^{O(\sqrt{\numu})}$ time on instances of size $\numu$. 
	We present an overview of known a priori approximation bounds for some commonly used uncertainty sets on two-stage ARO in Table \ref{tab:literature_bounds} of Appendix \ref{sec:literature_bounds}.
	
	To the best of our knowledge, Bertsimas et al. \cite{Bertsimas2011} are the only ones that provide a priori bounds for multi-stage ARO so far. 
	They show these bounds for piecewise constant policies using geometric properties of the uncertainty sets. 
	More specifically, they consider multi-stage uncertainty networks, where the uncertainty realization is taken from one of multiple independent uncertainty sets in each stage. 
	While the choice of the set selected in each stage may depend on the sets selected before, the uncertainty sets are otherwise independent. 
	Although this assumption is fairly general, it still leaves many commonly used sets where uncertainty is dependent over multiple stages uncovered. 
	Among others, uncovered sets include widely used hypersphere and budgeted uncertainty.
	
	As can be seen, previous work has predominantly focused on providing tighter approximation bounds for two-stage ARO, inevitably raising the question of whether similar bounds hold for multi-stage ARO as well. 
	This work contributes towards answering this question by extending many of the currently best-known a priori bounds on two-stage ARO to its multi-stage setting. 
	By so doing, we are - to the best of our knowledge - the first ones to provide a priori approximation bounds for the multi-stage ARO Problem (\ref{eq:problem_formulation}), where uncertainty sets can range over multiple stages. 
	Unless explicitly stated otherwise, we will always refer to a priori approximation bounds when we discuss approximation bounds in the remainder of this paper.
	
	\subsection{Our Contributions}
	\label{sec:contributions}
	
	With this work, we extend the existing literature in multiple ways, where our main contributions are as follows.
	
	\textbf{Tractable Piecewise Affine Policies for Multi-Stage ARO:}
	Motivated by piecewise affine policies for two-stage ARO \cite{BenTal2020}, we present a framework to construct policies that can be used to efficiently find good solutions for the multi-stage ARO Problem (\ref{eq:problem_formulation}). 
	Instead of solving the problem directly for uncertainty $\U$, we first approximate $\U$ by a dominating set $\hat{\U}$. 
	To do so, we define the concept of nonanticipative multi-stage domination and show that this new definition of domination fulfills similar properties to two-stage domination. 
	Based on this new definition, we then construct dominating sets $\hat{\U}$ such that solutions on $\hat{\U}$ can be found efficiently. 
	More specifically, we choose $\hat{\U}$ to be a polytope for which worst-case solutions can be computed by a linear program (LP) over its vertices. 
	In order to ensure nonanticipativity, which is the main challenge of this construction, we introduce a new set of constraints on the vertices that guarantee the existence of nonanticipative extensions from the vertex solutions to the full set $\hat{\U}$. 
	Finally, we show how to use the solution on the dominating set $\hat{\U}$ to construct a valid solution for the original uncertainty set $\U$.

	\textbf{Approximation Bounds:}
	To the best of our knowledge, we provide the first approximation bounds for the multi-stage Problem (\ref{eq:problem_formulation}) with general uncertainty sets. 
	More specifically, we show that our policies yield $O(\sqrt{\numu})$ approximations of fully adjustable policies. 
	While this bound is tight for our type of policies in general, we show that better bounds hold for many commonly used uncertainty sets.
	
	While our main contribution is to extend approximation bounds to multi-stage ARO, \mbox{Problem (\ref{eq:problem_formulation})} is further less restrictive than problems previously discussed in the literature on approximation bounds. 
	In addition to being restricted to two-stage ARO, previous work often assumed $\cvec$ and $\dvec$ to be non-negative \cite{Bertsimas2012, Bertsimas2015, BenTal2020}. 
	\mbox{Ben-Tal et al. \cite{BenTal2020}} additionally restricted parts of $\dmat$, i.e., they require the parts of $\dmat$ associated with the first stage decision to be non-negative. 
	Our policies do not need this assumption. 
	However, we show that Problem (\ref{eq:problem_formulation}) is unbounded whenever there is a feasible $\dvec$ with $\cvec^\intercal \dvec < 0$, due to the non-negativity of the right-hand side. 
	As a consequence, our policies do not readily extend to general maximization problems. 
	
	From a theoretical perspective, mainly asymptotic bounds are of interest. 
	In practice, however, also the exact factors of the approximation are important. 
	Throughout the paper, we thus always give the asymptotic, as well as the exact bounds. 
	We compare all our bounds to the previously best-known bounds for the two-stage setting given in Ben-Tal et al. \cite{BenTal2020} and show that our constructions yield both constant factor, as well as asymptotic improvements. 
	
	\textbf{Comparison with Affine Policies:}
	Using the newly found bounds, we show that no approximation bound for affine policies on hypersphere uncertainty exists that is better than the bound we show for our policies. 
	For budgeted uncertainty, on the other hand, we show that affine policies strictly dominate our piecewise affine policies. 
	These findings confirm results that have been reported for the two-stage variant, where affine policies do not perform well for hypersphere uncertainty \cite{Bertsimas2012}, but very well for budgeted uncertainty \cite{Housni2021}.
	
	\textbf{Improvement Heuristic:}
	Due to inherent properties of our policy construction, resulting solutions are overly pessimistic on instances where the impact on the objective varies significantly between different uncertainty dimensions. 
	To diminish this effect, we introduce an improvement heuristic that performs at least as well as our policies and that can be integrated into the LP used to construct our policies. 
	While these modifications come at the cost of higher solution times, they allow for significant objective improvements on some instance classes. 
	
	\textbf{Tightening Piecewise Affine Policies via Lifting:}
	We show that in the context of Problem (\ref{eq:problem_formulation}) the piecewise affine policies via lifting presented by Georghiou et \mbox{al. \cite{Georghiou2015}} yield equivalent solutions to affine policies. 
	To prevent this from happening, we construct tightened piecewise affine policies via lifting using insights from our piecewise affine policies. 
	These new policies integrate the approximative power of affine policies, and our piecewise affine policies and are guaranteed to perform at least as well as the individual policies they combine.
	
	\textbf{Numerical Evidence:}
	Finally, we present two sets of numerical experiments showing that our policies solve by orders of magnitude faster than the affine adjustable policies presented by Ben-Tal et al. \cite{BenTal2004}, the piecewise affine policies via lifting presented by Georghiou et al. \cite{Georghiou2015}, and the near-optimal piecewise affine policies by Bertsimas and Georghiou \cite{Bertsimas2015b}, while often yielding comparable or improving results. 
	First, we study a slightly modified version of the tests presented in Ben-Tal et al. \cite{BenTal2020}, allowing us to demonstrate our policies' scalability and the impact of our improvement heuristic. 
	Second, we focus on demand covering instances to demonstrate good performances of our policies for a problem that resembles a practical application. 
	We refer to our git repository (\href{https://github.com/tumBAIS/piecewise-affine-ARO}{https://github.com/tumBAIS/piecewise-affine-ARO}) for all material necessary to reproduce the numerical results outlined in this paper.
	
\begin{figure}
		\centering
		\includegraphics{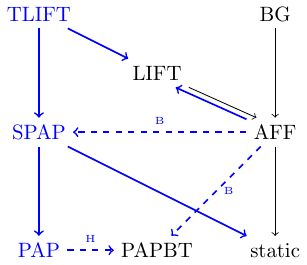}
		\caption{Relations between multi-stage ARO policies compared in this paper.
			An arc from a policy $P$ to another policy $P'$ states $Z_P \leq Z_{P'}$, where $Z_P$ and $Z_{P'}$ are optimal objective values for the ARO Problem (\ref{eq:problem_formulation}) solved with policy $P, P'$ respectively. 
			Dashed arcs only hold for hypersphere (H) or budgeted (B) uncertainty. 
			Relations proved for the first time in this paper are highlighted (blue, bold).
			The compared policies are:
			static policies (static);
			piecewise affine policies via domination by Ben-Tal et al. \cite{BenTal2020} (PAPBT);
			our piecewise affine policies via domination (PAP), c.f., Sections \ref{sec:framework} and \ref{sec:bounds};
			affine policies \cite{BenTal2004} (AFF);
			our piecewise affine policies with rescaling (SPAP), c.f., Section \ref{sec:rescaling};
			near-optimal piecewise affine policies \cite{Bertsimas2015b} (BG);
			piecewise affine policies via lifting \cite{Georghiou2015} (LIFT);
			our tightened piecewise affine policies via lifting (TLIFT), c.f., Section \ref{sec:lifing_combination}}
		\label{fig:relation_graph}
	\end{figure}

	\textbf{Comparison Against Closely Related Work:}
	Compared to the closely related work by Ben-Tal et al. \cite{BenTal2020}, who first introduced the concept of domination in the context of ARO, our contributions are multifold. 
	First, we extend domination-based piecewise affine policies to a wider class of problems by switching from a two-stage to a multi-stage setting and relaxing assumptions. 
	We discuss the structural reasons that make this extension non-trivial at the beginning of Section \ref{sec:framework}. 
	In addition to showing stronger approximation guarantees for our policies, we conduct comprehensive theoretical and numerical comparisons with other adaptable policies. 
	Based on these comparisons, we construct two new policies that mitigate weaknesses of domination and integrate its strength with the strength of other policies. 
	More specifically, the first policy integrates finding a good outer approximation of the uncertainty set in the optimization process. 
	The second policy integrates structural results from domination into lifting \mbox{policies, c.f., \cite{Georghiou2015}}. 
	As a result, we get a hierarchy of piecewise affine adjustable policies with provable relative performance guarantees. 
	We give an overview of all policies constructed in our work and their relative performance guarantees compared to other policies in Figure \ref{fig:relation_graph}.
	
	\vspace{\baselineskip}
	The rest of this paper is structured as follows. 
	In Section \ref{sec:framework}, we introduce our policies and elaborate on their construction. 
	In Section \ref{sec:bounds}, we present our approximation bounds for the multi-stage ARO \mbox{Problem (\ref{eq:problem_formulation})}. 
	We present an improvement heuristic for our policies in Section \ref{sec:rescaling}. 
	By using the results of Sections~\ref{sec:framework} and~\ref{sec:bounds}, we construct tightened piecewise affine policies via lifting in Section \ref{sec:lifing_combination}. 
	Finally, we provide numerical evidence for the performance of our policy compared to other state of the art policies in \mbox{Section \ref{sec:experiments}}. 
	\mbox{Section \ref{sec:conclusion}} concludes this paper with a brief reflection of our work and avenues for future research. 
	To keep the paper concise, we defer proofs that could possibly interrupt the reading flow to Appendices \ref{sec:proof_bounded_domination}-\ref{sec:proof_tightenedlifting}.

	\section{Framework For Piecewise Affine Multi-Stage Policies}
	\label{sec:framework}
	
	In this section, we present our piecewise affine framework for the multi-stage ARO Problem (\ref{eq:problem_formulation}). 
	The main rationale of our framework is to construct new uncertainty sets $\hat{\U}$ that dominate the original uncertainty sets $\U$. 
	With this, our framework follows a similar rationale as the two-stage framework from Ben-Tal et al. \cite{BenTal2020}.
	For a problem $Z_{AR}(\U)$ we construct $\hat{\U}$ in such a way that $Z_{AR}(\hat{\U})$ can be efficiently solved, and a solution of $Z_{AR}(\hat{\U})$ can be used to generate solutions for $Z_{AR}(\U)$.
	
	In this context, we note that one cannot straightforwardly apply the construction scheme used by  Ben-Tal et al. \cite{BenTal2020} due to nonanticipativity requirements. 
	More specifically,  Ben-Tal et al. \cite{BenTal2020} construct $\hat{\U}$ as polytopes, where it is well known that worst-case solutions always occur on extreme points, as any solution can be represented by convex combinations of extreme point solutions. 
	The construction of these convex combinations, however, is not guaranteed to be nonanticipative in the multi-stage setting. 
	To overcome this challenge, we incorporate nonanticipativity in the concept of uncertainty set domination and extend it to a multi-stage setting.
	
	\begin{definition}[Domination]
		\label{def:domination}
		Given an uncertainty set $\U\subseteq \mathbb{R}^{\numu}_+$, we say that $\U$ is \textit{dominated} by $\hat{\U}\subseteq \mathbb{R}^{\numu}_+$ if there is a \textit{domination function} $\umap\colon\U \to \hat{\U}$ with $\umap(\uvec)\geq\uvec$, and $\umap$ can be expressed as $\umap(\uvec) = \left(\umap^1(\underline{\uvec}^1), \dots, \umap^T(\underline{\uvec}^T)\right)$ where $\umap^t$ maps to the uncertainties in stage $t$ and depends on uncertainties up to that stage.
	\end{definition}
	
\begin{figure}[tp]
		\centering
		\includegraphics{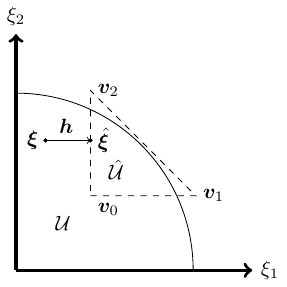}
		\caption{Two dimensional hypersphere uncertainty set $\U$ with (dashed) dominating set \mbox{$\hat{\U}$ (\ref{eq:dominating_U})} induced by the convex combination of vertices $\vvec_0,\vvec_1,\vvec_2$ and dominating function $\umap$ (\ref{eq:domination_function}) that maps a point $\uvec\in\U$ to a point $\hat{\uvec}\in\hat{\U}$}
		\label{fig:dominating_illustration}
	\end{figure}

	Intuitively, an uncertainty set $\hat{\U}$ dominates another set $\U$ if for every point $\uvec\in\U$ there is a point $\hat{\uvec}\in\hat{\U}$ that is at least as large in each component, i.e., $\hat{\uvec}\geq\uvec$. 
	Later we show that the dominating set $\hat{\U}$ can be constructed as the convex combination of $\numu+1$ vertices $\vvec_0,\dots,\vvec_\numu$.
	We also show how to construct dominating functions $\umap$ for these vertex induced dominating sets. 
	Figure \ref{fig:dominating_illustration} illustrates the hypersphere uncertainty set $\U=\left\{\uvec\in\mathbb{R}_+^{\numu} \middle\vert \sqnorm{\uvec}\leq 1\right\}$ together with our dominating set \mbox{$\hat{\U}$ (\ref{eq:dominating_U})} and the dominating function \mbox{$\umap$ (\ref{eq:domination_function})} for $\numu=2$.
	
	Due to the non-negativity of the problem's right-hand side, domination at most restricts the set of feasible solutions. 
	As a consequence, each feasible solution for a realization $\hat{\uvec}\in\hat{\U}$ is also a feasible solution for all realizations $\uvec\in\U$ that are dominated by $\hat{\uvec}$. 
	Using this property, we can derive piecewise affine policies for $Z_{AR}(\U)$ from solutions of $Z_{AR}(\hat{\U})$. 
	Since $\U$ is full-dimensional and down-monotone by \mbox{Assumption \ref{ass:uncertainty_assumption}}, there always exists a factor $\beta\geq0$ such that scaling $\U$ by $\beta$ contains $\hat{\U}$. 
	\mbox{Theorem \ref{thm:bounded_domination}} shows that with this factor $\beta$, solutions of problem $Z_{AR}(\hat{\U})$ are $\beta$-approximations for problem $Z_{AR}(\U)$. 
	It also shows that $Z_{AR}(\hat{\U})$ is unbounded exactly when $Z_{AR}(\U)$ is unbounded. 
	Thus, we assume for the remainder of this paper w.l.o.g. that both $Z_{AR}(\U)$ and $Z_{AR}(\hat{\U})$ are bounded.
	
	\begin{theorem}
		\label{thm:bounded_domination}
		Consider an uncertainty set $\U$ from Problem (\ref{eq:problem_formulation}) and a dominating set $\hat{\U}$.
		Let $\beta\geq1$ be such that $\forall\hat{\uvec}\in\hat{\U} \colon \frac{1}{\beta}\hat{\uvec}\in\U$.
		Moreover, let $Z_{AR}(\U)$ and $Z_{AR}(\hat{\U})$ be optimal values of Problem (\ref{eq:problem_formulation}).
		Then, either $Z_{AR}(\U)$ and $Z_{AR}(\hat{\U})$ are unbounded or
		$$
		0\leq Z_{AR}(\U) \leq Z_{AR}(\hat{\U}) \leq \beta \cdot Z_{AR}(\U).
		$$
	\end{theorem}%
	\noindent
	We present the proof for Theorem \ref{thm:bounded_domination} in Appendix \ref{sec:proof_bounded_domination}.
	
	In the remainder of this section, we demonstrate how the results of \mbox{Theorem \ref{thm:bounded_domination}} can be used to efficiently construct $\beta$-approximations for $Z_{AR}(\U)$. 
	Therefore, we show in Section \ref{sec:construction_dom_set} how to construct dominating polytopes $\hat{\U}$ and efficiently find solutions $Z_{AR}(\hat{\U})$ that comply with nonanticipativity requirements. 
	Then, we construct the dominating function $\umap\colon\U \to \hat{\U}$, which allows us to extend these solutions to solutions for $Z_{AR}(\U)$ in Section \ref{sec:construction_dom_fct}.
	
	\subsection{Construction of the Dominating Set}
	\label{sec:construction_dom_set}
	
	In the following, we construct a dominating set in the form of a polytope for which the worst-case solution can be efficiently found by solving a linear program on its vertices. 
	Specifically, for an uncertainty set $\U$, we consider dominating sets $\hat{\U}$ of the form
	\begin{equation}
		\hat{\U} := \convex(\vvec_0, \vvec_1, \dots, \vvec_{\numu})
		\label{eq:dominating_U}
	\end{equation}
	where for all $i\in\{0,\dots, \numu\}: \frac{1}{\beta} \vvec_i\in\U$ and for all $i\in\{1,\dots, \numu\}: \vvec_i = \vvec_0 + \rho_i \unitvec_i$ for some $\rho_i\in \mathbb{R}_+$. 
	We postpone the construction of the domination function $\umap$, the base vertex $\vvec_0$, and parameters $\rho_1,\dots, \rho_{\numu}$ to \mbox{Section \ref{sec:construction_dom_fct}} and first focus on the construction of solutions for $Z_{AR}(\hat{\U})$. 
	Here, we extend the notation on $\dvec$ and $\uvec$ introduced in \mbox{Section \ref{sec:problem_description}} to $\dvec_i$ and $\vvec_i$. 
	Consequently, $\dvec_i^t$ is the sub-vector of $\dvec_i$ corresponding to decisions made in stage $t$, and $\underline{\vvec_i}^t$ is the sub-vector of $\vvec_i$ corresponding to uncertainties up to stage $t$. 
	Then, the key component for our construction is LP (\ref{eq:lp})
	\begin{subequations}
		\label{eq:lp}
		\begin{align}
			Z_{LP}(\hat{\U}) = \min_{\dvec_0,\dots, \dvec_{\numu}} & z \label{eq:lpobj}\\
			\text{s.t.} \quad & z \geq \cvec^\intercal \dvec_i & \forall i \in \{0,\dots, \numu\} \label{eq:lpconstr_max}\\
			& \dmat \dvec_i \geq \umat \vvec_i + \uconst &  \forall i \in \{0,\dots, \numu\} \label{eq:lpconstr_constr}\\
			& \dvec_i^t = \dvec_j^t &
			\hspace{-.5cm}\forall i,j \in \{0,\dots, \numu\},
			t\in\{1,\dots, T\},
			\underline{\vvec_i}^t = \underline{\vvec_j}^t \label{eq:lpconstr_nonanticipativity}.
		\end{align}
	\end{subequations}
	
	Intuitively, the Objective (\ref{eq:lpobj}) together with Constraints (\ref{eq:lpconstr_max}) minimize the maximal cost over all vertex solutions $\dvec_i$. 
	Constraints (\ref{eq:lpconstr_constr}) ensure that each $\dvec_i$ is a feasible solution for the respective uncertainty vertex $\vvec_i$ of $\hat{\U}$. 
	Finally, Constraints (\ref{eq:lpconstr_nonanticipativity}) ensure nonanticipativity by forcing vertex solutions to be equal unless different uncertainties were observed. 
	With these constraints, we construct LP (\ref{eq:lp}) such that it is sufficient to find an optimal solution for $Z_{LP}(\hat{\U})$ in order to find an optimal solution for $Z_{AR}(\hat{\U})$. 
	
	\begin{lemma}
		\label{lem:lp_sol}
		Let $\hat{\U}$ be a dominating set as described in (\ref{eq:dominating_U}), $Z_{LP}(\hat{\U})$ be the solution of LP (\ref{eq:lp}), and $Z_{AR}(\hat{\U})$ be the solution of Problem (\ref{eq:problem_formulation}).
		Then the LP solution $(\dvec_i)$ on the vertices of $\hat{\U}$ can be extended to a solution on the full set $\hat{\U}$ and we find:
		$$
		Z_{LP}(\hat{\U}) = Z_{AR}(\hat{\U}).
		$$
	\end{lemma}
	\noindent
	We present the proof for Lemma \ref{lem:lp_sol} in Appendix \ref{sec:proof_lp_sol}.
	
	\subsection{Construction of the Domination Function}
	\label{sec:construction_dom_fct}
	
	In the previous section, we showed how to construct dominating sets $\hat{\U}$ such that $Z_{AR}(\hat{\U})$ can be solved efficiently. 
	In order for $\hat{\U}$ to be a valid dominating set for some uncertainty set $\U$, we additionally have to construct a nonanticipative dominating function $\umap\colon\U \to \hat{\U}$ according to Definition \ref{def:domination}. 
	Specifically, we use
	\begin{equation}
		\label{eq:domination_function}
		\umap(\uvec) := \left(\uvec - \vvec_0\right)_+ + \vvec_0
	\end{equation}
	where $(\cdot)_+$ is the element-wise maximum with $0$ and $\vvec_0$ is the base vertex from Definition (\ref{eq:dominating_U}). 
	By construction, $\umap$ maps each uncertainty realization $\uvec$ to its element-wise maximum with $\vvec_0$. 
	It directly follows that $\umap$ is nonanticipative, as each element in $\umap(\uvec)$ solely depends on the corresponding element in $\uvec$.
	
	Finally, we have to ensure that $\umap(\uvec)\in\hat{\U}$ for all $\uvec\in\U$. 
	We do so by choosing the base vertex $\vvec_0$ and parameters $\rho_1,\dots, \rho_{\numu}$ during the construction of $\hat{\U}$ appropriately. 
	Using $\lambda_i(\uvec) := \frac{\left( \left( \uvec - \vvec_0 \right)_+\right)_i}{\rho_i}$ with the convention $\frac{0}{0}=0$, we find
	$$
	\umap(\uvec) = \sum_{i=1}^{\numu} \lambda_i(\uvec) \vvec_i + \left(1-\sum_{i=1}^{\numu} \lambda_i(\uvec)\right) \vvec_0.
	$$
	By definition, any convex combination of $\vvec_i$ is contained in $\hat{\U}$. 
	Thus, $\umap$ is a valid domination function if and only if
	\begin{equation}
		\label{eq:validity_criterion}
		\max_{\uvec\in\U}\sum_{i=1}^{\numu} \lambda_i(\uvec) \leq 1.
	\end{equation}
	Condition (\ref{eq:validity_criterion}) gives a compact criterion to check the validity of dominating sets. 
	By doing so, it lays the basis for our optimal selection of the base vertex $\vvec_0$ and parameters $\rho_0, \dots, \rho_\numu$. 
	Checking Condition (\ref{eq:validity_criterion}) generally requires solving a convex optimization problem. 
	However, in Section \ref{sec:bounds} we show that for many commonly used special uncertainty sets, this problem can be significantly simplified, leading to low dimensional unconstrained minimization problems or even analytical solutions.
	
	Recall that we showed how to extend a solution $(\dvec_0,\dots,\dvec_\numu)$ of LP (\ref{eq:lp}) to the full set $\hat{\U}$ in the proof of Lemma \ref{lem:lp_sol}. 
	Combining this with $\umap$ and using Theorem \ref{thm:bounded_domination} we get a piecewise affine solution for $Z_{AR}(\U)$ by
	\begin{equation}
		\label{eq:piecewise_policy}
		\dvec(\uvec) = \sum_{i=1}^{\numu} \lambda_i(\uvec) \dvec_i + \left(1-\sum_{i=1}^{\numu} \lambda_i(\uvec)\right) \dvec_0
	\end{equation}
	that has an optimality bound of $\beta$.
	
	\subsection{Limitations}
	\label{sec:limitations}
	While our policies overcome the nontrivial challenge of nonanticipativity on extreme point solutions, they still rely on the ability to form convex combinations. 
	As integrality is not preserved by convex combinations, there is no natural way to extend our approach to integer or binary recourse decisions $\dvec$. 
	However, including non-adjustable integer or binary first-stage decisions in our framework is straightforward. 
	Also, it is not straightforwardly possible to incorporate uncertain recourse decisions, i.e., dependence of $\dmat$ on $\uvec$, into our approach, as worst-case realizations for problems with uncertain recourse are not necessarily extreme points of $\U$, see, e.g., \cite{Ayoub2016, Georghiou2020}.
	
	\section{Optimality Bounds for Different Uncertainty Sets}
	\label{sec:bounds}
	
	In the previous section, we demonstrated how to construct nonanticipative piecewise affine policies for the multi-stage Problem (\ref{eq:problem_formulation}). 
	On this basis, proving approximation bounds mainly depends on geometric properties of the uncertainty sets $\U$. 
	We first show approximation bounds for some commonly used permutation invariant uncertainty sets. 
	On these sets, the dominating sets are permutation invariant and we give closed-form constructions. 
	We then give approximation bounds for our piecewise affine policies on general uncertainty sets. 
	Finally, we demonstrate how the bounds of an uncertainty set $\U$ generalize to transformations of that set. 
	While in theory, mostly asymptotic bounds are of interest, in practice constant factors are important as well. 
	Thus we always state exact, as well as asymptotic bounds. 
	Table \ref{tab:optimality_bounds} gives an overview of all bounds that are explicitly proven in Propositions \ref{pro:hypersphere}, \ref{pro:budgeted}, \ref{pro:pball}, \ref{pro:ellipsoid} and \ref{pro:general} of this section. 
	We compare all our results against the results for the two-stage setting in Ben-Tal et al. \cite{BenTal2020} and show constant factor, as well as asymptotic improvements. 
	For budgeted and hypersphere uncertainty sets, we further compare the theoretical performance of our piecewise affine policies with affine adjustable policies.
	
	\begin{table}[tb]
		\caption{
			Performance bounds of the piecewise affine policy for different uncertainty sets}
		\label{tab:optimality_bounds}
		\newcounter{uncertaintyexampleid}
		\renewcommand{\theuncertaintyexampleid}{\Roman{uncertaintyexampleid}}
		\newcommand{\uncertaintyexample}{{\stepcounter{uncertaintyexampleid}\theuncertaintyexampleid}}%
		\setcounter{uncertaintyexampleid}{0}
			\centering
			\begin{tabular}{llll}
				\hline
				No. & Uncertainty Set $\U$ & Bound & Asymptotic Bound\\
				\hline
				\hline
				\uncertaintyexample
				& $\left\{\uvec\in\mathbb{R}_+^{\numu} \middle\vert \sqnorm{\uvec}\leq 1\right\}$ 
				& $\sqrt[4]{\numu}\frac{\sqrt{\numu-1}}{\sqrt{2\left(\numu-\sqrt{\numu}\right)}}$
				& $O(\sqrt[4]{\numu})$\\
				\uncertaintyexample
				& $\left\{\uvec\in[0,1]^{\numu} \middle\vert \norm{\uvec}_1\leq k\right\}$ 
				& $\frac{k(\numu-1)}{\numu+k(k-2)}$ 
				& $O\left(\min\left\{k, \frac{\numu}{k}\right\}\right)$ \\
				\uncertaintyexample
				& $\left\{\uvec\in\mathbb{R}_+^{\numu} \middle\vert \norm{\uvec}_p\leq 1\right\}$
				& $\left(2(\numu-1) + 2^p\right)^{\frac{p-1}{p^2}}p^{\frac{1-p}{p}}\left(p-1\right)^{\left(\frac{p-1}{p}\right)^2}$ 
				& $O\left(\numu^\frac{p-1}{p^2}\right)$ \\
				\uncertaintyexample
				& $\left\{\uvec\in\mathbb{R}_+^{\numu} \middle\vert \uvec^\intercal \boldsymbol{\Sigma}\uvec \leq 1\right\}$
				& $\begin{cases}
					\sqrt{\frac{1}{2} + \frac{a\numu + \sqrt{(1-a)\numu+a\numu^2}}{2(1-a)}} &  a \leq \numu^{-\frac{2}{3}} \\
					\frac{1}{\sqrt{a}} & a > \numu^{-\frac{2}{3}}
				\end{cases} $
				& $O\left(\numu^{\frac{1}{3}}\right)$ \\
				\uncertaintyexample
				& $\U\subset\mathbb{R}_+^\numu$
				& $2\sqrt{\numu}+1$ 
				& $O\left(\sqrt{\numu}\right)$ \\
				\hline
			\end{tabular}%
		
		\setcounter{uncertaintyexampleid}{0}
		\begin{minipage}{\textwidth}\footnotesize
			We prove specific bounds for uncertainty sets of the forms 
			\uncertaintyexample) hypersphere uncertainty; 
			\uncertaintyexample) budgeted uncertainty; 
			\uncertaintyexample) $p$-norm ball uncertainty, with $p\geq1$;
			\uncertaintyexample) ellipsoid uncertainty, with ${\boldsymbol{\Sigma}:=(1-a)\boldsymbol{1}+a\boldsymbol{J}}$ where $\boldsymbol{1}$ is the unity matrix and $\boldsymbol{J}$ the matrix of all ones;
			\uncertaintyexample) general uncertainty sets.%
		\end{minipage}
	\end{table}
	
	For permutation invariant uncertainty sets there exists an optimal choice of $\hat{\U}$ that is also permutation invariant. 
	More specifically, $\vvec_i$ simplifies to $\vvec_0 = \mu \unitvec$ and ${\vvec_i = \vvec_0+\rho \unitvec_i}$, for some $\mu, \rho$. 
	\begin{lemma}
		\label{lem:permutation_invariant}
		Let $\U$ be a permutation invariant uncertainty set.
		Then there exist $\mu, \rho$ such that the dominating uncertainty set $\hat{\U}$ spanned by $\vvec_0 = \mu \unitvec$ and ${\vvec_i = \vvec_0+\rho \unitvec_i}$ for $i\in\{1,\dots,\numu\}$ there is no other dominating set $\hat{\U}'$ constructed as in (\ref{eq:dominating_U}) with a smaller approximation factor $\beta$.
	\end{lemma}
	\noindent
	We present the proof for Lemma \ref{lem:permutation_invariant} in Appendix \ref{sec:proof_lem_permutation_invariant}.
	
	With these simplifications, Condition (\ref{eq:validity_criterion}) becomes
	\begin{equation}
		\label{eq:rot_inv_validity_criterion}
		\frac{1}{\rho}\max_{\uvec\in\U}\sum_{i=1}^\numu\left(\uscal_i-\mu\right)_+ \leq 1.
	\end{equation}
	With the permutation invariance of the problem, Ben-Tal et al. \cite{BenTal2020} show that for any $\mu$ there exists a $j\leq\numu$, such that the maximization problem in (\ref{eq:rot_inv_validity_criterion}) has a solution that is constant on the first $j$ components and zero on all others components. 
	
	\begin{lemma}[Lemma 4 in  Ben-Tal et al. \cite{BenTal2020}]
		\label{lem:rotational_invariate_worst_case}
		Let $\gamma(j)$ be the maximal average value of the first $j$ components of any $\uvec\in\U$
		$$
		\gamma(j) := \frac{1}{j}\max_{\uvec\in\U}\sum_{i=1}^j \uscal_i.
		$$
		Then for each $\mu$ there exists an optimal solution $\uvec^*$ for the maximization problem in Equation (\ref{eq:rot_inv_validity_criterion}) that has the form
		$$
		\uvec^* = \sum_{i=1}^j \gamma(j)\unitvec_i,
		$$
		for some $j\leq\numu$.
	\end{lemma}
	
\begin{figure}[ht]
		\centering
		\includegraphics[width=\textwidth]{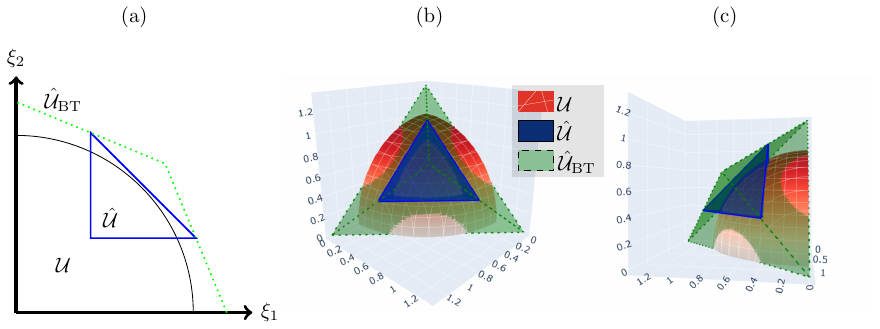}
		\caption{Comparison of our dominating set $\hat{\U}$ (blue, solid frame) and the dominating set $\hat{\U}_\text{BT}$ proposed in  Ben-Tal et al. \cite{BenTal2020} (green, dashed frame) for the hypersphere uncertainty set $\U$ in $\numu = 2$ (a) and $\numu=3$ (b), (c) uncertainty dimensions}
		\label{fig:hypersphere_comparison}
	\end{figure}

	\textbf{Hypersphere Uncertainty:}
	We first use Lemma \ref{lem:rotational_invariate_worst_case} to find a new dominating set for hypersphere uncertainty. 
	By doing so, we find a new approximation bound that improves the bound of $\sqrt[4]{\numu}$ provided in  Ben-Tal et al. \cite{BenTal2020} by a factor of $\sqrt{\frac{\sqrt{m}+1}{2\sqrt{m}}}$, which for large $\numu$ converges towards $\frac{1}{\sqrt{2}}$. 
	While this improvement is irrelevant for the asymptotic complexity of the problem, the new formulation of $\hat{\U}$ does make a difference in practice. 
	In Figure \ref{fig:hypersphere_comparison} we illustrate the improvement of our dominating set $\hat{\U}$ over the dominating set $\hat{\U}_{\text{BT}}$ proposed in  Ben-Tal et al. \cite{BenTal2020} for hypersphere uncertainty sets in two and three uncertainty dimensions. 
	Note, that our sets $\hat{\U}$ are fully contained in the sets $\hat{\U}_{\text{BT}}$, and all extreme points of $\hat{\U}_{\text{BT}}$ are located outside of $\hat{\U}$.
	This implies $Z_{AR}(\hat{\U})\leq Z_{AR}(\hat{\U}_{\text{BT}})$ for hypersphere uncertainty. 
	The formal proof follows from straightforward convex containment and is left for brevity.
	
	\begin{proposition}[Hypersphere]
		\label{pro:hypersphere}
		Consider the hypersphere uncertainty set $\U=\left\{\uvec\in\mathbb{R}_+^{\numu} \middle\vert \sqnorm{\uvec}\leq 1\right\}$.
		Then a solution for $Z_{AR}(\hat{\U})$ where $\hat{\U}$ is constructed using Criterion (\ref{eq:rot_inv_validity_criterion}) with
		\begin{align*}
			\mu &= \frac{1}{2\sqrt[4]{\numu}}, &
			\rho &= \frac{\sqrt[4]{\numu}}{2},
		\end{align*}
		gives a $\beta = \sqrt{\frac{\sqrt{m}+1}{2}}$ approximation for problem $Z_{AR}(\U)$.
	\end{proposition}
	\noindent
	We present the proof for Proposition \ref{pro:hypersphere} in Appendix \ref{sec:proof_hypersphere}.
	
	We can also use this improved performance bound to show that affine adjustable policies cannot yield better bounds than piecewise affine adjustable policies for $\numu\geq 153$. 
	This is because there are instances of Problem (\ref{eq:problem_formulation}) with hypersphere uncertainty where affine adjustable policies perform at least $\frac{4}{5}\left(\sqrt[4]{\numu} - \frac{1}{\sqrt[4]{\numu}}\right)$ worse than an optimal policy. 
	We formalize these results in Proposition \ref{pro:performance_aff_ub}. 
	Note, that these better bounds do not imply that piecewise affine adjustable policies always yield better results than affine adjustable policies for hypersphere uncertainty.
	
	\begin{proposition}
		\label{pro:performance_aff_ub}
		Affine adjustable policies cannot achieve better performance bounds than $\frac{4}{5}\left(\sqrt[4]{\numu} - \frac{1}{\sqrt[4]{\numu}}\right)$ for Problem (\ref{eq:problem_formulation}) with hypersphere uncertainty, even for $\cvec, \dvec, \dmat$ being non negative and $\dmat$ being a 0,1-matrix.
	\end{proposition}
	\noindent
	We present the proof for Proposition \ref{pro:performance_aff_ub} in Appendix \ref{sec:proof_performance_aff_ub}.
	
	\textbf{Budgeted Uncertainty:}
	Next, we tighten the bounds for budgeted uncertainty sets.
	Proposition \ref{pro:budgeted} shows that our new bound is given by $\beta = \frac{k(\numu-1)}{\numu+k(k-2)}$.
	Using
	\begin{equation*}
		\begin{aligned}
			&\frac{\beta}{k} &=& \frac{k(\numu-1)}{k(\numu+k(k-2))} = \frac{\numu-1}{\numu-1+k^2-2k+1} = \frac{\numu-1}{\numu-1+(k-1)^2} \leq 1 \\
			&\text{and}\span\span\\
			&\frac{\beta}{\frac{\numu}{k}} &=& \frac{k^2(\numu-1)}{\numu(\numu+k(k-2))} = \frac{k^2(\numu-1)}{k^2(\numu-1) + k^2 + \numu^2 - 2 k \numu} \\
			&&=& \frac{k^2(\numu-1)}{k^2(\numu-1) + (\numu-k)^2} \leq 1,
		\end{aligned}
	\end{equation*}
	we show $\beta\leq\min(k,\frac{\numu}{k})$, which matches the bound for the two-stage problem variant in  Ben-Tal et al. \cite{BenTal2020}. 
	As $\frac{\beta}{k}$ is decreasing in $k$ and $\frac{\beta}{\frac{\numu}{k}}$ is increasing in $k$, we obtain a maximum improvement for $k=\frac{\numu}{k} \Leftrightarrow k=\sqrt{\numu}$. 
	At this point the improvement of the bound reaches a factor of $\frac{1}{2}$.
	
	\begin{proposition}[Budget]
		\label{pro:budgeted}
		Consider the budgeted uncertainty set $\U=\left\{\uvec\in[0,1]^{\numu} \middle\vert \norm{\uvec}_1\leq k\right\}$ for some $k\in\{1,\dots, \numu\}$. 
		Then a solution for $Z_{AR}(\hat{\U})$ where $\hat{\U}$ is constructed using Criterion (\ref{eq:rot_inv_validity_criterion}) with
		\begin{align*}
			\mu &= \frac{k(k-1)}{\numu+k(k-2)}, &
			\rho &= \frac{k(\numu-k)}{\numu+k(k-2)},
		\end{align*}
		gives a $\beta = \frac{k(\numu-1)}{\numu+k(k-2)}$ approximation for problem $Z_{AR}(\U)$.
	\end{proposition}
	\noindent
	We present the proof for Proposition \ref{pro:budgeted} in Appendix \ref{sec:proof_budgeted}.
	
	Note, that there is no result analogous to Proposition \ref{pro:performance_aff_ub} for budgeted uncertainty as Housni and Goyal \cite{Housni2021} showed that affine policies are in $O\left(\frac{\log(\numu)}{\log\log(\numu)}\right)$ for two-stage problems with non-negative $\cvec, \dvec, \dmat$. 
	Furthermore, our piecewise affine policies are strictly dominated by affine policies for integer budgeted uncertainty.
	
	\begin{proposition}
		\label{pro:aff_leq_pap}
		Consider Problem (\ref{eq:problem_formulation}) with budgeted uncertainty and an integer budget. 
		Let $Z_{PAP}$ be the optimal value found by our piecewise affine policy and $Z_{AFF}$ be the optimal value found by an affine policy. Then
		$$
		Z_{AFF} \leq Z_{PAP}.
		$$
	\end{proposition}
	\noindent
	We present the proof for Proposition \ref{pro:aff_leq_pap} in Appendix \ref{sec:proof_aff_leq_pap}.
	
	\textbf{Norm Ball Uncertainty:}
	In a similar manner as before, we construct new dominating sets for $p$-norm ball uncertainty and tighten the bound in Ben-Tal et al. \cite{BenTal2020} by a factor of $2^{-1+\frac{1}{p}-\frac{1}{p^2}}p^{\frac{1}{p}}(p-1)^{\frac{1}{p^2}-\frac{1}{p}}$ for sufficiently large $\numu$. 
	This factor is always smaller than one and converges to $\frac{1}{2}$ for large $p$. 
	\begin{proposition}[$p$-norm ball]
		\label{pro:pball}
		Consider the $p$-norm ball uncertainty set $\U=\left\{\uvec\in\mathbb{R}_+^{\numu} \middle\vert \norm{\uvec}_p\leq 1\right\}$ with $p > 1$. 
		Then a solution for $Z_{AR}(\hat{\U})$ where $\hat{\U}$ is constructed using Criterion (\ref{eq:rot_inv_validity_criterion}) with
		\begin{align*}
			\mu &= 2^{\frac{1}{p}}\left(2(m-1)+2^p\right)^{-\frac{1}{p^2}}p^{-1}\left(p-1\right)^{\frac{1}{p} + \left(1-\frac{1}{p}\right)^2}, \\
			\rho &= 2^{\frac{1}{p}-1}\left(2(m-1)+2^p\right)^{\frac{1}{p}-\frac{1}{p^2}}p^{-1}\left(p-1\right)^{\left(1-\frac{1}{p}\right)^2}
		\end{align*}
		gives a 
		\begin{align*}
			\beta 
			\leq &\left(2(\numu-1) + 2^p\right)^{\frac{1}{p}-\frac{1}{p^2}}p^{\frac{1}{p}-1}\left(p-1\right)^{\left(1-\frac{1}{p}\right)^2} \\
			\leq& \left((2\numu)^{\frac{1}{p}-\frac{1}{p^2}} + 2^{1-\frac{1}{p}}\right)p^{\frac{1}{p}-1}\left(p-1\right)^{\left(1-\frac{1}{p}\right)^2} 
			= O(\numu^{\frac{1}{p}-\frac{1}{p^2}})
		\end{align*}
		approximation for problem $Z_{AR}(\U)$.
	\end{proposition}
	\noindent
	We present the proof for Proposition \ref{pro:pball} in Appendix \ref{sec:proof_pball}.
	
	\textbf{Ellipsoid Uncertainty:} 
	For the permutation invariant ellipsoid uncertainty set $\left\{\uvec\in\mathbb{R}_+^{\numu} \middle\vert \uvec^\intercal \boldsymbol{\Sigma}\uvec \leq 1\right\}$ with $\numu > 1$, \mbox{$\boldsymbol{\Sigma}:=\boldsymbol{1}+a(\boldsymbol{J}-\boldsymbol{1})$}, $a\in[0,1]$, $\boldsymbol{1}$ being the unity matrix, and $\boldsymbol{J}$ being the matrix of all ones, we construct dominating sets via a case distinction on the size of $a$. 
	While for large $a$ already a scaled simplex gives a good approximation, we construct the dominating set for small $a$ more carefully. 
	By doing so, we improve the previously best known asymptotic bound for the two-stage problem variant of $O(\numu^\frac{2}{5})$ \cite{BenTal2020} to $O(\numu^\frac{1}{3})$. 
	Note, that for $a=0$ our bounds converge to the bounds of hypersphere uncertainty in Proposition \ref{pro:hypersphere} and for $a=1$ towards an exact representation.
	\begin{proposition}[Ellipsoid]
		\label{pro:ellipsoid}
		Consider the ellipsoid uncertainty set $\U=\left\{\uvec\in\mathbb{R}_+^{\numu} \middle\vert \uvec^\intercal \boldsymbol{\Sigma}\uvec \leq 1\right\}$ with $\numu > 1$ and \mbox{$\boldsymbol{\Sigma}:=\boldsymbol{1}+a(\boldsymbol{J}-\boldsymbol{1})$} for $a\in[0,1]$. 
		Here $\boldsymbol{1}$ is the unity matrix and $\boldsymbol{J}$ is the matrix of all ones.
		Then a solution for $Z_{AR}(\hat{\U})$ where $\hat{\U}$ is constructed using Criterion (\ref{eq:rot_inv_validity_criterion}) with
		\begin{align*}
			\mu &= \frac{1}{2\sqrt[4]{(1-a)^3\numu + (1-a)^2a\numu^2}}, &
			\rho &=  \frac{1}{4(1-a)\mu} &
			\text{if } a \leq \numu^{-\frac{2}{3}},\\
			\mu &= 0, &
			\rho &=  \frac{1}{\sqrt{a}} & 
			\text{if } a > \numu^{-\frac{2}{3}},\\
		\end{align*}
		gives a 
		$$\beta = \begin{cases}
			\sqrt{\frac{1}{2}\left(1 + \frac{1}{1-a}\left(a\numu + \sqrt{(1-a)\numu+a\numu^2}\right)\right)} & \text{if } a \leq \numu^{-\frac{2}{3}} \\
			\frac{1}{\sqrt{a}} & \text{if } a > \numu^{-\frac{2}{3}}
		\end{cases} = O(\numu^{\frac{1}{3}})
		$$
		approximation for problem $Z_{AR}(\U)$.
	\end{proposition}
	\noindent
	We present the proof for Proposition \ref{pro:ellipsoid} in Appendix \ref{sec:proof_ellipsoid}.
	
	\textbf{General Uncertainty Sets:}
	After having shown specific bounds for some commonly used permutation invariant uncertainty sets, we now give a general bound that holds for all uncertainty sets that fulfill the assumptions of Problem (\ref{eq:problem_formulation}). 
	We show that any uncertainty set can be dominated within an approximation factor of $\beta=2\sqrt{\numu}+1$, which improves the bound in Ben-Tal et al. \cite{BenTal2020} by a factor of $\frac{1}{2}$. 
	As shown in Ben-Tal et al. \cite{BenTal2020} this approximation bounds is asymptotically tight, when using pure domination techniques. 
	More precisely, for any polynomial number of vertices the budgeted uncertainty set with $k=\sqrt{\numu}$ cannot be dominated with some $\beta$ better than $O(\sqrt{\numu})$.
	\begin{proposition}[General Uncertainty]
		\label{pro:general}
		Consider any uncertainty set $\U\subseteq [0,1]^\numu$ that is convex, full-dimensional with $\unitvec_i\in\U$ for all $i\in\{1,\dots,\numu\}$ and down-monotone.
		Then, there always exists a dominating uncertainty set $\hat{\U}$ of the form in (\ref{eq:dominating_U}) that dominates $\U$ by at most a factor of $\beta = 2\sqrt{\numu}+1$.
	\end{proposition}
	\noindent
	We present the proof for Proposition \ref{pro:general} in Appendix \ref{sec:proof_general}.
	
	\textbf{Stagewise Uncertainty Sets:}
	In general, our policies do not require stagewise uncertainty. 
	However, the existence of such a structure can be utilized in the construction of dominating uncertainty sets, leading to approximation bounds that depend linearly on the stagewise approximation bounds. 
	\begin{proposition}[Stagewise]
		\label{pro:stagewise}
		Let $\U = \U_1\times\dots\times\U_T$ be a stagewise independent uncertainty set and for each $\U_t$, let $\hat{\U}_t$ be a dominating set constructed as in (\ref{eq:dominating_U}). Let $\beta_t$ be the approximation factor for $\hat{\U}_t$, and 
		let $\beta'_t = \min\{\beta'\colon \frac{1}{\beta'}\unitvec\in\U_t\}$ be the constant approximation factor for set $\U_t$.
		Then for any partition $\mathcal{T}_1\cup\mathcal{T}_2 = \{1,\dots,T\}$, $\mathcal{T}_1\cap\mathcal{T}_2 = \emptyset$ of the stages, there exists a dominating set $\hat{\U}$ for $\U$ with approximation factor
		$$   
		\beta  \leq
		\max\left(\sum_{t\in\mathcal{T}_1} \beta_t \,,\; \max_{t\in\mathcal{T}_2}\beta'_t\right).
		$$
	\end{proposition}
	\noindent
	We present the proof for Proposition~\ref{pro:stagewise} in Appendix \ref{sec:proof_stagewise}.
	
	\textbf{Transformed Uncertainty Sets:}
	We note that by the right-hand side $\umat \uvec + \uconst$ of Problem (\ref{eq:problem_formulation}) any positive affine transformation of uncertainty sets $\U$ can be dominated by the same affine transformation of the dominating set $\hat{\U}$. 
	As the approximation bounds do not depend on $\umat$ and $\uconst$, the bounds for the transformed set are the same as for the original set. 
	One well-known uncertainty type covered by these transformations is scaled ellipsoidal uncertainty $\left\{\uvec \middle\vert \sum_{i=1}^\numu w_i \uscal_i^2 \leq 1\right\}$, which was first proposed by Ben-Tal
	and Nemirovski \cite{BenTal1999}. 
	These sets can be constructed via transformations from hypersphere uncertainty sets with a diagonal matrix $\umat$ with $D_{ii} = \frac{1}{\sqrt{w_i}}$. 
	Scaled ellipsoidal uncertainty has been applied to many robust optimization problems, including portfolio optimization \cite{BenTal1999}, supply chain contracting \cite{BenTal2005}, network design \cite{Mudchanatongsuk2008}, and facility location \cite{Baron2011}.
	
	Another widely used class that is partially covered by these positive affine transformations are factor-based uncertainties given by sets ${\U=\left\{\umat \boldsymbol{z} + \uconst \middle\vert \boldsymbol{z}\in\U^z\right\}}$. 
	In these sets, uncertainties affinely depend on a set of factors $\boldsymbol{z}$ that are drawn from a factor uncertainty set $\U^z$. 
	Problems that were solved using such uncertainty sets include, among others, portfolio optimization \cite{Goldfarb2003, Ghaoui2003} and multi-period inventory management \cite{See2010, Ang2012}.
	In contrast to the general factor sets, that have no limitations on $\umat$ and $\uconst$, our approach is restricted to positive factor matrices which allows only for positive correlations between uncertainties. 
	Nevertheless, even this subset of factor-based uncertainties has wide applicational use. 
	As an intuitive example, one could consider component demands where the factors are demands for finished products.
	
	\section{Re-scaling Expensive Vertices}
	\label{sec:rescaling}
	
	On some instances, uncertainties do not have a uniform impact on the objective. 
	While accounting for high values in some of the uncertainty dimensions might drastically influence the objective, accounting for high values in others might barely have an impact. 
	This effect can be crucial in our problem setting, because the creation of dominating sets may overemphasize single uncertainty dimensions by up to a factor of $\beta$ by design. 
	Accordingly, it might thus be beneficial to dominate these uncertainty dimensions more carefully on instances where a few critical uncertainties cause almost all the cost. 
	In this context, we show that it is possible to shrink critical vertices $\vvec_i$ at the cost of slightly shifting all other vertices towards their direction.
	
\begin{figure}[tp]
		\centering
		\includegraphics[width=\textwidth]{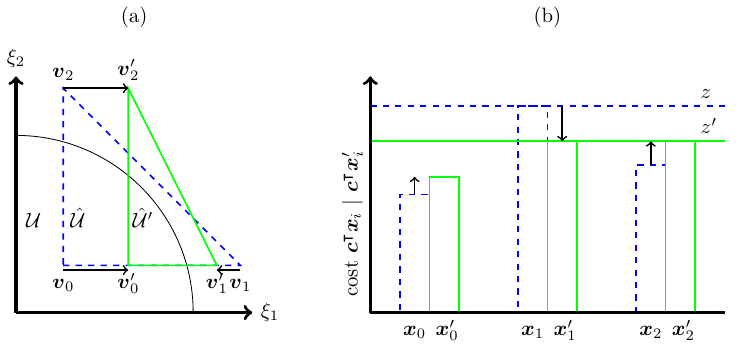}     
		\caption{Re-scaling of the expensive vertex $\vvec_1$ in a dominating set $\hat{\U}$ for uncertainty set $\U$. 
			(a) shows the change of dominating set $\hat{\U}$ (blue, dashed) and it's vertices $\vvec_0, \vvec_1, \vvec_2$ to the new re-scaled dominating set $\hat{\U}'$ (green, solid) with vertices $\vvec'_0, \vvec'_1, \vvec'_2$ for $s_1=0.5, s_2=0$.
			(b) shows the costs for the vertex solutions $\dvec_i$ with maximal cost $z$ (blue, dashed) compared to the costs for the re-scaled vertex solutions $\dvec'_i$ with maximal cost $z'$ (green, solid)
		}
		\label{fig:rescaling_illustration}
	\end{figure}

	We illustrate the re-scaling process following from Lemma \ref{lem:vertex_shift} in Figure \ref{fig:rescaling_illustration}. 
	In the depicted example we shrink the critical vertex $\vvec_1$ and shift the two remaining vertices $\vvec_0, \vvec_2$ towards uncertainty dimension $\uscal_1$. 
	As a consequence, the cost of the vertex solution $\dvec_1$ decreases, while the costs of the solutions $\dvec_0, \dvec_2$ increase. 
	The cost reduction on the critical vertex $\vvec_1$ leads to a reduction of the worst case vertex cost $z$ which by the construction of LP (\ref{eq:lp}) corresponds to an overall improvement of the objective function. 
	Note, that the dominating set $\hat{\U}$ used in the example is not an optimal choice for the hypersphere uncertainty set depicted. 
	However, the effect would barely be visible in two dimensions without this sub-optimal choice.
	
	\begin{lemma}
		\label{lem:vertex_shift}
		Let $\hat{\U}:=\convex(\vvec_0,\dots, \vvec_\numu)$ be a dominating set for $\U$.
		Let $\svec\in [0,1]^\numu$ be a vector of scales.
		Then $\hat{\U}':=\convex(\vvec'_0,\dots, \vvec'_\numu)$ with $v'_{ij} := s_j(1-v_{ij}) + v_{ij}$ is also a dominating set for $\U$.
	\end{lemma}
	\noindent
	We present the proof for Lemma \ref{lem:vertex_shift} in Appendix \ref{sec:proof_vertex_shift}.
	
	The two extreme cases for the modified dominating sets from Lemma \ref{lem:vertex_shift} are given by $\svec=\nullvec$ and $\svec=\unitvec$. 
	While the dominating set does not change for $\svec=\nullvec$, all vertices become the unit vector consisting of ones in every component for $\svec=\unitvec$. 
	Intuitively, increasing $s_i$ leads to shifting the $i$\textsuperscript{th} component towards one. 
	In the same way as we constructed our dominating sets in (\ref{eq:dominating_U}) this shift of the $i$\textsuperscript{th} component towards one increases the value for all $\vvec_j$ with $j\neq i$ and decrease it for $j=i$. 
	Note that this Lemma is not limited to uncertainty sets that are constructed as described in (\ref{eq:dominating_U}), but holds for any dominating set created as a convex combination of vertices.
	
	The transformation used in Lemma \ref{lem:vertex_shift} is linear, which allows us to add $\svec$ as a further decision variable to the second constraint of LP (\ref{eq:lp}) for a given $\hat{\U}$. 
	As $\svec=\nullvec$ gives the original dominating set, any optimal solution found with these additional decision variables is at least as good as a non-modified solution. 
	Thus, all performance bounds shown in Section \ref{sec:bounds} also hold for these re-scaled piecewise affine policies. 
	Note that Proposition \ref{pro:aff_leq_pap} also extends to the re-scaled uncertainty set; thus, all re-scaled piecewise affine policies are strictly dominated by affine policies for integer budgeted uncertainty.
	
	Adding $\svec$ to the LP increases its size, which in practice will often lead to an increase of solution times. 
	To limit the increase of model size, it is possible to add only those $s_i$ where one expects $s_i>0$, as not adding a variable $s_i$ is equivalent to fixing $s_i = 0$. 
	Those $s_i$ with $s_i>0$ correspond to the critical uncertainty dimensions, and an experienced decision maker with sufficient knowledge of the problem might be able to identify them a priori.

	\section{Piecewise Affine Policies via Liftings}
	\label{sec:lifing_combination}
	
	Georghiou et al. \cite{Georghiou2015} propose piecewise affine policies via liftings. 
	In this section, we strengthen these policies by using the insights from the policies constructed in Sections~\ref{sec:framework} and~\ref{sec:bounds}. 
	
	To construct piecewise affine policies via liftings, in the context of \mbox{Assumption \ref{ass:uncertainty_assumption}}, we first choose $r_i-1$ breakpoints
	\begin{equation*}
		0 < z_1^i < z_2^i < \dots < z_{r_i-1}^i < 1,
	\end{equation*}
	for each uncertainty dimension $i\in\{1,2, \dots, \numu\}$.
	For ease of notation let $z_0^i:=0$, $z_{r_i}^i:=1$.
	With these breakpoints, we define the lifting operator $L\colon \mathbb{R}^{\numu} \to \mathbb{R}^{\numu^L}$ ,where $\numu^L:=\sum_{i=1}^\numu r_i$, componentwise by
	\begin{equation*}
		L_{i,j}(\uvec) := \left(\min(\uscal_i - z^i_{j-1}, z^i_j - z^i_{j-1})\right)_+.
	\end{equation*}
	Further, we define the linear retraction operator $R\colon \mathbb{R}^{\numu^L}\to\mathbb{R}^{\numu}$ componentwise by
	\begin{equation*}
		R_i(\uvec^L) := \sum_{j=1}^{r_i} \uscal^L_{i,j}.
	\end{equation*}
	Here $\uscal^L_{i, j}$ are the components of the $\numu^L$ dimensional vector
	$$
	\uvec^L = (\uscal_{1,1}^L, \uscal_{1,2}^L,  \dots,\uscal_{1, r_1}^L, \uscal_{2, 1}^L, \dots, \uscal_{\numu, 1}^L, \dots, \uscal_{\numu, r_\numu}^L)^\intercal \in \mathbb{R}^{\numu^L}.
	$$
	Note that $R\circ L \colon \U \to \U$ is the identity.
	Finally, Georghiou et al. \cite{Georghiou2015} construct a lifted uncertainty set $\U^L$ via
	\begin{equation}
		\label{eq:lifteduncertainty}
		\U^L := \left\{
		\begin{aligned}
			&\uvec^L\in\mathbb{R}^{\numu^L}_+ \colon R(\uvec^L) \in \U, \\
			&\frac{\uscal^L_{i, j+1}}{z^i_{j+1}-z^i_{j}} \leq 
			\frac{\uscal^L_{i, j}}{z^i_{j}-z^i_{j-1}} &\forall i\in\{1,\dots, \numu\}, j\in\{1, \dots, r_i-1\}\\
			&\uscal^L_{i, 1} \leq z^i_1 & \forall i\in\{1,\dots, \numu\}
		\end{aligned}
		\right\}.
	\end{equation}
	Uncertainty set~(\ref{eq:lifteduncertainty}) is an outer approximation of the lifting $L(\U)\subseteq\U^L$ and omits $R(\U^L) = \U$.
	Replacing the uncertainty in Problem (\ref{eq:problem_formulation}) with this lifted uncertainty set yields the lifted adjustable problem
	\begin{equation}
		\label{eq:lifted_formulation}
		\begin{aligned}
			Z^L_{AR}(\U^L) = &\min_{\dvec(\uvec^L)} \max_{\uvec^L\in\U^L} && \cvec^\intercal \dvec(\uvec^L)\\
			&\st && \dmat \dvec(\uvec^L) \geq \umat R(\uvec^L) + \uconst & \forall \uvec^L\in\U^L.
		\end{aligned}
	\end{equation}
	Limiting $\dvec$ to affine policies in the lifted space, yields piecewise affine policies in the original space, which give tighter approximations than affine policies in the original space, i.e., $Z^L_{AFF}(\U^L)\leq Z_{AFF}(\U)$ \cite{Georghiou2015}. 
	However, $\U^L$ is not a tight outer approximation of $L(\U)$, leading to little or no improvements over affine policies on some instances \cite{Georghiou2015, Bertsimas2015b}.
	In fact, we show that in the framework of ARO, the piecewise affine policies induced by lifted uncertainty (\ref{eq:lifteduncertainty}) are equivalent to classical affine policies, in the sense that for any optimal feasible lifted affine policy, there is an affine policy with the same objective value and vice versa.
	\begin{proposition}
		\label{pro:lift_eq_aff}
		Let $Z_{AFF}$ be the optimal objective for affine policies on $Z_{AR}(\U)$ and let $Z_{LIFT}$ be the optimal objective value for lifted affine policies on $Z^L_{AR}(\U^L)$. Then
		$$
		Z_{AFF} = Z_{LIFT}.
		$$
	\end{proposition}
	\noindent
	We present the proof for Proposition \ref{pro:lift_eq_aff} in Appendix \ref{sec:proof_lift_eq_aff}.

	We overcome this shortcoming in the construction of lifted affine policies using our results on dominating uncertainty sets from Section \ref{sec:framework}. 
	Consider the lifting with one break-point $v_{0i}$ per uncertainty dimension.
	Here $v_{0i}$ is the $i$\textsuperscript{th} component of the base vector $\vvec_0$ from Section \ref{sec:framework}.
	Then the lifting operator $L$ becomes
	\begin{equation*}
		L_{ij}(\uvec) = \begin{cases}
			\min(\uscal_i, v_{0i}) &\text{  if } j=1\\
			\left(\umap(\uvec) - \vvec_0\right)_i &\text{  if } j=2.
		\end{cases}
	\end{equation*}
	With this construction of $L$ it is easy to verify that by Condition (\ref{eq:validity_criterion}) each $\uvec^L\in L(\U)$ satisfies
	$
	{\sum_{i=1}^\numu \frac{\uscal^L_{i,2}}{\rho_i} \leq 1}.
	$
	Accordingly, we can tighten the lifted uncertainty set $\U^L$ and get the new lifted uncertainty set 
	\begin{equation}
		\label{eq:tightenedlifteduncertainty}
		\hat{\U}^L := \left\{
		\uvec^L \in \U^L, 
		\sum_{i=1}^\numu \frac{\uscal^L_{i,2}}{\rho_i} \leq 1
		\right\}.
	\end{equation}
	By construction $\hat{\U}^L$ is an outer approximation of $L(\U)$ and we have $L(\U)\subseteq \hat{\U}^L \subseteq \U^L$ and $R(\hat{\U}^L) = \U$.
	Thus, affine policies on the lifted problem with uncertainty set $\hat{\U}^L$ yield valid piecewise affine policies for the original problem. 
	Furthermore, the construction of $\hat{\U}^L$ guarantees that the lifted policies yield tighter approximations than our piecewise affine policies via domination and classical lifted policies with the same breakpoints. 
	Consequently, all approximation bounds for piecewise affine policies via domination also hold for the strengthened piecewise affine policies via lifting. 
	
	\begin{proposition}
		\label{pro:tightenedlifting}
		Let $Z_{LIFT}$ be the optimal objective value found by the lifting policies with breakpoints $\vvec_0$ and lifted uncertainty set $\U^L$ defined in (\ref{eq:lifteduncertainty}), $Z_{TLIFT}$ be the optimal objective value found by the lifting policies with breakpoints $\vvec_0$ and tightened lifted uncertainty set $\hat{\U}^L$ defined in (\ref{eq:tightenedlifteduncertainty}), and $Z_{SPAP}$ be the optimal objective value found by the piecewise affine policies with re-scaling described in Section \ref{sec:rescaling}. Then
		\begin{align*}
			Z_{TLIFT} & \leq Z_{LIFT}, & Z_{TLIFT} & \leq Z_{SPAP}.
		\end{align*}
	\end{proposition}
	\noindent
	We present the proof for Proposition \ref{pro:tightenedlifting} in Appendix \ref{sec:proof_tightenedlifting}.

	\section{Numerical Experiments}
	\label{sec:experiments}
	
	In this section, we present two numerical experiments to compare the performance of our piecewise affine policies for the different constructions of $\hat{\U}$ and our tightened piecewise affine policies via lifting with the performance of other policies from the literature. 
	We compare the performance in terms of both objective value and computational time.
	
	We run both of the following tests with hypersphere uncertainty sets and budgeted uncertainty sets. 
	In the experiments of \mbox{Ben-Tal et al. \cite{BenTal2020}} piecewise affine policies performed particularly well compared to affine adjustable policies for hypersphere uncertainty and relatively bad for budgeted uncertainty. 
	Accordingly, considering these two uncertainty types gives a good impression of the benefits and limitations of piecewise affine policies. 
	Additionally, this experimental design allows us to analyze whether or not the new formulations with the tighter bounds presented in Propositions \ref{pro:hypersphere} \mbox{and \ref{pro:budgeted}} have a significant impact in practice.
	
	In our studies, we compare the following policies: 
	the affine policies described in Ben-Tal et al. \cite{BenTal2004} (AFF), 
	the constant policies resulting from a dominating set $\hat{\U}=\{\unitvec\}$ with only a single point which by down-monotonicity corresponds to a box (BOX), 
	the near-optimal piecewise affine policies with two pieces proposed in Bertsimas and Georghiou \cite{Bertsimas2015b} (BG), 
	our piecewise affine policies constructed as described in Propositions \ref{pro:hypersphere} and \ref{pro:budgeted} (PAP), 
	the piecewise affine policies constructed as described in Propositions 1 \mbox{and 5} in Ben-Tal et al. \cite{BenTal2020} (PAPBT), 
	our piecewise affine policies with the vertex re-scaling heuristic described in Section \ref{sec:rescaling} (SPAP), 
	and our tightened piecewise affine policies via lifting described in Section \ref{sec:lifing_combination} (TLIFT). 
	Note, that piecewise affine policies via lifting from Georghiou et al. \cite{Georghiou2015} are implicitly included in the comparison by Proposition \ref{pro:lift_eq_aff}. 
	In Table \ref{tab:experiment_policies} we give an overview of all policies compared in our experiments.
	
	\begin{table}[htb]
		\caption{Overview of policies compared in the experiments}
		\label{tab:experiment_policies}
		\centering
		\begin{tabular}{lll}
			Policy & Dominating Set or Policy Function & \\
			\hline
			BOX & $\hat{\U}_\text{BOX} = \{\unitvec\}$ & constant policy \\
			AFF & $\dvec^t(\underline{\uvec}^t) = \boldsymbol{P}^t \underline{\uvec}^t + \boldsymbol{q}^t$ & affine adjustable policy \cite{BenTal2004} \\
			BG & $\dvec^t(\underline{\uvec}^t) = \max_i \overline{\boldsymbol{P}}_i^t \underline{\uvec}^t + \overline{\boldsymbol{q}}^t_i - \max_i \underline{\boldsymbol{P}}_i^t \underline{\uvec}^t + \underline{\boldsymbol{q}}^t_i$ & near-optimal pap \cite{Bertsimas2015b}\\
			PAPBT & $\hat{\U}_\text{BT}=\beta\convex(\unitvec_1, \dots, \unitvec_m, \vvec)$ & literature dominating set \cite{BenTal2020}\\
			PAP & $\hat{\U} = \convex(\vvec_0, \vvec_1, \dots, \vvec_{\numu})$ & our dominating set (c.f., Section \ref{sec:construction_dom_set})\\
			SPAP & $\hat{\U} = \convex(\vvec'_0, \vvec'_1, \dots, \vvec'_{\numu})$ & PAP with re-scaling (c.f., Section \ref{sec:rescaling})\\
			TLIFT & AFF on $\hat{\U}^L$ & tightened pap via lifting (c.f., Section \ref{sec:lifing_combination})\\
			\hline
		\end{tabular}%
	\end{table}

	For all studies, we used Gurobi Version 9.5 on a 6 core 3.70 GHz i7 8700K processor using a single core per instance.
	
	\subsection{Gaussian Instances}
	\label{sec:gaussian_instances}
	
	We base our first set of benchmark instances on the experiments of \mbox{Ben-Tal et al. \cite{BenTal2020}} and Housni and Goyal \cite{Housni2021}.
	Accordingly, we generate instances of Problem (\ref{eq:problem_formulation}) by choosing $\numu = \numc = \numd$, $\uconst = \nullvec$, $\umat = \unitytary_{\numu}$ and generate $\cvec, \dmat$ randomly as
	\begin{align*}
		\cvec &= \unitvec + \alpha \boldsymbol{g},\\
		\dmat &= \unitytary + \boldsymbol{G}.
	\end{align*}
	Here, $\unitvec$ is the vector of all ones, $\unitytary$ is the identity matrix, $\boldsymbol{g}$ and $\boldsymbol{G}$ are randomly generated by independent and identically distributed (i.i.d.) standard gaussians, and $\alpha$ is a parameter that increases the asymmetry of the problem. 
	More specifically, $\boldsymbol{G}$ is given by $G_{ij} = \lvert Y_{ij} \rvert/\sqrt{m}$ and $\boldsymbol{g}$ is given by $g_i = \lvert y_i \rvert$, where $Y_{ij}$ and $y_i$ are i.i.d. standard gaussians. 
	Uncertainties $\uvec$ and decision variables $\dvec$ are split into $\lfloor\sqrt{m}\rfloor$ stages where the $i$\textsuperscript{th} decision always belongs to the same stage as the $i$\textsuperscript{th} uncertainty. 
	For the budgeted uncertainty sets we use a budget of $k=\sqrt{\numu}$. 
	We consider values of $\numu = i^2$ for $i\in\{2,\dots,10\}$ and values of $\alpha$ in $\{0,0.1,0.5,1,5\}$. 
	For each pair of $\numu, \alpha$, we consider $30$ instances. 
	To make the results more comparable, we scale all objective values presented by the constant policies results (i.e., $Z_\cdot / Z_{BOX}$) and report averages over all solved instances. 
	For each parameter pair $\numu, \alpha$, we only consider those policies that found solutions on at least $75\%$ of instances within a hard solution time limit of $4$ hours. 
	Additionally, we present all results on a logarithmic scale and artificially lower bound the scale for solution times by $0.01s$ to make the effects on higher solution times more visible.
	
	\begin{figure}[tb]
		\includegraphics[width=\textwidth]{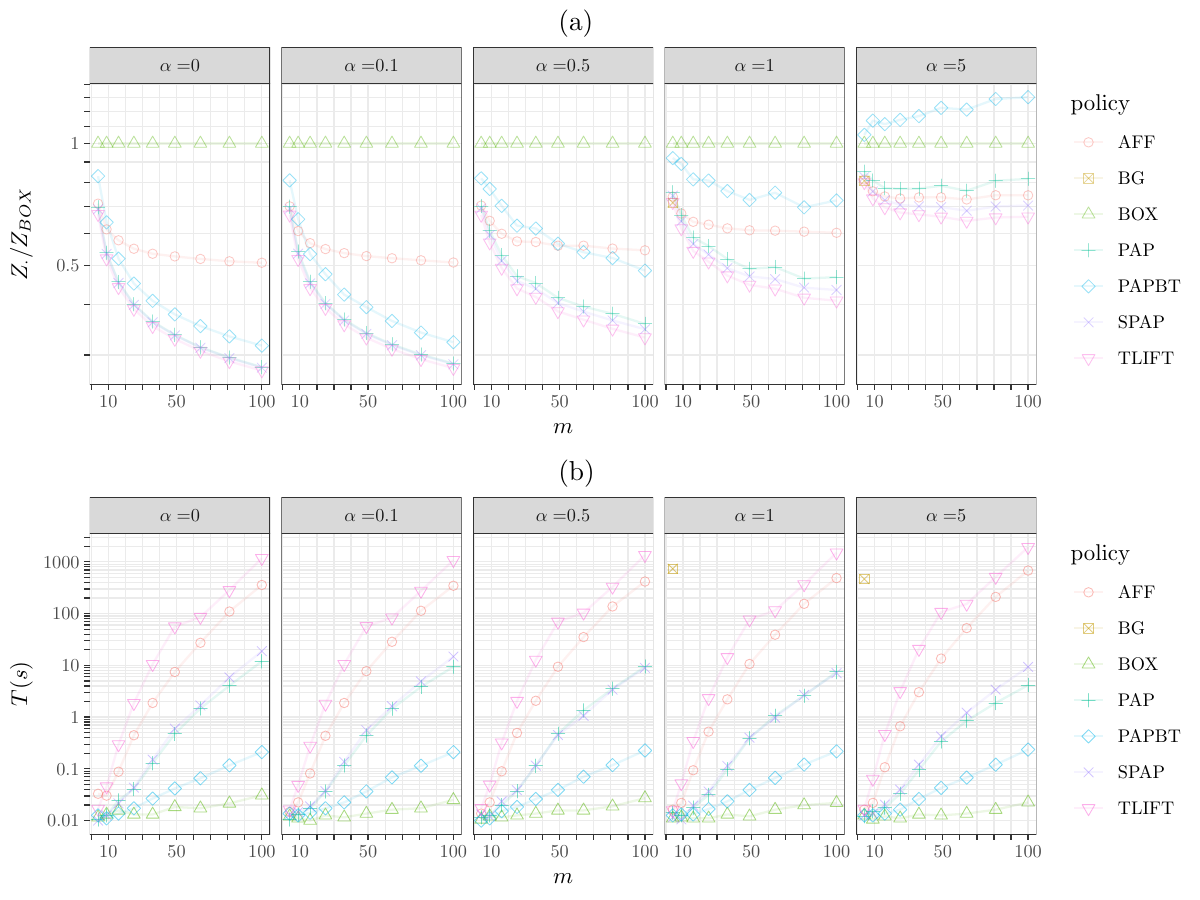}
		\caption{Relative objective values (a) and solution times (b) for different policies on gaussian instances with hypersphere uncertainty} \label{fig:hypersphere_performance}
	\end{figure}
	
	Figure \ref{fig:hypersphere_performance} shows the performances and solution time results on hypersphere uncertainty sets for the different policies. 
	First, note that BG only finds solutions within the time limit for the smallest instances yielding objective values comparable or marginally better to TLIFT. 
	For the other policies we observe that piecewise adjustable policies perform significantly better than affine adjustable policies for small values of $\alpha$. 
	The improvement increases for larger values of $m$ reaching almost a factor of $2$ for $m=100$ on our policies PAP, SPAP, and TLIFT. 
	As expected, the performance of PAP and SPAP for small values of $\alpha$ is almost indistinguishable due to the construction. 
	Additionally, we find that TLIFT only yields marginal improvements over PAP and SPAP for small values of $\alpha$. 
	For larger values of $\alpha$ the improvements of the piecewise affine adjustable policies vanish and TLIFT starts to improve over SPAP. 
	The two policies without re-scaling (PAP and PAPBT) perform even worse than AFF for $\alpha=5$. 
	More severely, PAPBT even performs worse than BOX, which already is a worst-case policy. 
	Only SPAP and TLIFT achieve better results than affine adjustable policies for all values of $\alpha$. 
	
	While solution times for all policies except BOX grow exponentially in the instance size, domination-based piecewise affine adjustable policies are by orders of magnitude faster than classical affine adjustable policies (AFF) and piecewise affine polices via lifting (TPAP). 
	These solution time improvements exceed a factor of $100$ for piecewise affine policies PAP and SPAP and a factor of $1,000$ for PAPBT. 
	Also, solution times of domination-based piecewise affine adjustable policies are barely influenced by values of $\alpha$. 
	This is not the case for AFF and TLIFT, which take longer to solve for increasing $\alpha$. 
	While this effect is not easily visible in Figure \ref{fig:hypersphere_performance} due to the logarithmic scale, the solution time difference for AFF and TLIFT between $\alpha=0$ and $\alpha=5$ reaches up to a factor of two on large instances.
	
\begin{figure}[tb]
		\includegraphics[width=\textwidth]{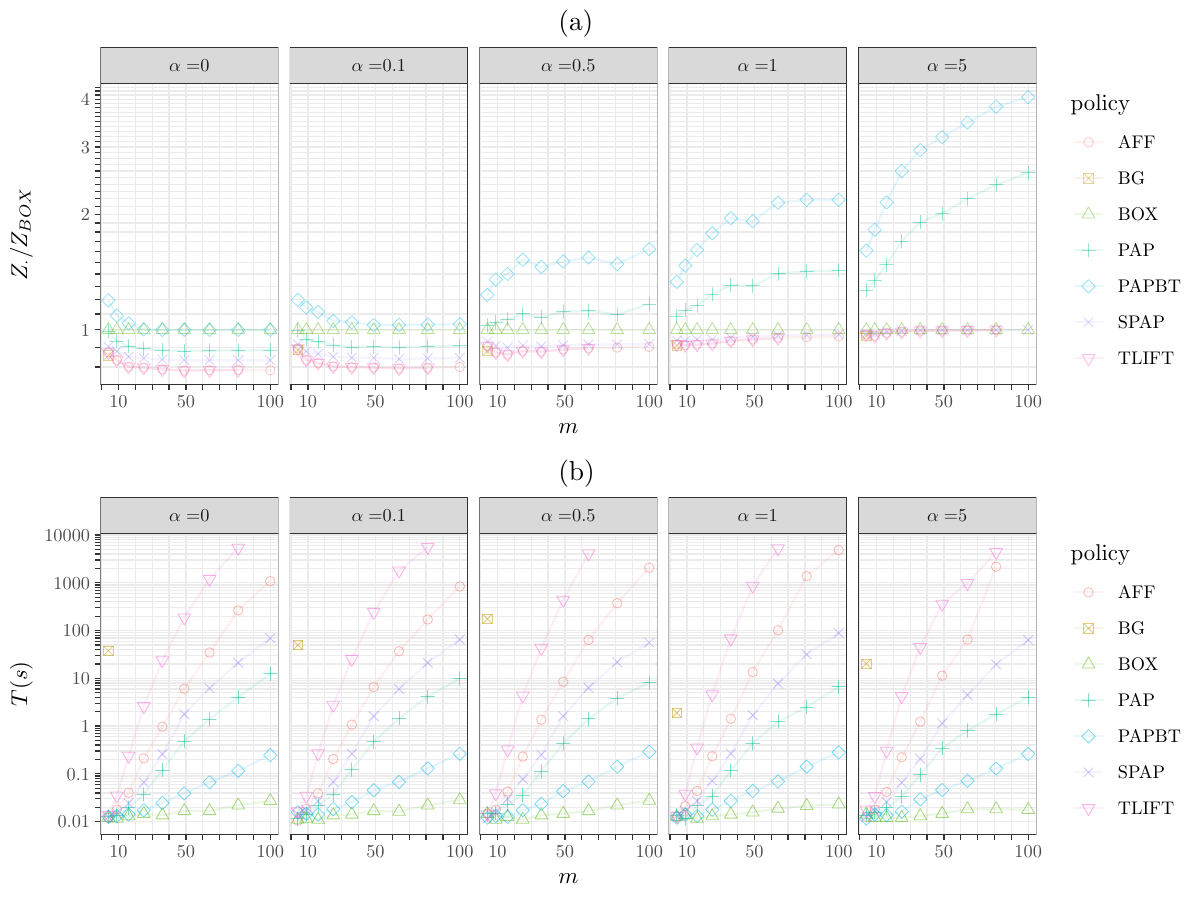}
		\caption{Relative objective values (a) and solution times (b) for different policies on gaussian instances with budgeted uncertainty} \label{fig:budgeted_performance}
	\end{figure}

	Figure \ref{fig:budgeted_performance} shows the performances and solution time results on budgeted uncertainty sets. 
	We observe that for budgeted uncertainty sets domination-based piecewise affine policies perform slightly worse than affine policies throughout all instances. 
	This observation nicely demonstrates that our theoretical and experimental results are aligned, as the worse performance is perfectly explained by \mbox{Proposition \ref{pro:aff_leq_pap}}, which shows that affine policies strictly dominate our piecewise affine policies. 
	Again, we observe that for higher values of $\alpha$, PAP and PAPBT perform even worse than BOX. 
	However, the solution values for SPAP stay within $5\%$ of the affine solution values throughout all instances. 
	We further observe that TLIFT yields the same objective values as AFF throughout all instances. 
	Only BG yields slightly better solutions than AFF on some instances. 
	Solution times behave similarly to solution times on hypersphere uncertainty, confirming that piecewise affine policies are found by orders of magnitude faster on different instances and uncertainty types. 
	Only for BG solution times improve significantly, suggesting that BG is highly dependent on the shape of the uncertainty sets.
	
	Solution time and performance results for $\alpha=0$ align with the results found by Ben-Tal et al. \cite{BenTal2020} for the two-stage problem variant. 
	This demonstrates that our generalized piecewise affine policies do not only extend all theoretical performance bounds, but also achieve comparable numerical results in a multi-stage setting. 
	However, by breaking the symmetry by increasing $\alpha$, we show that pure domination-based piecewise adjustable policies perform poorly on highly asymmetric instances and re-scaling (SPAP) or tightened lifting (TLIFT) constitute good techniques to overcome this shortcoming.
	
	\subsection{Demand Covering Instances}
	\label{sec:demand_instances}
	
	For the second set of test instances, we consider the robust demand covering problem with non-consumed resources and uncertain demands. 
	The problem has various applications, among others in the domains of appointment scheduling, production planning, and dispatching and is especially relevant for the optimization of service levels. 
	Our instances consist of $\numu^l$ locations, $\numu^p$ planning periods, and $\numu^e$ execution periods per planning period. 
	In each execution period $t$, an uncertain demand $\uscal_{lt}$ has to be covered at each location $l$. 
	To do so, the decision maker can buy a fixed number of resources $R$ at a unit cost of $c^R$ in the first stage and then distributes these $R$ resources among the locations at the beginning of each planning period. 
	If a demand cannot be met with the resources assigned to a location, the decision maker will either delay the demand to the next period or redirect it to another location. 
	In either case a fraction $q^d_{tl}\in[0,1]$ or $q^r_{tll'}\in[0,1]$ of the demand is lost. 
	Each unit of lost demand causes costs of $c^D$. 
	Mathematically, the robust demand covering problem with non-consumed resources and uncertain demands is given by the robust LP (\ref{eq:demand_covering_problem}), where parameters and variables are summarized in Table \ref{tab:notation}.
	
	\begin{subequations}
		\label{eq:demand_covering_problem}
		\small
		\begin{align}
			\min \quad& c^R R + c^D\left(\sum_{t\in\mathcal{T},l\in\mathcal{L}} q^d_{t} s^d_{tl} + \sum_{t\in\mathcal{T},l\neq l'\in\mathcal{L}} q^r_{ll'}s^r_{tll'}\right) \label{eq:dem_covering_objective}\span\\
			\text{s.t.}\quad& r_{p(t)l} + s^d_{tl} - (1-q^d_{t-1})s^d_{(t-1) l} \nonumber\\&+ \sum_{l'\in\mathcal{L},l'\neq l}\left( s^r_{tll'} - (1-q^r_{ll'})s^r_{tl'l}\right) \geq \uscal_{tl} & \forall\uvec\in\U, \forall t\in\mathcal{T}, l\in\mathcal{L} \label{eq:dem_covering_constr}\\*
			& \sum_{l\in\mathcal{L}} r_{pl} \leq R & \forall p\in\mathcal{P} \label{eq:dem_covering_resource}\\*
			& R, \boldsymbol{r}, \boldsymbol{s}^r, \boldsymbol{s}^d \geq 0
		\end{align}
	\end{subequations}
	
	\begin{table}[htb]
		\caption{
			Notation for the robust demand covering problem with non-consumed resources and uncertain demands} 
		\label{tab:notation}
		\centering
		\begin{tabular}{ll}
			\hline
			\multicolumn{2}{l}{Sets \& Elements}\\
			$l\in\mathcal{L}$ & locations \\
			$p\in\mathcal{P}$ & planning periods \\
			$t\in\mathcal{T}$ & all ($\numu^p\numu^e$) execution periods\\
			\multicolumn{2}{l}{Parameters}\\
			$c^R$ & cost of buying one unit of resources\\
			$c^D$ & cost of loosing one unit of demand\\
			$q^d_{t}\in[0,1]$ &  fraction of demand lost when delayed in period $t$\\
			$q^r_{ll'}\in[0,1]$ & fraction of demand lost when redirected form location $l$ to $l'$\\
			\multicolumn{2}{l}{Uncertainties}\\
			$p(t)$ & planning period containing execution period $t$\\
			$\uscal_{tl}$ & demand at location $l\in\mathcal{L}$ in execution period $t\in\mathcal{T}$ \\
			\multicolumn{2}{l}{Decision variables}\\
			$R$ & number of resources bought (decided in the very first stage)\\
			$r_{pl}$ & number of resources assigned to location $l\in\mathcal{L}$ in planning period $p\in\mathcal{P}$ \\ 
			& (decided at the beginning of planning period $p$) \\
			$s^r_{tll'}$ & demand redirected from location $l$ to location $l'$(decided in decision stage $t$)\\
			$s^d_{tl}$ & demand delayed to the next period (decided in decision stage $t$)\\
			\hline
		\end{tabular}%
	\end{table}
	
	Here Objective (\ref{eq:dem_covering_objective}) minimizes the sum of resource costs and lost demand costs due to delay and relocation. 
	Constraints (\ref{eq:dem_covering_constr}) ensure that all demands are fulfilled, delayed, or relocated, and Constraints (\ref{eq:dem_covering_resource}) upper bound the allocated resources in each planning period by the total number of available resources $R$. 
	
	In most real-world applications of the demand covering problem, some of the demand will be revealed before the actual demand occurs, e.g. due to already existing contracts, sign-ups, orders, or due to forecasting. 
	To incorporate the increase in knowledge over time, we assume an uncertainty vector of the form $\uvec = \uvec^c+\uvec^p+\uvec^e$. 
	Here, $\uvec^c$ is constant and known before the first stage decision, $\uvec^p$ is revealed before each planning period and $\uvec^e$ accounts for the short-term uncertainties revealed before each execution period. 
	Specifically, we assume that demands are given by
	\begin{equation}
		\label{eq:demand_covering_instance_demand}
		\uscal_{lt} = d_{lt}\left(1+\uscal^p_{lt}+\frac{1}{2}\uscal^e_{lt}\right),
	\end{equation}
	where $d_{lt}$ is the base demand for location $l$ in execution period $t$. 
	Here the uncertainty vector $(\uvec^p,\uvec^e)$ is taken from one base uncertainty set $\U^{B}$ of dimension $\numu =2\numu^l\numu^p\numu^e$.
	Note, that short-term uncertainties have a less severe effect on demands than uncertainties known in advance. 
	The resulting demand uncertainty set
	\begin{equation*}
		\U=\left\{\uvec \,\middle\vert\, \uscal_{lt} = d_{lt}\left(1+\uscal^p_{lt}+\frac{1}{2}\uscal^e_{lt}\right), (\uvec^p,\uvec^e)\in\U^{B} \right\}
	\end{equation*}
	is $m/2$ dimensional, where in our experiments $\U^B$ is either a hypersphere uncertainty set, or a budgeted uncertainty set with budget $\sqrt{m}$.
	
	To construct our instances we draw $\numu^l\in\{2,4,6,8,10\}$ locations at uniform random integer positions in the square $\left[0, 2\lfloor\sqrt{\numu^l}\rfloor + 1\right]^2$. 
	In each of the $\numu^p\in\{1,3,5,7\}$ planning periods, we consider $\numu^e=8$ execution periods corresponding to the hours in a working day. 
	We assume that a fraction $q^d_{tl}=0.1$ of the demand is lost when deferred to a later execution period and consider a doubled loss rate ($q^d_{tl}=0.2$) when demand is deferred to another planning period. 
	Similarly, a fraction of the demand is lost when assigned to another location. 
	We assume this fraction to be correlated to the distance and given by $q^r_{ll'}:=\min\left(1, 0.02\cdot\text{dist}(l,l')\right)$. 
	We draw the base demands $d_{lt}$ uniformly from the normal distribution $\mathcal{N}(10,4)$.
	Finally, we set ${c^R=1}$ and choose $c^D\in\{0.1,0.25,0.5\}$.
	For each combination of $\numu^l, \numu^p$ we consider $45$ instances, where we use each possible value for $c^D\in\{0.1,0.25,0.5\}$ in a third of these instances.
	
	To analyze practical expected objectives, we also report a simulated average objective that the respective policies achieved on $500$ randomly drawn uncertainty realizations $\uvec$, in addition to the robust objective value. %
	We give a detailed description of how uniform uncertainty realizations can efficiently be sampled from the budgeted and hypersphere uncertainty sets in Appendix \ref{sec:drawing_realizations}. 
	We again scale the results by the results achieved by constant policies and use logarithmic scales. 
	For each instance size, we only consider those policies that found solutions on at least $75\%$ of instances within a hard solution time limit of $2$ hours. 
	
	First, we observe that BG did not solve any instance within the time limit 
	which can be attributed to the fact that our demand covering instances are significantly larger than our gaussian instances. 
	
	\begin{figure}[tb]
		\includegraphics[width=\textwidth]{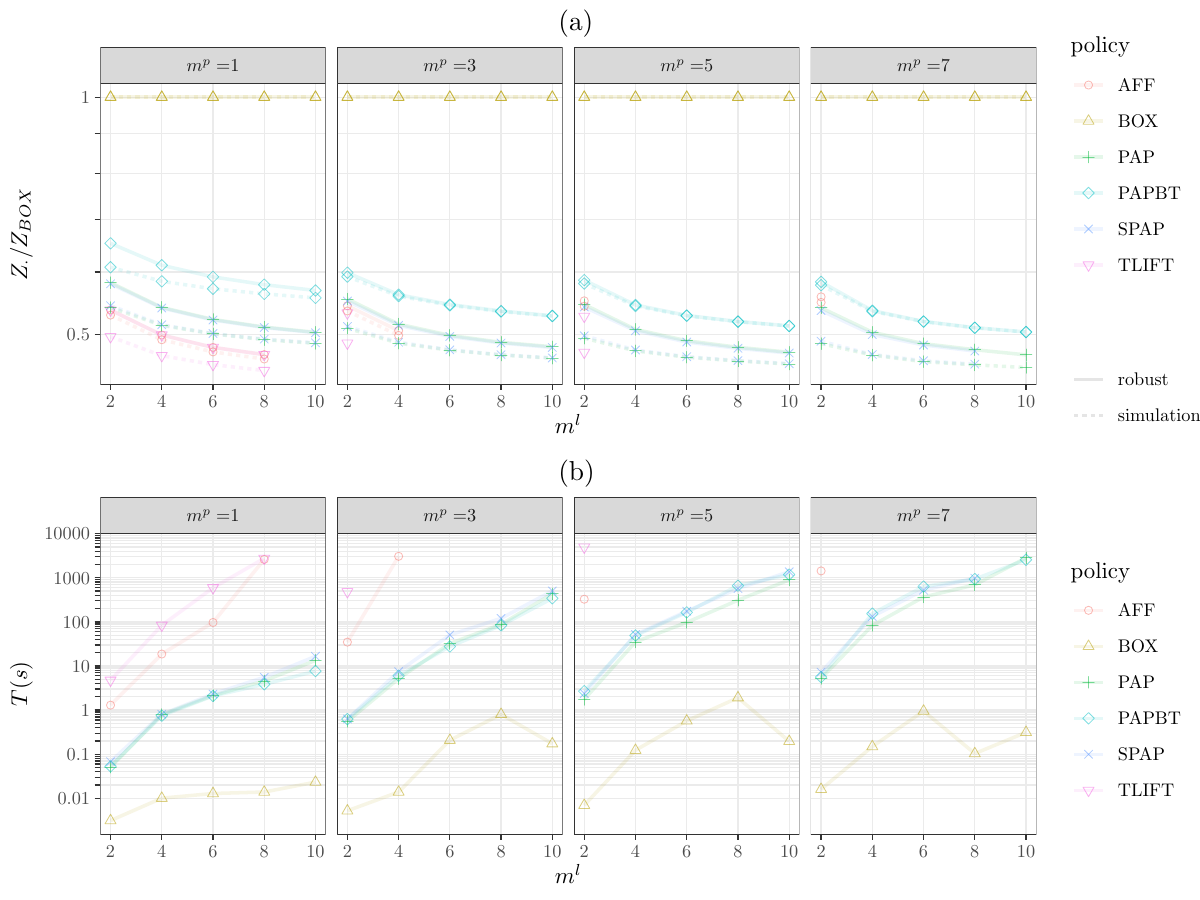}
		\caption{Relative objective values (a) and solution times (b) for different policies on demand covering instances with hypersphere uncertainty} \label{fig:hypersphere_demand_covering}
	\end{figure}
	
	For the remaining policies, Figure \ref{fig:hypersphere_demand_covering} shows the performance and solution time results on demand covering instances with hypersphere uncertainty. 
	Compared to our previous experiment (see \mbox{Section \ref{sec:gaussian_instances}}) we no longer observe the strong objective improvements of piecewise affine policies over affine adjustable policies. 
	Still, our piecewise affine formulations give similar results as affine policies, and we observe small improvements on instances with a larger number of planning periods, with TLIFT yielding strict improvements for $\numu^p\geq 3$.   
	On the simulated realizations, improvements of PAP and SPAP over affine adjustable policies can already be seen for $\numu^p\geq 3$, which might be of interest for a decision maker with practical interest beyond worst-case solutions.
	
	For the solution times, we observe similar improvements to the ones observed on the gaussian instances. 
	Still, all domination-based piecewise affine policies can be found by orders of magnitude faster than affine policies and TLIFT. 
	Interestingly, PAP is solved similarly fast as PAPBT on these instances, while still achieving up to $15\%$ better objective values on all instances. 
	The largest instance that could be solved by affine adjustable policies within two hours consisted of $320$ uncertainty variables and $700$ decision variables, while the largest instance solved by PAPs was more than three times larger with $1,120$ uncertainty variables and $6,230$ decision variables.
	
	\begin{figure}[tb]
		\includegraphics[width=\textwidth]{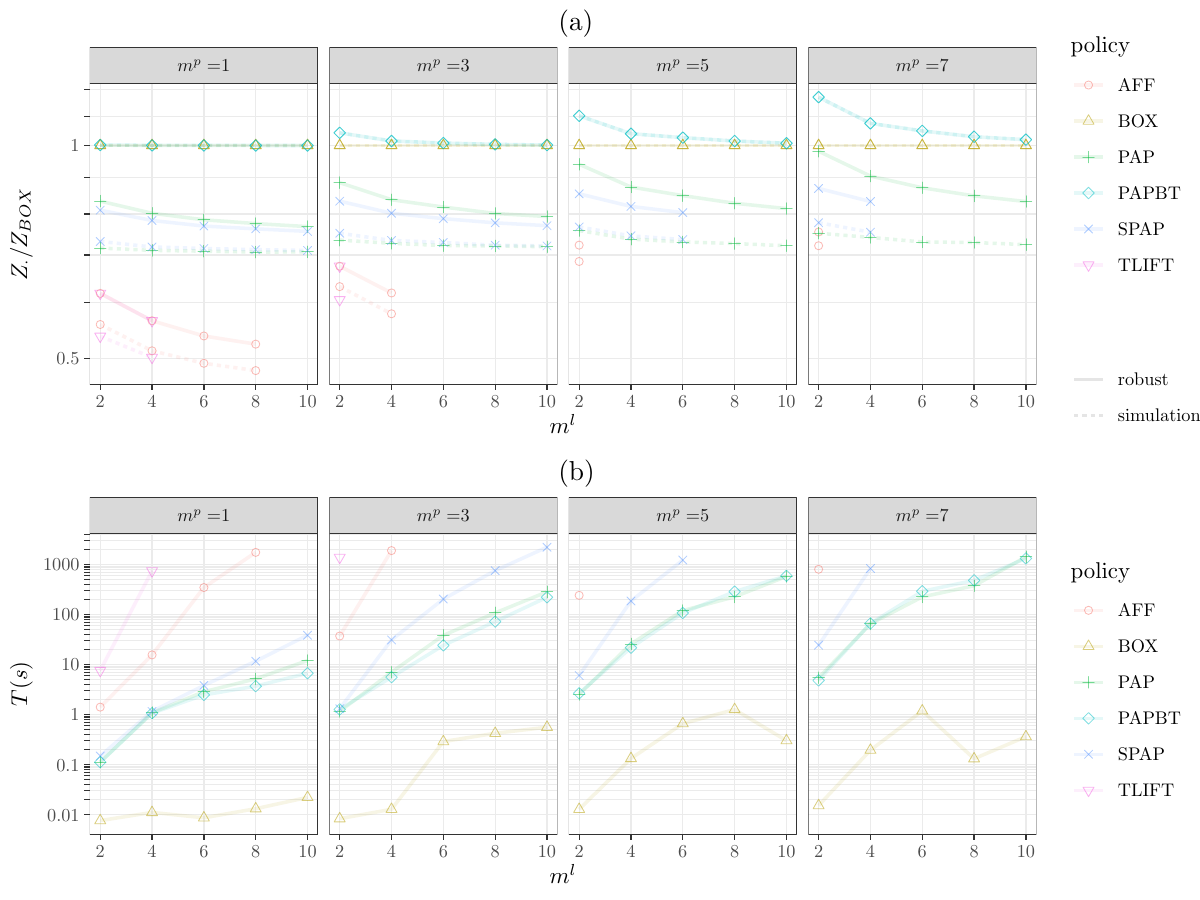}
		\caption{Relative objective values (a) and solution times (b) for different policies on demand covering instances with budgeted uncertainty} \label{fig:budgeted_demand_covering}
	\end{figure}
	
	Figure \ref{fig:budgeted_demand_covering} shows the results on demand covering instances with budgeted uncertainty sets. 
	For budgeted uncertainty sets, domination-based piecewise affine policies perform worse than affine adjustable policies throughout all instances, which again can be explained by Proposition \ref{pro:aff_leq_pap}. 
	Notably, PAPBT even performs worse than constant policies on most instances. 
	In this setting, domination-based piecewise affine policies remain better only from a solution time perspective, as they still solve by orders of magnitude faster than affine policies. 
	Also, TLIFT does not yield any improvements over AFF. 
	
	While in this set of experiments piecewise affine policies do not show the same improvements in the objective over affine adjustable policies, they still perform slightly better with hypersphere uncertainty on larger instances. 
	Also, they still solve by orders of magnitude faster, which makes them an attractive alternative for large-scale optimization in practice.
	
	\subsection{Discussion}
	
	In the experiments presented in Sections \ref{sec:gaussian_instances} and \ref{sec:demand_instances}, some results, e.g. the solution time improvements of piecewise affine policies over affine policies, are consistent throughout all instances. 
	However, other results strongly depend on the set of benchmark instances used. 
	In the following, we give an explanation for the strong solution time improvements, discuss two of the main deviations that we observe between our results on gaussian instances and demand covering instances, and give intuitions of why these differences occur.
	
	\textbf{Size of Robust Counterparts:}
	Throughout all experiments, we see strong solution time improvements of piecewise affine policies over affine policies - including affine policies on lifted uncertainty TLIFT. 
	These can be explained by their respective robust counterparts. 
	The robust counterpart for piecewise affine policies is given by LP (\ref{eq:lp}). 
	For the robust counterparts of affine policies, we refer to Ben-Tal et al. \cite{BenTal2004}. 
	Counterparts of piecewise affine policies have $O(\numd\numu)$ variables compared to $O(\numd\numu + \numc\numu)$ variables for affine policies and both counterparts have $O(\numc\numu)$ constraints. 
	More critically, constraints for piecewise affine policies contain at most $O(\numd)$ variables each and feature a block structure, which solvers use to significantly speed up the solution process. 
	This block structure is only connected by the nonanticipativity constraints. 
	In the robust counterpart of affine policies on the other hand, $O(\numc)$ constraints have up to $O(\numd\numu)$ variables resulting in a denser constraint matrix and the lack of a block structure. 
	Moreover, the robust counterpart for affine policies on hypersphere uncertainty sets is no longer linear. 
	Instead, a quadratic program has to be solved, which tends to be computationally more challenging.
	
	\textbf{Solution Times of PAPBT:}
	In the experiments we see that PAPBT solves by a factor of $10$ to $50$ faster than PAP on gaussian instances, but both find solutions similarly fast on-demand covering instances. 
	This can be explained by the construction of PAPBT and the structure of the instances' constraints. 
	On the gaussian instances, $\dmat$ and $\dvec$ are non-negative and $\umat$ is the unit matrix. 
	In the construction of dominating sets $\hat{\U}$ for PAPBT most of the vertices are chosen to be scaled unit vectors. 
	As a consequence, most Constraints (\ref{eq:lpconstr_constr}) in LP (\ref{eq:lp}) have a zero right-hand side, such that they trivially hold. 
	Consequently, these constraints can be eliminated, which reduces the total number of constraints by a factor of $O(\numu)$.
	On demand covering instances, however, $\dmat$ contains negative entries. 
	Thus, no constraints trivially hold and no constraints can be removed. 
	For a decision maker who is primarily interested in fast policies, this gives a good criterion on when PAPBT can improve solution times and when no such improvements can be expected.
	
	\textbf{Performance Differences Between Gaussian and Demand Covering Instances:}
	We observe that the strong performance improvements of piecewise affine policies over affine policies on gaussian instances with hypersphere uncertainty do not transfer to our demand covering instances. 
	This suggests that the relative performance between piecewise affine policies and affine policies significantly depends on the structure of the problem at hand. 
	An intuitive explanation for this lies in the policies' construction.
	
	Recall that piecewise affine policies derive solutions by finding vertex solutions $\dvec_i$ that can be extended to a full solution. 
	Thereby, $\dvec_0$ focuses on finding a good solution for uncertainty realizations where all uncertainties take equal values, and each $\dvec_i$ focuses on finding a good recourse to uncertainty $\uscal_i$. 
	Thus, good results can be expected when there are 
	(a) synergy effects that can be utilized by $\dvec_0$, and 
	(b) good universal recourse decisions for each uncertainty $\uscal_i$ that do not depend on the realization of other uncertainty dimensions and can be exploited by $\dvec_i$. 
	On the other hand, affine policies directly find solutions on the original uncertainty set. 
	In doing so, they do not depend as strongly on good universal recourse decisions as piecewise affine policies do. 
	However, they also lack the ability to use synergy effects in the way vertex solutions $\dvec_0$ do.
	
	The gaussian instances used fulfill both of the properties that are favorable for piecewise affine policies. 
	By being based on a unity matrix, $\dmat$ has relatively large values along the diagonal, leading to the existence of good universal recourse decisions. 
	Additionally, the relatively small non-negative entries on the non-diagonals lead to synergy effects for uncertainty realizations with many small values. 
	Demand covering instances, however, do not fulfill these properties. 
	The question of how to redirect demand optimally heavily depends on the demand observed at other locations. 
	Also, the only synergistic effects that can be used solely emerge when multiple demands occur at the same location in a single planning period.
	
	On general instances in practice, we would thus not expect to see the same performance improvements that could be observed on our gaussian benchmark instance. 
	Still, piecewise affine policies find solutions by orders of magnitude faster than affine policies and achieve good results throughout all benchmark instances with hypersphere uncertainty. 
	Additionally, Properties (a) and (b) give intuitive criteria on when strong objective improvements over affine policies can be expected.
	
	\section{Conclusion}
	\label{sec:conclusion}
	
	In this work, we presented piecewise affine policies for multi-stage adjustable robust optimization. 
	We construct these policies by carefully approximating uncertainty sets with a dominating polytope, which yields a new problem that we efficiently solve with a linear program. 
	By making use of the problem's structure, we then extend solutions for the new problem with approximated uncertainty to solutions for the original problem. 
	We show strong approximation bounds for our policies that extend many previously best-known bounds for two-stage ARO to its multi-stage counterpart. 
	By doing so, we contribute towards closing the gap between the state of the art for two-stage and multi-stage ARO. 
	To the best of our knowledge, the bounds we give are the first bounds shown for the general multi-stage ARO Problem. 
	Furthermore, our bounds yield constant factor as well as asymptotic improvements over the state-of-the-art bounds for the two-stage problem variant.
	
	In two numerical experiments, we find that our policies find solutions by a factor of 10 to 1,000 faster than affine adjustable policies, while mostly yielding similar or even better results. 
	Especially for hypersphere uncertainty sets our new policies perform well and sometimes even outperform affine adjustable policies up to a factor of two. 
	We observe particularly high improvements on instances that exhibit certain synergistic effects and allow for universal recourse decisions. 
	However, on some instances where few uncertainty dimensions have a high impact on the objective, pure piecewise affine policies perform particularly badly by design, sometimes even worse than constant policies. 
	To mitigate this shortcoming, we present an improvement heuristic that significantly improves the solution quality by re-scaling the critical uncertainty dimension. 
	Furthermore, we construct new tightened piecewise affine policies via lifting that integrate the two frameworks of piecewise affine policies via domination and piecewise affine policies via lifting and combine their approximative power. 
	
	While this work extends most of the best-known approximation results for a relatively general class of ARO problems from the two-stage to the multi-stage setting, it remains an open question whether other strong two-stage ARO results can be generalized to multi-stage ARO in a similar manner. 
	Answering this question remains an interesting area for further research. 
	In this context, binary and uncertain recourse decisions remain particularly relevant challenges. 
	Our analysis in Section~\ref{sec:limitations} has shown that the extension of our policies to encompass these recourse decision types is not straightforward. 
	Nevertheless, exploring the integration of our methodology into the established approaches of piecewise constant policies and k-adaptability, which have proven to be effective in these cases, appears as a promising starting point for future work. 
	Another interesting area for future research is the extension of piecewise affine policies and the concept of domination to adjustable data-driven and distributionally robust optimization. 
	More specifically, we believe that one can obtain tractable data-driven policies by directly fitting the polyhedral uncertainty sets used to construct our policies from data.

	\begin{appendices}
		
		\section{Comparison of Literature Approximation Bounds}
		\label{sec:literature_bounds}
		
		Table \ref{tab:literature_bounds} summarizes existing approximation bounds for some commonly used uncertainty sets in two-stage ARO. 
		We want to point out that in addition to being less restrictive, our bounds for multi-stage ARO presented in Table \ref{tab:optimality_bounds} outperform all these bounds except the ones by Housni and Goyal \cite{Housni2021}, which require the significant restrictive assumption of $\dmat, \dvec, \cvec$ being non-negative. 
		\begin{table}[h]
			\caption{Approximation Bounds for Two-Stage ARO}
			\label{tab:literature_bounds}
			\renewcommand{\theuncertaintyexampleid}{\Roman{uncertaintyexampleid}}
			\newcommand{\uncertaintyexample}{{\stepcounter{uncertaintyexampleid}\theuncertaintyexampleid}}%
			\setcounter{uncertaintyexampleid}{0}
				\centering
				\begin{tabular}{llllll}
					\hline
					No.&  \cite{Bertsimas2011} & \cite{Bertsimas2012} & \cite{Bertsimas2015} & \cite{BenTal2020} & \cite{Housni2021}\\
					\hline
					\hline
					\uncertaintyexample
					& $1+\sqrt{\numu}$ 
					& $-$
					& $\frac{1}{2}(1+\sqrt{\numu})$
					& $\sqrt[4]{\numu}$
					& $-$
					\\
					\uncertaintyexample
					& $1+\frac{\numu}{k}$ 
					& $-$
					& $\frac{k^2+k\numu}{k^2 + \numu}$
					& $\min(k, \frac{\numu}{k})$
					& $\frac{8\log \numu}{\log\log \numu} + 1$
					\\
					\uncertaintyexample
					& $1+\sqrt[p]{\numu}$
					& $-$
					& $\frac{\numu^{\frac{p-1}{p}}+\numu}{\numu^{\frac{p-1}{p}}+\numu^\frac{1}{p}}$
					& $\frac{2}{p}(p-1)^{\frac{p-1}{p}} \numu^\frac{p-1}{p^2}$
					& $-$
					\\
					\uncertaintyexample
					& $-$
					& $-$
					& $-$
					& $\left(\frac{a}{2} + \frac{\sqrt{1-a}}{\sqrt[4]{a\numu^2 + (1-a)\numu}}\right)^{-1}$
					& $-$
					\\
					\uncertaintyexample
					& $\numu$
					& $4\sqrt{\numu}$
					& $-$
					& $4\sqrt{\numu}$
					& $-$
					\\
					\hline
					\hline
					Policy
					& static
					& affine
					& affine
					& piecewise affine
					& affine
					\\
					\hline
				\end{tabular}%
			\setcounter{uncertaintyexampleid}{0}
			
			\begin{minipage}{\textwidth}\footnotesize
				We compare bounds for uncertainty sets of the forms 
				\uncertaintyexample) hypersphere uncertainty $\left\{\uvec\in\mathbb{R}_+^{\numu} \middle\vert \sqnorm{\uvec}\leq 1\right\}$; 
				\uncertaintyexample) budgeted uncertainty $\left\{\uvec\in[0,1]^{\numu} \middle\vert \norm{\uvec}_1\leq k\right\}$ ; 
				\uncertaintyexample) $p$-norm ball uncertainty $\left\{\uvec\in\mathbb{R}_+^{\numu} \middle\vert \norm{\uvec}_p\leq 1\right\}$, with $p\geq1$;
				\uncertaintyexample) ellipsoid uncertainty $\left\{\uvec\in\mathbb{R}_+^{\numu} \middle\vert \uvec^\intercal \boldsymbol{\Sigma}\uvec \leq 1\right\}$, with ${\boldsymbol{\Sigma}:=(1-a)\boldsymbol{1}+a\boldsymbol{J}}$ where $\boldsymbol{1}$ is the unity matrix and $\boldsymbol{J}$ the matrix of all ones;
				\uncertaintyexample) general uncertainty sets $\U\subset\mathbb{R}_+^\numu$.
				Note, that some of the results require more restrictive assumptions on $\dmat, \dvec, \cvec$ than our Problem (\ref{eq:problem_formulation}).
				Most prominently \cite{Housni2021} requires $\dmat, \dvec, \cvec$ to be non-negative.
			\end{minipage}
		\end{table}
		
		\section{Proof of Theorem \ref{thm:bounded_domination}}
		\label{sec:proof_bounded_domination}
		\begin{proof}
			We split the proof in two parts. 
			First, we handle the cases where at least $Z_{AR}(\U)$ or $Z_{AR}(\hat{\U})$ is negative and show that this already implies that both $Z_{AR}(\U)$ and $Z_{AR}(\hat{\U})$ are unbounded. 
			Second, we assume $Z_{AR}(\U), Z_{AR}(\hat{\U})\geq0$ are bounded and prove that in this case, the desired performance bounds hold.
			
			\textbf{Part 1:}
			Let $\tilde{\U}\in\{\U, \hat{\U}\}$ such that $Z_{AR}(\tilde{\U})<0$.
			Then, there exists a solution $\tilde{\dvec}$ such that
			$$
			\max_{\uvec\in\tilde{\U}} \cvec^\intercal \tilde{\dvec}(\uvec)<0.
			$$
			
			Now, assume that $Z_{AR}(\U)$ or $Z_{AR}(\hat{\U})$ is bounded and let $\bar{\U}\in\{\U, \hat{\U}\}$ such that $Z_{AR}(\bar{\U})$ is bounded with an optimal solution $\bar{\dvec}$. 
			Arbitrarily fix one $\tilde{\uvec}\in\tilde{\U}$ and consider the constant vector $\tilde{\dvec}(\tilde{\uvec})$. 
			Then $\bar{\dvec}(\bar{\uvec}) + \tilde{\dvec}(\tilde{\uvec})$ is a feasible solution for $Z_{AR}(\bar{\U})$ as for all $\bar{\uvec}\in\bar{\U}:$
			$$
			\dmat \left(\bar{\dvec}(\bar{\uvec}) + \tilde{\dvec}(\tilde{\uvec})\right) 
			\overset{\text{(a)}}{\geq} \umat\left(\bar{\uvec} + \tilde{\uvec}\right) + 2\uconst 
			\overset{\text{(b)}}{\geq} \umat\bar{\uvec} +  \uconst.
			$$
			Here, (a) holds because $\tilde{\dvec}, \bar{\dvec}$ are feasible solutions and (b) holds as $\tilde{\uvec}, \umat, \uconst$ are non-negative. 
			For the objective we then find
			$$
			\max_{\bar{\uvec}\in\bar{\U}} \cvec^\intercal \left(\bar{\dvec}(\bar{\uvec}) + \tilde{\dvec}(\tilde{\uvec})\right)
			\leq \max_{\bar{\uvec}\in\bar{\U}} \cvec^\intercal \bar{\dvec}(\bar{\uvec}) + \max_{\uvec\in\tilde{\U}} \cvec^\intercal \tilde{\dvec}(\uvec)
			< Z_{AR}(\bar{\U}).
			$$
			This is a contradiction to $\bar{\dvec}$ being a minimal solution. 
			Thus, $Z_{AR}(\bar{\U})$ cannot be bounded. 
			We have thus shown that if any of $Z_{AR}(\U)$, $Z_{AR}(\hat{\U})$ are negative, both have to be unbounded.
			
			\textbf{Part 2:}
			Assume $Z_{AR}(\U), Z_{AR}(\hat{\U}) \geq 0$ are bounded. 
			Let $\hat{\dvec}$ be an optimal solution for $Z_{AR}(\hat{\U})$. 
			Furthermore, let $\umap \colon \U \to \hat{\U}$ be the domination function from \mbox{Definition \ref{def:domination}}. 
			We claim that ${\tilde{\dvec}:=\hat{\dvec}\circ\umap}$ is a feasible solution for $Z_{AR}(\U)$. 
			First, we see that by the definition of $\umap$ the solution $\tilde{\dvec}$ fulfills the nonanticipativity requirements. 
			Specifically, we have
			$$
			\tilde{\dvec}(\uvec) = \hat{\dvec}(\umap(\uvec)) = \left(\hat{\dvec}^1\left(\umap^1(\underline{\uvec}^1)\right), \dots, \hat{\dvec}^T\left(\umap^T(\underline{\uvec}^T)\right)\right)
			$$
			where the decisions in stage $t$ depend only on the uncertainty revealed up to that stage. 
			For the constraints, we find
			$$
			\dmat \tilde{\dvec}(\uvec) = \dmat \hat{\dvec}(\umap(\uvec)) 
			\overset{\text{(a)}}{\geq} \umat\umap(\uvec) + \uconst 
			\overset{\text{(b)}}{\geq} \umat\uvec + \uconst.
			$$
			Here, (a) follows from the feasibility of $\hat{\dvec}$ and (b) follows from the definition of $\umap$ and $\umat$ being non-negative. 
			Thus, $\tilde{\dvec}$ is a well-defined feasible solution for $Z_{AR}(\U)$ and we have
			$$
			Z_{AR}(\U) \leq \max_{\uvec\in\U} \cvec^\intercal \tilde{\dvec}(\uvec) 
			\leq \max_{\hat{\uvec}\in\hat{\U}} \cvec^\intercal\hat{\dvec}(\hat\uvec) = Z_{AR}(\hat{\U}).
			$$
			
			For the other direction let $\dvec^*$ be an optimal solution for $Z_{AR}(\U)$. 
			Then, for all $\hat{\uvec}\in\hat{\U}$ we have $\frac{1}{\beta}\hat{\uvec}\in\U$ by definition of $\beta$. 
			We define $\tilde{\dvec}(\hat{\uvec}):=\beta\dvec^*\left(\frac{1}{\beta} \hat{\uvec}\right)$ and find
			$$
			\dmat \tilde{\dvec}(\hat{\uvec}) = \beta\dmat\dvec^*\left(\frac{1}{\beta}\hat{\uvec}\right)
			\overset{\text{(a)}}{\geq} \beta\left(\frac{1}{\beta}\umat\hat{\uvec} + \uconst\right) 
			\overset{\text{(b)}}{\geq} \umat\hat{\uvec} + \uconst
			$$
			where (a) follows from the feasibility of $\dvec^*$ and (b) from the non-negativity of $\uconst$. 
			Thus, $\tilde{\dvec}$ is a well-defined feasible solution for $Z_{AR}(\hat{\U})$ and we have
			$$
			Z_{AR}(\hat{\U}) \leq \max_{\hat{\uvec}\in\hat{\U}} \cvec^\intercal \tilde{\dvec}(\hat{\uvec}) 
			\leq \beta\max_{\uvec\in\U} \cvec^\intercal\dvec^*(\uvec) = \beta Z_{AR}(\U).
			$$
			Having shown both inequalities this concludes the proof.
		\end{proof}
		
		\section{Proof of Lemma \ref{lem:lp_sol}}
		\label{sec:proof_lp_sol}
		\begin{proof}
			We begin by showing $Z_{AR}(\hat{\U}) \geq Z_{LP}(\hat{\U})$.
			Let $\dvec$ be an optimal solution of $Z_{AR}(\hat{\U})$.
			Then, $(\dvec_i, z)$ is a valid solution for $Z_{LP}(\hat{\U})$, where $\dvec_i, z$ are defined by
			$$
			\begin{aligned}
				&\dvec_i := \dvec(\vvec_i) & \forall i\in\{0,\dots, \numu\},\\
				& z:=\max_{i\in\{0,\dots, \numu\}} c^\intercal \dvec_i.
			\end{aligned}
			$$
			The first set of LP constraints (\ref{eq:lpconstr_max}) holds by definition of $z$.
			The second set of constraints (\ref{eq:lpconstr_constr}) holds as $\dvec$ is a solution of $Z_{AR}(\hat{\U})$ and $\vvec_i\in\hat{\U}$, and the last set of constraints (\ref{eq:lpconstr_nonanticipativity}) holds by nonanticipativity of $\dvec$ and the definition of $\vvec_i$.
			Furthermore, we find
			$$
			Z_{LP}(\hat{\U}) \leq z = \max_{i\in\{0,\dots, \numu\}} c^\intercal \dvec_i \leq \max_{\uvec\in\hat{\U}} c^\intercal \dvec(\uvec) = Z_{AR}(\hat{\U}).
			$$
			
			For the other direction let $(\dvec_i, z)$ be an optimal solution for $Z_{LP}(\hat{\U})$.
			Define $\lambda_i(\uvec)$ for each $\uvec\in\hat{\U}$ by:
			$$
			\lambda_i(\uvec) := 
			\frac{(\uvec - \vvec_0)_i}{\rho_i}.
			$$
			For ease of notation we will in the following drop the explicit dependence on $\uvec$ and use $\lambda_i$.
			We directly find
			$$
			\sum_{i=1}^{\numu} \lambda_i \leq 1
			$$
			for all $\uvec\in\hat{\U}$, as by definition of $\hat{\U}$ each $\uvec$ is a convex combination of $\vvec_0,\dots,\vvec_{\numu}$. 
			Also, note that $\lambda_i$ only depends on the $i$\textsuperscript{th} component of the uncertainty vector $\uvec$. 
			Using this we can define a nonanticipative decision vector $\dvec(\uvec)=\left(\dvec^1(\underline{\uvec}^1), \dots, \dvec^T(\underline{\uvec}^T)\right)$ by
			$$
			\dvec^t(\underline{\uvec}^t) := \sum_{i\in \underline{I}^t} \lambda_i \dvec_i^t + \left(1-\sum_{i\in \underline{I}^t} \lambda_i\right)\dvec_0^t
			$$
			where $\underline{I}^t:=\left\{i\in\{1,\dots,\numu\}\middle\vert \underline{\unitvec_i}^t \neq \nullvec\right\}$ is the index set corresponding to the uncertainties up to stage $t$.
			
			By the nonanticipativity constraints (\ref{eq:lpconstr_nonanticipativity}) we have $\dvec_i^t = \dvec_0^t$ for all $1\leq t \leq T$ and $i\in\{1,\dots,\numu\}\setminus \underline{I}^t$. 
			Thus,
			$$
			\dvec^t(\underline{\uvec}^t) = \sum_{i\in \underline{I}^t} \lambda_i \dvec_i^t + \left(1-\sum_{i\in \underline{I}^t} \lambda_i\right)\dvec_0^t = \sum_{i=1}^{\numu} \lambda_i \dvec_i^t + \left(1- \sum_{i=1}^{\numu} \lambda_i\right)\dvec_0^t,
			$$
			which shows that $\dvec$ is a convex combination of $\dvec_0,\dots, \dvec_{\numu}$ by
			$$
			\dvec(\uvec) = \sum_{i=1}^{\numu} \lambda_i \dvec_i + \left(1- \sum_{i=1}^{\numu} \lambda_i\right)\dvec_0.
			$$
			Using this, we find
			\begin{align*}
				\dmat \dvec(\uvec) &= \sum_{i=1}^{\numu} \lambda_i \dmat\dvec_i + \left(1- \sum_{i=1}^{\numu} \lambda_i\right)\dmat\dvec_0 \\
				&\overset{\text{(a)}}{\geq} \sum_{i=1}^{\numu} \lambda_i \umat\vvec_i + \left(1- \sum_{i=1}^{\numu} \lambda_i\right)\umat\vvec_0 + \uconst
				\overset{\text{(b)}}{=} \umat\uvec + \uconst.
			\end{align*}
			Here, (a) holds because all $x_i$ are valid solutions for (\ref{eq:lp}), and (b) holds as by definition of $\lambda_i$ we have $\uvec = \sum_{i=1}^{\numu} \lambda_i \vvec_i + \left(1- \sum_{i=1}^{\numu} \lambda_i\right)\vvec_0$.
			
			Now we find the missing inequality
			$$
			\begin{aligned}
				Z_{AR}(\hat{\U}) 
				&\overset{\text{(a)}}{\leq} \max_{\uvec\in\hat{\U}} \cvec^\intercal \dvec(\uvec) 
				= \max_{\uvec\in\hat{\U}} \sum_{i=1}^{\numu} \lambda_i \cvec^\intercal\dvec_i + \left(1- \sum_{i=1}^{\numu} \lambda_i\right)\cvec^\intercal\dvec_0\\
				&\overset{\text{(b)}}{\leq} \max_{\uvec\in\hat{\U}} \max_{i\in\{0,\dots, \numu\}} \cvec^\intercal\dvec_i
				= \max_{i\in\{0,\dots, \numu\}} \cvec^\intercal\dvec_i
				= Z_{LP}(\hat{\U}),
			\end{aligned}
			$$
			where we use for (a) that $\dvec$ is a valid solution for $Z_{AR}(\hat{\U})$. 
			For (b), we use that $\dvec$ is a convex combination of $\dvec_0,\dots,\dvec_\numu$ with convex coefficients $(1-\sum_{i=1}^\numu\lambda_i), \lambda_1, \dots,\lambda_\numu$. 
			The convex sum is upper bounded by its largest summand $\cvec^\intercal\dvec_i$. 
			We can thus take the maximum over all these $\cvec^\intercal\dvec_i$, which equals the LP solution value $Z_{LP}(\hat{\U})$.
		\end{proof}
		
		\section{Proof of Lemma \ref{lem:permutation_invariant}}
		\label{sec:proof_lem_permutation_invariant}
		\begin{proof}
			Let $\mathfrak{S}$ be the set of all permutations on $\{1,\dots,\numu\}$ and for any vector $\uvec$ let $\sigma(\uvec)$ be the vector with components permuted according to $\sigma$. 
			Let $\mathfrak{S}_{ji} \subseteq \mathfrak{S}$ be the subset of all permutations mapping component $j$ to component $i$. 
			Let $\hat{\U}$ be a dominating uncertainty set for $\U$ of the form in (\ref{eq:dominating_U}), such that the approximation factor $\beta$ is minimal. 
			Let $\vvec_0$, $\rho_1,\dots,\rho_\numu$ be the parameters defining the vertices of $\hat{\U}$. 
			Define $\mu := \frac{1}{\numu} \unitvec^\intercal \vvec_0, \rho:=\frac{1}{\numu}\sum_{i=1}^{\numu} \rho_i$. 
			Then for each $i\in\{1,\dots,\numu\}$
			\begin{align*}
				\mu \unitvec + \rho \unitvec_i 
				&\overset{\text{(a)}}{=} \frac{1}{\lvert\mathfrak{S}\rvert} \sum_{\sigma\in\mathfrak{S}} \sigma(\vvec_0) + \frac{1}{\numu} \sum_{j=1}^\numu \rho_j \unitvec_i 
				\overset{\text{(b)}}{=} \sum_{j=1}^\numu \left( \frac{1}{\lvert\mathfrak{S}\rvert} \sum_{\sigma\in\mathfrak{S}_{ji}} \sigma(\vvec_0) + \frac{1}{\numu} \rho_j \unitvec_i \right)\\
				&\overset{\text{(c)}}{=} \sum_{j=1}^\numu \left( \frac{1}{\lvert\mathfrak{S}\rvert}\sum_{\sigma\in\mathfrak{S}_{ji}} \sigma(\vvec_0) + \frac{1}{\numu\lvert\mathfrak{S}_{ji}\rvert} \sum_{\sigma\in\mathfrak{S}_{ji}} \rho_j \unitvec_i \right) \\
				&\overset{\text{(d)}}{=} \frac{1}{\lvert\mathfrak{S}\rvert} \sum_{j=1}^\numu  \sum_{\sigma\in\mathfrak{S}_{ji}} \left(\sigma(\vvec_0) + \rho_j \unitvec_i \right)
				\overset{\text{(e)}}{=} \frac{1}{\lvert\mathfrak{S}\rvert} \sum_{j=1}^\numu  \sum_{\sigma\in\mathfrak{S}_{ji}} \sigma\left(\vvec_0 + \rho_j \unitvec_j\right),
			\end{align*}
			where (a) follows from symmetry and the invariance of the average component under permutations, 
			(b) follows from $\mathfrak{S} = \bigcup_{j=1}^\numu \mathfrak{S}_{ji}$ for all $i$,
			(c) follows from ${\sum_{\sigma\in\mathfrak{S}_{ji}} \frac{1}{\lvert\mathfrak{S}_{ji}\rvert}=1}$,
			(d) follows from $\lvert\mathfrak{S}_{ji}\rvert = \frac{1}{\numu} \lvert\mathfrak{S}\rvert$, and
			(e) follows from the linearity of permutations.
			For each $i\in\{1,\dots,\numu\}$ we now conclude
			\begin{align*}
				&\forall j\in\{1,\dots,\numu\}\colon \frac{1}{\beta}\left(\vvec_0 + \rho_j \unitvec_j\right) \overset{\text{(a)}}{\in} \U \\
				\Rightarrow & \forall j\in\{1,\dots,\numu\}\colon\sigma\left(\frac{1}{\beta}\left(\vvec_0 + \rho_j \unitvec_j\right)\right) \overset{\text{(b)}}{\in} \U \\
				\Rightarrow & \frac{1}{\beta} \left(\mu \unitvec + \rho \unitvec_i\right) 
				= \frac{1}{\lvert\mathfrak{S}\rvert} \sum_{j=1}^\numu  \sum_{\sigma\in\mathfrak{S}_{ji}} \sigma\left(\frac{1}{\beta}\left(\vvec_0 + \rho_j \unitvec_j\right)\right) \overset{\text{(c)}}{\in} \U.
			\end{align*}
			Here (a) follows from $\hat{\U}$ being a valid dominating set with approximation factor $\beta$, 
			(b) follows from permutation invariance of $\U$, and 
			(c) from convexity of $\U$. 
			Furthermore, $\frac{1}{\beta} \mu \unitvec \in \U$ as $\U$ is down-monotone. 
			Thus the convex set induced by vertices $\mu \unitvec, {\mu \unitvec + \rho \unitvec_1}, \dots, \mu \unitvec + \rho \unitvec_\numu$ is also a valid dominating set for $\U$ with approximation factor $\beta$.
		\end{proof}
		
		\section{Proof of Proposition \ref{pro:hypersphere}}
		\label{sec:proof_hypersphere}
		\begin{proof}
			Our main idea is to find $\mu, \rho$ such that (\ref{eq:rot_inv_validity_criterion}) is always fulfilled and $\beta$ is minimized. 
			Recall that $\beta$ is given by the minimal value such that $\frac{1}{\beta}\vvec_i\in\U$ for all $\vvec_i$. 
			In the case of the hypersphere uncertainty set, this is given when $\sqnorm{\vvec_i}\leq\beta^2$. 
			As for each $i\in\{1,\dots,\numu\}$, $\vvec_{i}$ is given by $\vvec_i=\vvec_0 + \rho\unitvec_i$ and $\vvec_0 = \mu \unitvec$ we find: 
			\begin{equation}
				\label{eq:pre_beta_hypersphere}
				\sqnorm{\vvec_0}\leq\sqnorm{\vvec_{i}} = (\numu-1)\mu^2 + (\mu+\rho)^2 = \numu\mu^2+2\mu\rho+\rho^2.
			\end{equation}
			This is the term for $\beta^2$ that we want to minimize in the following.
			
			Using Lemma \ref{lem:rotational_invariate_worst_case} the left-hand side in (\ref{eq:rot_inv_validity_criterion}) becomes:
			$$
			\frac{1}{\rho}\max_{j\in\{1,\dots, \numu\}} \sum_{i=1}^j(\gamma(j)-\mu)_+ 
			\overset{\text{(a)}}{=} \frac{1}{\rho}\max_{j\in\{1,\dots, \numu\}} j(\gamma(j)-\mu)_+ 
			\overset{\text{(b)}}{=} \frac{1}{\rho}\max_{j\in\{1,\dots, \numu\}} j(\gamma(j)-\mu).
			$$
			For (a) we use that all summands are the same and for (b) we assume w.l.o.g. that $\max_{j\in\{1,\dots, \numu\}} (\gamma(j)-\mu)\geq 0$, as otherwise (\ref{eq:rot_inv_validity_criterion}) trivially holds.
			
			Using the property $\sqnorm{\uvec} \leq 1$ of the hypersphere uncertainty set, we find $\gamma(j)=\frac{1}{\sqrt{j}}$ for all $j$.
			Substituting $\gamma(j)$ in the above equation we find Property (\ref{eq:rot_inv_validity_criterion}) to become:
			$$
			\max_{j\in\{1,\dots, \numu\}} \sqrt{j}-j\mu \leq \rho.
			$$
			The maximum of the left hand side is taken for $j=\frac{1}{4\mu^2}$ and for any optimal choice of $\mu, \rho$ the inequality will be tight. 
			Thus, we find $\rho = \frac{1}{4\mu}$.
			
			Substituting $\rho$ in Equation (\ref{eq:pre_beta_hypersphere}), we get:
			$$
			\beta = \sqrt{\numu \mu^2 + \frac{1}{2} + \frac{1}{16\mu^2}}.
			$$
			This is minimized by $\mu = \frac{1}{2\sqrt[4]{\numu}}$, which gives $\rho=\frac{\sqrt[4]{\numu}}{2}$ and concludes the proof.
		\end{proof}
		
		\section{Proof of Proposition \ref{pro:performance_aff_ub}}
		\label{sec:proof_performance_aff_ub}
		\begin{proof}
			Consider the following two-stage instance of Problem (\ref{eq:problem_formulation}):
			\begin{equation*}
				\begin{aligned}
					Z (\U) =& \min_{\alpha, \dvec(\uvec)} \max_{\uvec\in\U} && \sqrt{m} \alpha + \unitvec^\intercal\dvec(\uvec)\\
					&\text{s.t.} &&\alpha \unitvec + \dvec \geq \uvec & \forall \uvec\in\U\\
					&&&\alpha, \dvec \geq 0,
				\end{aligned}
			\end{equation*}
			where $\U=\{\uvec \vert \sqnorm{\uvec}\leq1\}\subset[0,1]^\numu$ is the hypersphere uncertainty set, $\alpha\in\mathbb{R}$ is the first stage decision variable and $\dvec\in\mathbb{R}^\numu$ are the second stage decisions that may depend on the uncertainty \mbox{realization $\uvec$}.
			
			We begin by finding a lower bound for affine policies on this problem class. 
			By symmetry of the problem there always exists an optimal affine adjustable solution $\dvec (\uvec)$ of the form
			$$
			x_i = a \uscal_i + b\sum_{j\neq i} \uscal_j + c,
			$$
			for some $a,b,c\in\mathbb{R}$.
			To see this let $\dvec^*, \alpha$ be an optimal affine solution to the above problem and define 
			$$
			\dvec(\uvec):= \frac{1}{\lvert\mathfrak{S}\rvert} \sum_{\sigma \in \mathfrak{S}} \sigma\left(\dvec^*\left(\sigma^{-1}(\uvec)\right)\right),
			$$
			where $\mathfrak{S}$ is the set of all permutations on $\numu$ elements.
			Then
			\begin{align*}
				\alpha \unitvec + \dvec(\uvec) 
				&= \alpha \unitvec + \frac{1}{\lvert\mathfrak{S}\rvert} \sum_{\sigma \in \mathfrak{S}} \sigma\left(\dvec^*\left(\sigma^{-1}(\uvec)\right)\right) 
				\overset{\text{(a)}}{=} \frac{1}{\lvert\mathfrak{S}\rvert} \sum_{\sigma \in \mathfrak{S}} \sigma\left(\alpha \unitvec + \dvec^*\left(\sigma^{-1}(\uvec)\right)\right) \\
				&\overset{\text{(b)}}{\geq} \frac{1}{\lvert\mathfrak{S}\rvert} \sum_{\sigma \in \mathfrak{S}} \sigma\left(\sigma^{-1}(\uvec)\right) = \uvec,
			\end{align*}
			where (a) follows from the permutation invariance of $\alpha\unitvec$ and (b) follows from permutation invariance of $\U$ and feasibility of $\dvec^*$.
			Furthermore,
			\begin{align*}
				\max_{\uvec\in\U} \unitvec^\intercal \dvec(\uvec) 
				&= \frac{1}{\lvert\mathfrak{S}\rvert}\max_{\uvec\in\U} \unitvec^\intercal \sum_{\sigma \in \mathfrak{S}} \sigma\left(\dvec^*\left(\sigma^{-1}(\uvec)\right)\right) 
				\overset{\text{(a)}}{=} \frac{1}{\lvert\mathfrak{S}\rvert}\max_{\uvec\in\U} \sum_{\sigma \in \mathfrak{S}} \unitvec^\intercal \dvec^*\left(\sigma^{-1}(\uvec)\right) \\
				&\overset{\text{(b)}}{\leq} \frac{1}{\lvert\mathfrak{S}\rvert} \sum_{\sigma \in \mathfrak{S}} \max_{\uvec\in\U} \unitvec^\intercal \dvec^*\left(\sigma^{-1}(\uvec)\right)
				\overset{\text{(c)}}{=} \frac{1}{\lvert\mathfrak{S}\rvert} \sum_{\sigma \in \mathfrak{S}} \max_{\uvec\in\U} \unitvec^\intercal \dvec^*\left(\uvec\right) = \max_{\uvec\in\U} \unitvec^\intercal \dvec^*\left(\uvec\right).
			\end{align*}
			Here, (a) follows from permutation invariance of $\unitvec$, 
			(b) follows from the subadditivity of a maximization function, and
			(c) follows from permutation invariance of $\U$.
			Non-negativity of $\dvec$ follows from $\nullvec\in \U$.
			Thus $\dvec, \alpha$ is an optimal feasible solution.
			Additionally for any $\overline{\sigma}\in\mathfrak{S}$,
			\begin{align*}
				\dvec(\overline{\sigma}(\uvec)) 
				&= \frac{1}{\lvert\mathfrak{S}\rvert} \sum_{\sigma \in \mathfrak{S}} \sigma\left(\dvec^*\left(\sigma^{-1}\circ\overline{\sigma}(\uvec)\right)\right)
				= \overline{\sigma} \left(\frac{1}{\lvert\mathfrak{S}\rvert} \sum_{\sigma \in \mathfrak{S}} \overline{\sigma}^{-1}\circ\sigma\left(\dvec^*\left(\sigma^{-1}\circ\overline{\sigma}(\uvec)\right)\right) \right) \\
				&\overset{\text{(a)}}{=} \overline{\sigma} \left(\frac{1}{\lvert\mathfrak{S}\rvert} \sum_{\sigma \in \mathfrak{S}} \sigma\left(\dvec^*\left(\sigma^{-1}(\uvec)\right)\right)\right) 
				= \overline{\sigma} (\dvec(\uvec)),
			\end{align*}
			where (a) follows from $(\overline{\sigma}^{-1}\circ \cdot )$ being an automorphism on $\mathfrak{S}$.
			Let $\sigma_{ij}$ be the permutation switching components $i$ and $j$.
			Then $\dvec(\sigma_{ij}(\nullvec)) = \sigma_{ij}(\dvec(\nullvec))$ implies $x_i(\nullvec) = x_j(\nullvec)$,
			$\dvec(\sigma_{ij}(\unitvec_j)) = \sigma_{ij}(\dvec(\unitvec_j))$ implies $x_i(\unitvec_i) = x_j(\unitvec_j)$ and $x_k(\unitvec_i) = x_k(\unitvec_j)$ for $i\neq k\neq j$, and
			$\dvec(\sigma_{ij}(\unitvec_k)) = \sigma_{ij}(\dvec(\unitvec_k))$ for $i\neq k\neq j$ implies $x_i(\unitvec_k) = x_j(\unitvec_k)$.
			This guarantees the claimed structure.
			
			Let $\uvec\in\U$ and for some $i\leq\numu$ define the vector $\uvec'$ by ${\uvec':=\uvec - \uscal_i\unitvec_i}$.
			Then $\uvec'$ is also in $\U$ and by $x_i(\uvec')\geq 0$ we find that every feasible solution has to fulfill
			\begin{equation}
				\label{eq:bc_geq_0}
				b\sum_{j\neq i} \uscal_j + c\geq0 
			\end{equation}
			for all $\uvec\in\U$.
			Additionally, by $\alpha\unitvec + \dvec(\unitvec_i) \geq \unitvec_i$ every feasible solution also fulfills
			\begin{equation}
				\label{eq:a_lb}
				\alpha + a + c \geq 1 .
			\end{equation}
			Using this we obtain a lower bound for the maximum over $\U$ in the objective function as follows
			\begin{equation*}
				\begin{aligned}
					&\max_{\uvec\in\U} \sqrt{\numu} \alpha + \sum_i\left( a\uscal_i + b\sum_{j\neq i} \uscal_j + c\right)\\
					&\overset{\text{(a)}}{\geq}\max_{\uvec\in\U} \sqrt{\numu} \alpha + a\sum_i \uscal_i = \max_{\uvec\in\U} \sqrt{\numu} \alpha + a\unitvec^\intercal\uvec\\
					&\overset{\text{(b)}}{=} \sqrt{\numu} \alpha + \sqrt{\numu} a
					\overset{\text{(c)}}{\geq}\sqrt{\numu} \alpha + \sqrt{\numu} (1- \alpha - c) = \sqrt{\numu} (1-c).
				\end{aligned}   
			\end{equation*}
			Here (a) follows from Equation (\ref{eq:bc_geq_0}), (b) follows as the maximum is taken for ${\uvec = \frac{1}{\sqrt{\numu}}\unitvec}$, and (c) follows from Equation (\ref{eq:a_lb}).
			On the other hand, using $\nullvec\in\U$ the maximum can also be lower bounded by
			\begin{equation*}
				\sqrt{\numu} \alpha + \sum_i c\\
				= \sqrt{\numu} \alpha + \numu c \geq \numu c.
			\end{equation*}
			Thus we find
			\begin{equation*}
				\begin{aligned}
					Z_{AFF}(\U) &= \min_{\alpha, a, b, c} \max_{\uvec\in\U} \sqrt{\numu} \alpha + \sum_i\left( a\uscal_i + b\sum_{j\neq i} \uscal_j + c\right)\\
					&\geq \min_{c} \max\left(\sqrt{\numu}(1-c), \numu c\right) \overset{\text{(a)}}{=} \sqrt{\numu} - \frac{\numu}{\numu + \sqrt{\numu}} \geq \sqrt{\numu} - 1,
				\end{aligned}
			\end{equation*}
			where (a) holds as the maximum is taken at equality of the two terms, which is given by $c=\frac{\sqrt{\numu}}{\numu + \sqrt{\numu}}$.
			
			Having found a lower bound for affine policies, we will now give an upper bound for optimal policies. 
			Consider the policy given by $\alpha = \frac{1}{\sqrt[4]{\numu}}$ and $\dvec = (\uvec-\alpha\unitvec)_+$. 
			By construction this is feasible and we find
			\begin{equation*}
				\begin{aligned}
					Z_{OPT}(\U) &\leq \max_{\uvec\in\U} \sqrt{\numu}\alpha + \unitvec^\intercal(\uvec-\alpha\unitvec)_+\\
					&\overset{(a)}{=} \max_{0\leq k\leq \numu} \sqrt[4]{m} + k \left(\frac{1}{\sqrt{k}}-\frac{1}{\sqrt[4]{\numu}}\right) \\&= \max_{0\leq k\leq \numu} \sqrt[4]{m} + \sqrt{k} - \frac{k}{\sqrt[4]{\numu}}
					\overset{(b)}{=}  \frac{5}{4}\sqrt[4]{\numu}.
				\end{aligned}
			\end{equation*}
			Here (a) follows from Lemma \ref{lem:rotational_invariate_worst_case} and (b) follows as the maximum over $k$ is taken for $k=\frac{\sqrt{\numu}}{4}$. 
			To finish the proof, we find the optimality ratio, which is given by
			\begin{equation*}
				\frac{Z_{AFF}(\U)}{Z_{OPT}(\U)} \geq \frac{\sqrt{\numu}-1}{5\sqrt[4]{\numu}/4} = \frac{4}{5}\left(\sqrt[4]{\numu} - \frac{1}{\sqrt[4]{\numu}}\right).
			\end{equation*}
		\end{proof}
		
		\section{Proof of Proposition \ref{pro:budgeted}}
		\label{sec:proof_budgeted}
		\begin{proof}
			We proceed analogously as in the proof for Proposition \ref{pro:hypersphere}. 
			Again, the idea is to find $\mu, \rho$ such that (\ref{eq:rot_inv_validity_criterion}) is always fulfilled and $\beta$ is minimized. 
			Using the definition of $\U$ and $\vvec_i$ we find $\frac{1}{\beta}\vvec_i\in\U$ to be fulfilled when
			\begin{equation}
				\label{eq:proof_pro_budget_beta_def}
				\beta = \min\left( \frac{\numu}{k}\mu + \frac{1}{k}\rho, \mu+\rho\right).
			\end{equation}
			Here the first term of the maximization follows from $\frac{1}{\beta}\unitvec^\intercal \vvec_i \leq k$ and the second term from $\frac{1}{\beta}\unitvec_i^\intercal \vvec_i \leq 1$.
			By Lemma \ref{lem:rotational_invariate_worst_case} the left hand side in (\ref{eq:rot_inv_validity_criterion}) becomes:
			$$
			\frac{1}{\rho}\max_{j\in\{1,\dots, \numu\}} j(\gamma(j)-\mu).
			$$
			Using the properties $\norm{\uvec}_1\leq k$ and $\uscal_i\leq 1$ of the budgeted uncertainty set we find $\gamma(j)=\frac{k}{j}$ for all $j\geq k$ and $\gamma(j)=1$ for all $j\leq k$.
			Inserting $\gamma(j)$ into the above, we find Property (\ref{eq:rot_inv_validity_criterion}) to become:
			$$
			\max_{j\in\{1,\dots, \numu\}} j\left(\min\left(1,\frac{k}{j}\right)-\mu\right) \leq \rho.
			$$
			As for $\mu\geq1$ the vector $\vvec_0$ would already dominate every uncertainty realization in $\U$, we can assume $\mu\leq1$.
			With this, the maximum on the left-hand side is taken for $j=k$ and we find the minimal possible $\rho$ to be $\rho = k-k\mu$.
			
			Using this $\rho$ together with $\beta$ from Equation (\ref{eq:proof_pro_budget_beta_def}), we obtain:
			$$
			\beta =\min\left( \left(\frac{\numu}{k}-1\right)\mu + 1, k-(k-1)\mu\right).
			$$
			As the left term in the minimization is linearly increasing in $\mu$ and the right term is linearly decreasing in $\mu$, the minimum is taken at equality.
			This is the case when $\mu=\frac{k(k-1)}{\numu+k(k-2)}$, which gives $\rho= \frac{k(\numu-k)}{\numu+k(k-2)}$.
		\end{proof}
		
		\section{Proof of Proposition \ref{pro:aff_leq_pap}}
		\label{sec:proof_aff_leq_pap}
		To prove Proposition \ref{pro:aff_leq_pap} we first introduce a more general form of piecewise affine policies.
		Recall that affine policies are given by
		${\dvec^t(\underline{\uvec}^t) = \boldsymbol{P}^t \underline{\uvec}^t + \boldsymbol{q}^t}$,
		with decision variables $\boldsymbol{P}^t$ and $\boldsymbol{q}^t$.
		By introducing additional decision variables $\uvecbase$ of the same dimension as the uncertainty set, we extend affine policies to \textit{fully piecewise affine policies} (FPAP)
		$$
		{\dvec^t(\underline{\uvec}^t) = \boldsymbol{P}^t \left(\underline{\uvec}^t-\underline{\uvecbase}^t\right)_+ + \boldsymbol{q}^t},
		$$
		where we use the same notation for $\uvecbase$, as for $\uvec$.
		By extending the decision matrices $\boldsymbol{P}^t$ with zeroes and concatenating them to a single decision matrix
		$$
		\boldsymbol{P} = \left(
		\begin{matrix}
			\boldsymbol{P}^1 & \boldsymbol{0} & \boldsymbol{0} & \cdots & \boldsymbol{0}\\
			\span \boldsymbol{P}^2 & \boldsymbol{0} & \cdots & \boldsymbol{0}\\
			\span\span\vdots\span\span\\
			\span\span\boldsymbol{P}^T\span\span
		\end{matrix}\right),
		$$
		we get the compact expressions 
		$
		{\dvec(\uvec) = \boldsymbol{P} \uvec + \boldsymbol{q}}
		$
		and
		$
		{\dvec(\uvec) = \boldsymbol{P} \left(\uvec-\uvecbase\right)_+ + \boldsymbol{q}}
		$
		for affine and fully piecewise affine policies, respectively.
		
		By our initial assumption $\U\subseteq[0,1]^\numu$, we can w.l.o.g. assume $\nullvec\leq\uvecbase\leq\unitvec$.
		First, consider a policy $\pmat, \uvecbase, \qvec$ with $\uscalbase_i<0$.
		Then replacing $\uscalbase_i$ by $0$ and $\qvec$ by $\qvec':= \qvec - \uscalbase_i \pmat_{:i}$, where $\pmat_{:i}$ is the $i$\textsuperscript{th} column of $\pmat$, yields an identical policy on the uncertainty set, as $\left(\uscal_i-\uscalbase_i\right)_+ = \uscal_i-\uscalbase_i$ for all $\uvec\in\U$.
		Similarly, if $\uscalbase_i>1$, we can replace it by $1$ without changing the policy, as $\left(\uscal_i-\uscalbase_i\right)_+ = \left(\uscal_i-1\right)_+ = 0$ for all $\uvec\in\U$.
		
		We can show that the fully piecewise affine adjustable policies defined above are in fact a generalization of the piecewise affine policies we introduce in this work.
		
		\begin{lemma}
			\label{lem:fpap_leq_pap}
			Consider Problem (\ref{eq:problem_formulation}).
			Let $Z_{PAP}$ be the optimal value found by our piecewise affine policy and let $Z_{FPAP}$ be the optimal value found by a fully piecewise affine policy where the right hand side uncertainty is replaced by $\umat\left((\uvec-\uvecbase)_++\uvecbase\right) + \uconst$.
			Then
			$$
			Z_{FPAP} \leq Z_{PAP}.
			$$
		\end{lemma}
		
		\begin{proof}
			We prove the inequality by showing that every feasible piecewise affine policy is a feasible fully piecewise affine policy.
			Let $\vvec_0,\dots,\vvec_\numu$ be the chosen vertices of the piecewise affine policy and let $\dvec_0,\dots,\dvec_\numu$ be the vertex solutions received from solving LP (\ref{eq:lp}).
			Then the piecewise affine policy is given by Equation (\ref{eq:piecewise_policy}) and we find
			\begin{align*}
				\dvec(\uvec) &= \sum_{i=1}^{\numu} \lambda_i(\uvec) \dvec_i + \left(1-\sum_{i=1}^{\numu} \lambda_i(\uvec)\right) \dvec_0 \\
				&= \sum_{i=1}^{\numu} \lambda_i(\uvec) \left(\dvec_i-\dvec_0\right) + \dvec_0 \\
				&\overset{\text{(a)}}{=} \sum_{i=1}^{\numu} \frac{\left( \left( \uvec - \vvec_0 \right)_+\right)_i}{\rho_i} \left(\dvec_i-\dvec_0\right) + \dvec_0\\
				&\overset{\text{(b)}}{=} \left(\frac{1}{\rho_1}\left(\dvec_1-\dvec_0\right),\dots,\frac{1}{\rho_\numu}\left(\dvec_\numu-\dvec_0\right)\right) \left( \uvec - \vvec_0 \right)_+ + \dvec_0.
			\end{align*}
			Here (a) follows from the definition of $\lambda_i$ in Section \ref{sec:construction_dom_fct} and (b) follows by replacing the sum with a matrix-vector multiplication.
			By the nonanticipativity Constraint (\ref{eq:lpconstr_nonanticipativity}) the vector $\dvec_i-\dvec_0$ is zero for all uncertainty dimensions associated to later stages.
			Thus setting
			\begin{align*}
				\boldsymbol{P} &:= \left(\frac{1}{\rho_1}\left(\dvec_1-\dvec_0\right),\dots,\frac{1}{\rho_\numu}\left(\dvec_\numu-\dvec_0\right)\right) \\
				\uvecbase &:= \vvec_0 \\
				\boldsymbol{q} &:= \dvec_0
			\end{align*}
			yields a valid fully piecewise affine policy that is identical to the initial piecewise affine policy.
			Even with the modified right hand side, as demanded in the Lemma, the constructed fully piecewise affine policy fulfills all constraints, as by construction of the policy, and definition of the dominating function $\umap$:
			$$
			\dmat \dvec (\uvec) \geq \umat\umap(\uvec)+\uconst = \umat\left((\uvec-\uvecbase)_++\uvecbase\right) + \uconst.
			$$
		\end{proof}
		
		Even though fully piecewise affine policies are slightly more flexible than affine policies in general, both are equivalent on Problem (\ref{eq:problem_formulation}) under some additional assumptions.
		
		\begin{lemma}
			\label{lem:aff_eq_fpap}
			Consider Problem (\ref{eq:problem_formulation}) with budgeted uncertainty and integer budget.
			Let $Z_{FPAP}$ be the optimal value found by a fully piecewise affine policy with modified right hand side uncertainty as in Lemma \ref{lem:fpap_leq_pap} and let $Z_{AFF}$ be the optimal value found by an affine policy.
			Then
			$$
			Z_{AFF} = Z_{FPAP}.
			$$
		\end{lemma}
		
		For the proof of Lemma \ref{lem:aff_eq_fpap} we first need the following helpful result that allows us to construct a robust counterpart of Problem (\ref{eq:problem_formulation}) with fully piecewise affine adjustable policies in a similar manner as robust counterparts for affine policies are constructed.
		
		\begin{lemma}
			\label{lem:piecewise_affine_dual}
			Let $\cvec\in\mathbb{R}^\numu$, $\uvecbase\in[0,1]^\numu$ and $k\in\mathbb{N}_+$. Then the following two problems yield the same objective.
			\begin{equation}
				\label{eq:dual_lem:lhs_problem}
				\begin{aligned}
					\max_{\uvec}\quad & \cvec^\intercal\left(\uvec-\uvecbase\right)_+\\
					\st\quad& \unitvec^\intercal\uvec\leq k \\
					& \nullvec\leq\uvec\leq\unitvec 
				\end{aligned}
			\end{equation}
			and
			\begin{equation}
				\label{eq:dual_lem:rhs_problem}
				\begin{aligned}
					\min_{\alpha,\betavec}\quad & \unitvec^\intercal\betavec + k\alpha\\
					\st\quad& \betavec + \alpha\unitvec \geq \diag(\unitvec-\uvecbase)\cvec\\
					& \alpha\geq0, \quad\betavec\geq\nullvec.
				\end{aligned}
			\end{equation}
		\end{lemma}
		
		\begin{proof}
			Let $I^k$ be the indices of the largest $k$ values $c_i(1-\uscalbase_i)$ and let $\unitvec_{I^k}$ be the indicator vector of $I^k$.
			Then for Problem (\ref{eq:dual_lem:lhs_problem}) we find the following set of inequalities.
			\begin{align}
				&\max\left\{\cvec^\intercal\left(\uvec-\uvecbase\right)_+
				\middle\vert \unitvec^\intercal\uvec\leq k
				, \nullvec\leq\uvec\leq\unitvec\right\} \label{eq:dual_proof:lhs_probl1}\\
				\overset{\text{(a)}}{\leq}&\max\left\{\cvec^\intercal\diag\left(\unitvec-\uvecbase\right)\uvec
				\middle\vert \unitvec^\intercal\uvec\leq k
				, \nullvec\leq\uvec\leq\unitvec\right\} \label{eq:dual_proof:lin_probl}\\
				\overset{\text{(b)}}{\leq}& \cvec^\intercal\diag\left(\unitvec-\uvecbase\right)\unitvec_{I^k}\nonumber\\
				\overset{\text{(c)}}{\leq}& \cvec^\intercal\left(\unitvec_{I^k}-\uvecbase\right)_+\nonumber\\
				\overset{\text{(d)}}{\leq}& \max\left\{\cvec^\intercal\left(\uvec-\uvecbase\right)_+
				\middle\vert \unitvec^\intercal\uvec\leq k
				, \nullvec\leq\uvec\leq\unitvec\right\} \label{eq:dual_proof:lhs_probl2}
			\end{align}
			Here (a) follows from $\left(\uvec-\uvecbase\right)_+\leq \diag\left(\unitvec-\uvecbase\right)\uvec$.
			Inequality (b) follows as $e_{I^k}$ is an optimal choice for Problem (\ref{eq:dual_proof:lin_probl}).
			For (c), observe that $\diag\left(\unitvec-\uvecbase\right)\unitvec_{I^k}=\left(\unitvec_{I^k}-\uvecbase\right)_+$.
			Finally, (d) follows, as $\unitvec_{I^k}$ is a feasible choice for Problem (\ref{eq:dual_proof:lhs_probl2}).
			As (\ref{eq:dual_proof:lhs_probl1}) and (\ref{eq:dual_proof:lhs_probl2}) are the same problem all inequalities are in fact equalities.
			
			To conclude the proof note that Problem (\ref{eq:dual_proof:lin_probl}) is a linear optimization problem, where taking the dual yields Problem (\ref{eq:dual_lem:rhs_problem}).
		\end{proof}
		
		Having this result, we now show Lemma \ref{lem:aff_eq_fpap}.
		
		\begin{proof}[Proof of Lemma \ref{lem:aff_eq_fpap}]
			To prove this result we show that the robust counterparts of the two policies on Problem (\ref{eq:problem_formulation}) with integer budget uncertainty are equivalent.
			Let $k\in\mathbb{N}_+$ be the budget.
			
			Then Problem (\ref{eq:problem_formulation}) with affine policies becomes
			\begin{equation*}
				\begin{aligned}
					Z_{AFF}(\U) = &\min_{\pmat,\qvec} \quad\max_{\uvec}
					\left(\cvec^\intercal \pmat\uvec 
					\middle\vert
					\nullvec\leq\uvec\leq\unitvec, \unitvec^\intercal\uvec\leq k\right) 
					+ \cvec^\intercal\qvec\\
					\text{s.t.}\quad& \max_{\uvec}\left(\left(\umat_i-\dmat_i\pmat\right)\uvec\middle\vert\nullvec\leq\uvec\leq\unitvec, \unitvec^\intercal\uvec\leq k\right)\leq \dmat_i\qvec-d_i & \forall i\in\{1,\dots,\numc\},
				\end{aligned}
			\end{equation*}
			where $\dmat_i$ and $\umat_i$ are the $i$\textsuperscript{th} rows of $\dmat$ and $\umat$, respectively, and $d_i$ is the $i$\textsuperscript{th} entry of $\uconst$.
			Dualizing the suproblems we find this to be equivalent to
			\begin{equation}
				\label{eq:roc_aff}
				\begin{aligned}
					Z_{AFF}(\U) = &\min_{\pmat,\qvec,\alpha^0,\dots,\alpha^\numc,\betavec^0,\dots,\betavec^\numc} \unitvec^\intercal\betavec^0+k\alpha^0
					+ \cvec^\intercal\qvec\\
					\text{s.t.}\quad& \betavec^0+\alpha^0\unitvec \geq \left(\cvec^\intercal\pmat\right)^\intercal\\
					& \unitvec^\intercal\betavec^i+k\alpha^i\leq\dmat_i\qvec - d_i & \forall i\in\{1,\dots,\numc\}\\
					& \betavec^i+\alpha^i\unitvec \geq \left(\umat_i-\dmat_i\pmat\right)^\intercal & \forall i\in\{1,\dots,\numc\}\\
					&\alpha^i\geq 0, \quad \betavec^i\geq\nullvec & \forall i\in\{1,\dots,\numc\}.
				\end{aligned}
			\end{equation}
			
			Similarly Problem (\ref{eq:problem_formulation}) with fully piecewise affine policies and modified right hand side becomes
			\begin{equation*}
				\begin{aligned}
					Z_{FPAP}(\U) = &\min_{\pmat,\qvec,\uvecbase} \quad\max_{\uvec}
					\left(\cvec^\intercal \pmat(\uvec-\uvecbase)_+ 
					\middle\vert
					\nullvec\leq\uvec\leq\unitvec, \unitvec^\intercal\uvec\leq k\right) 
					+ \cvec^\intercal\qvec\span\\
					\st\quad& 
					\max_{\uvec}\left(\left(\umat_i-\dmat_i\pmat\right)(\uvec-\uvecbase)_+\middle\vert\nullvec\leq\uvec\leq\unitvec, \unitvec^\intercal\uvec\leq k\right) 
					\\&\quad\quad\quad\quad\quad\quad\leq \dmat_i\qvec-\umat_i\uvecbase-d_i
					& \quad \forall i\in\{1,\dots,\numc\}\\
					&\nullvec\leq\uvecbase\leq\unitvec.
				\end{aligned}
			\end{equation*}
			Using Lemma \ref{lem:piecewise_affine_dual} on the subproblems we find this to be equivalent to 
			\begin{subequations}
				\label{eq:roc_spap}
				\begin{align}
					Z_{FPAP}(\U) = &\min_{\pmat,\qvec,\uvecbase,\alpha^0,\dots,\alpha^\numc,\betavec^0,\dots,\betavec^\numc} \unitvec^\intercal\betavec^0+k\alpha^0
					+ \cvec^\intercal\qvec\\
					\text{s.t.}\quad& \betavec^0+\alpha^0\unitvec \geq \left( \cvec^\intercal\pmat\diag(\unitvec-\uvecbase)\right)^\intercal\label{eq:proof_lem_spap_eq_aff:constr_first_dual_const}\\
					& \unitvec^\intercal\betavec^i+k\alpha^i\leq\dmat_i\qvec- \umat_i\uvecbase-d_i& \forall i\in\{1,\dots,\numc\}\label{eq:proof_lem_spap_eq_aff:constr_dual_obj}\\
					& \betavec^i+\alpha^i\unitvec \geq \left(\left(\umat_i - \dmat_i\pmat\right)\diag(\unitvec-\uvecbase)\right)^\intercal & \forall i\in\{1,\dots,\numc\} \label{eq:proof_lem_spap_eq_aff:constr_dual_constr}\\
					&\alpha^i\geq 0, \quad \betavec^i\geq\nullvec & \forall i\in\{1,\dots,\numc\}\\
					&\nullvec\leq\uvecbase\leq\unitvec.
				\end{align}
			\end{subequations}
			We now show that for each solution of Problem (\ref{eq:roc_spap}) there is an equivalent solution with $\uvecbase = \nullvec$.
			Let $\pmat,\qvec,\uvecbase,\alpha^0,\dots,\alpha^\numc,\betavec^0,\dots,\betavec^\numc$ be a feasible solution to Problem (\ref{eq:roc_spap}).
			We replace $\pmat$ by \mbox{$\pmat':=\pmat\diag(\unitvec-\uvecbase)$}, $\betavec^i$ by $\betavec'^i:=\betavec^i+\diag(\uvecbase)\umat_i^\intercal$ for $i>0$ and $\uvecbase$ by $\nullvec$ and leave all other decision variables untouched.
			As none of the variables in the objective function was changed, it is sufficient to show that the new solution is again feasible.
			Feasibility of Constraint (\ref{eq:proof_lem_spap_eq_aff:constr_first_dual_const}) directly follows from $\pmat'\diag(\unitvec) = \pmat\diag(\unitvec-\uvecbase)$.
			For Constraint (\ref{eq:proof_lem_spap_eq_aff:constr_dual_obj}), we find
			\begin{align*}
				&\unitvec^\intercal\betavec'^i+k\alpha^i \overset{\text{(a)}}{=} \unitvec^\intercal\diag(\uvecbase)\umat_i^\intercal + \unitvec^\intercal\betavec^i+k\alpha^i \\
				\overset{\text{(b)}}{\leq}& \unitvec^\intercal\diag(\uvecbase)\umat_i^\intercal  + \dmat_i \qvec- \umat_i\uvecbase-d_i \overset{\text{(c)}}{=} \dmat_i \qvec-d_i,
			\end{align*}
			where (a) follows from the definition of $\betavec'^i$, (b) from the feasibility of the original solution and (c) from $\unitvec^\intercal\diag(\uvecbase)\umat_i^\intercal = \umat_i\uvecbase$.
			Similarly, Constraint (\ref{eq:proof_lem_spap_eq_aff:constr_dual_constr}) is fulfilled by
			\begin{align*}
				& \betavec'^i+k\alpha^i\unitvec \overset{\text{(a)}}{=}
				\diag(\uvecbase)\umat_i^\intercal + \betavec^i+k\alpha^i\unitvec 
				\\ \overset{\text{(b)}}{\geq} &
				\diag(\uvecbase)\umat_i^\intercal + \left(\left(\umat_i - \dmat_i\pmat\right)\diag(\unitvec-\uvecbase)\right)^\intercal \overset{\text{(c)}}{=} (\umat_i - \dmat_i\pmat')^\intercal,
			\end{align*}
			where (a) and (b) again follow by the definition of $\betavec'^i$ and feasibility of the initial problem and (c) follows by definition of $\pmat'$.
			Finally, $\betavec'^i = \betavec^i+\diag(\uvecbase)\umat_i^\intercal\geq \betavec^i \geq \nullvec$ holds as $\uvecbase$ and $\umat$ are both non-negative.
			Thus we can w.l.o.g. set $\uvecbase = \nullvec$ in Problem (\ref{eq:roc_spap}).
			In doing so, we transform Problem (\ref{eq:roc_spap}) into \mbox{Problem (\ref{eq:roc_aff})}.
			Thus the robust counterparts of affine policies and fully piecewise affine policies with modified right hand side are equivalent.
		\end{proof}
		
		Having shown these results, Proposition \ref{pro:aff_leq_pap} directly follows by combining Lemmas \ref{lem:fpap_leq_pap} and \ref{lem:aff_eq_fpap}.
		Note, that Lemma \ref{lem:piecewise_affine_dual}, which is crucial for the construction of tractable reformulations of fully piecewise affine policies, heavily depends on the linear structure of budgeted uncertainty and the existence of integer optimal solutions.
		Thus the results cannot easily be extended to other settings and in general one cannot hope for tractable reformulations of FPAPs.
		
		\section{Proof of Proposition \ref{pro:pball}}
		\label{sec:proof_pball}
		\begin{proof}
			Analogously to the proofs of Propositions \ref{pro:hypersphere} and \ref{pro:budgeted}, we use Condition (\ref{eq:rot_inv_validity_criterion}) and \mbox{Lemma \ref{lem:rotational_invariate_worst_case}} to find a valid choice of $\mu$ and $\rho$.
			Using the property $\norm{\unitvec}_p \leq 1$ of the $p$-norm ball uncertainty set, we find $j \gamma(j)^p \leq 1$ for all $j\in\{1,\dots,\numu\}$ and thus
			\begin{equation*}
				\gamma(j) = j^{-\frac{1}{p}}.
			\end{equation*}
			Inserting $\gamma(j)$ in Condition (\ref{eq:rot_inv_validity_criterion}), 
			a valid choice of $\rho$ is given by
			\begin{equation}
				\label{eq:proof_pball_rho}
				\max_{j\in\{1,\dots,\numu\}} j \left(j^{-\frac{1}{p}} - \mu\right) \overset{\text{(a)}}{\leq} \mu^{1-p} p^{-p} (p-1)^{p-1} =: \rho,
			\end{equation}
			where (a) follows, as the maximum is taken at $j=(p\mu)^{-p}(p-1)^p$.
			Using the definition of $\U$ and $\vvec_i$ we find $\frac{1}{\beta}\vvec_i\in\U$ to be fulfilled when
			\begin{align*}
				\beta^p =& (\numu-1)\mu^p + (\mu+\rho)^p 
				\overset{\text{(a)}}{\leq} \left(\numu-1+2^{p-1}\right)\mu^p + 2^{p-1}\rho^p \\
				\overset{\text{(b)}}{=}& \left(\numu-1+2^{p-1}\right)\mu^p + 2^{p-1}\mu^{p-p^2} p^{-p^2} (p-1)^{p^2-p}. 
			\end{align*}
			Here (a) follows from Jensen's inequality, as $x^p$ is a convex function for ${x\geq0, p>1}$, and (b) holds by our choice of $\rho$ in Equation (\ref{eq:proof_pball_rho}).
			We now choose $\mu$ to minimize the right-hand side upper bound for $\beta$ and find the minimum to be realized at ${\mu = 2^{\frac{1}{p}}\left(2(m-1)+2^p\right)^{-\frac{1}{p^2}}p^{-1}\left(p-1\right)^{\frac{1}{p} + \left(1-\frac{1}{p}\right)^2}}$.
			Using this $\mu$ in Equation (\ref{eq:proof_pball_rho}), we find the desired $\rho$ and we get the desired upper bound for $\beta$ by taking the $p$\textsuperscript{th} root of
			\begin{align*}
				\beta^p 
				\leq& \left(2(m-1)+2^p\right)^{1-\frac{1}{p}}p^{-p}\left(p-1\right)^{1 + p\left(1-\frac{1}{p}\right)^2} 
				\\&+ \left(2(m-1)+2^p\right)^{1-\frac{1}{p}}p^{p^2-p}\left(p-1\right)^{1-p + (p-p^2)\left(1-\frac{1}{p}\right)^2} p^{-p^2} (p-1)^{p^2-p}\\
				=& \left(2(m-1)+2^p\right)^{1-\frac{1}{p}}p^{-p}\left(\left(p-1\right)^{1 + p\left(1-\frac{1}{p}\right)^2} + \left(p-1\right)^{1-2p + (p-p^2)\left(1-\frac{1}{p}\right)^2 + p^2}\right)\\
				=& \left(2(m-1)+2^p\right)^{1-\frac{1}{p}}p^{-p}\left(\left(p-1\right)^{p-1+\frac{1}{p}} + \left(p-1\right)^{p -2 + \frac{1}{p}}\right)\\
				=& \left(2(m-1)+2^p\right)^{1-\frac{1}{p}}p^{-p}\left(p-1\right)^{p-2+\frac{1}{p}} (1+ (p-1)) \\
				=& \left(2(m-1)+2^p\right)^{1-\frac{1}{p}}p^{1-p}\left(p-1\right)^{p-2+\frac{1}{p}}.
			\end{align*}
		\end{proof}
		
		\section{Proof of Proposition \ref{pro:ellipsoid}}
		\label{sec:proof_ellipsoid}
		\begin{proof}
			We first show that our choice of $\mu, \rho$ yields a valid dominating set using Condition (\ref{eq:rot_inv_validity_criterion}) and \mbox{Lemma \ref{lem:rotational_invariate_worst_case}}.
			Using the properties $\uvec^\intercal \boldsymbol{\Sigma}\uvec \leq 1$ and $\boldsymbol{\Sigma}=\boldsymbol{1}+a(\boldsymbol{J}-\boldsymbol{1})$ of the ellipsoid uncertainty set, we find $(1-a)j \gamma(j)^2 + aj^2\gamma(j)^2 \leq 1$ and thus
			\begin{equation*}
				\gamma(j) = \frac{1}{\sqrt{aj^2+(1-a)j}}.
			\end{equation*}
			Inserting $\gamma(j)$ in Condition (\ref{eq:rot_inv_validity_criterion}) we have to choose $\rho$ such that
			\begin{equation*}
				\rho \geq \max_{j\in\{1,\dots,\numu\}} j \left(\frac{1}{\sqrt{aj^2 + (1-a)j}} - \mu\right).
			\end{equation*}
			For the right-hand side term we find
			\begin{align*}
				&j \left(\frac{1}{\sqrt{aj^2 + (1-a)j}} - \mu\right) \overset{\text{(a)}}{\leq} j \left(\frac{1}{\sqrt{aj^2}}\right) = \frac{1}{\sqrt{a}} \\
				\intertext{and}
				&j \left(\frac{1}{\sqrt{aj^2 + (1-a)j}} - \mu\right) \leq j \left(\frac{1}{\sqrt{(1-a)j}} - \mu\right) \overset{\text{(b)}}{\leq} \frac{1}{4(1-a)\mu},
			\end{align*}
			where for (a) we use $\mu \geq 0, a\leq 1$ and (b) follows as the term is maximized for $j = \frac{1}{4(1-a)\mu^2}$.
			Thus choosing $\rho = \frac{1}{\sqrt{a}}$ or $\rho = \frac{1}{4(1-a)\mu}$ always yields a valid dominating set.
			Using the definition of $\U$ and $\vvec_i$ we find $\frac{1}{\beta}\vvec_i\in\U$ to be fulfilled when
			$$
			\beta^2 = \rho^2 + 2(1-a+a\numu)\rho\mu + (a\numu^2+(1-a)\numu)\mu^2.
			$$
			For $\mu=0, \rho=\frac{1}{\sqrt{a}}$ this results in $\beta = \rho = \frac{1}{\sqrt{a}}$.
			Fixing $\rho = \frac{1}{4(1-a)\mu}$ yields
			\begin{align*}
				\beta^2 = \frac{1}{16(1-a)^2\mu^2} + \frac{1}{2}\left(1+\frac{a}{1-a}\numu\right) + (a\numu^2+(1-a)\numu)\mu^2,
			\end{align*}
			which is minimized for $\mu = \frac{1}{2\sqrt[4]{(1-a)^3\numu + (1-a)^2a\numu^2}}$ resulting in
			$$
			\beta^2 = \frac{1}{2}\left(1 + \frac{1}{1-a}\left(a\numu + \sqrt{(1-a)\numu+a\numu^2}\right)\right).
			$$
			Finally, we close the proof by showing the desired asymptotic approximation bounds.
			The case $a>\numu^{-\frac{2}{3}}$ directly follows from $\beta = \frac{1}{\sqrt{a}} \leq \frac{1}{\sqrt{\numu^{-\frac{2}{3}}}} = O(\numu^{\frac{1}{3}})$ and
			for $a\leq\numu^{-\frac{2}{3}}$ we find
			\begin{align*}
				\beta = 
				& \sqrt{\frac{1}{2}\left(1 + \frac{1}{1-a}\left(a\numu + \sqrt{(1-a)\numu+a\numu^2}\right)\right)}\\
				\overset{\text{(a)}}{\leq}
				& \sqrt{\frac{1}{2}\left(1 + \frac{1}{1-\numu^{-\frac{2}{3}}}\left(\numu^{\frac{1}{3}} + \sqrt{\numu+\numu^{\frac{4}{3}}}\right)\right)} \\
				\overset{\text{(b)}}{=}
				& O(\numu^{\frac{1}{3}}),
			\end{align*}
			where for (a) we use $a\leq\numu^{-\frac{2}{3}}$ and for (b) we use $m\geq2$.
		\end{proof}
		
		\section{Proof of Proposition \ref{pro:general}}
		\label{sec:proof_general}
		We prove Proposition \ref{pro:general} by explicitly constructing a dominating set fulfilling the desired properties.
		To find a good choice of $\rho_1, \dots,\rho_\numu, \vvec_0$ for the construction of $\hat{\U}$ in (\ref{eq:dominating_U}), we use the iterative approach described in \mbox{Algorithm \ref{alg:general_uncertainty}}.
		Intuitively, we increase the base vertex $\vvec_0$, as well as an upper bound for $\rho_i$, in each iteration.
		We bound the approximation factor $\beta$, by bounding the maximal number of iterations in the algorithm.
		Algorithm \ref{alg:general_uncertainty} is a refinement of Algorithm 1 in  Ben-Tal et al. \cite{BenTal2020} and uses a different updating step.
		This modified updating step leads to a less aggressive increase of the base vertex $\vvec_0$ ultimately improving the approximation bound by a factor of $\frac{1}{2}$.
		\begin{algorithm}
			\caption{Generating $\vvec_0$ for general uncertainty sets $\U$}\label{alg:general_uncertainty}
			\begin{algorithmic}[1]
				\State $j = 0$
				\State $\vvec^0 = \nullvec$
				\While{$\max_{\uvec\in\U} \unitvec^\intercal (\uvec-\vvec^j)_+ > j+1$}
				\State $\uvec^j \in \argmax_{\uvec\in\U} \unitvec^\intercal (\uvec-\vvec^j)_+$
				\State $\vvec^{j+1} = \vvec^j + (\uvec^j-\vvec^j)_+$
				\State $j=j+1$
				\EndWhile
				\State\Return $\beta = 2j+1$, $\vvec_0 = \vvec^j$
			\end{algorithmic}
		\end{algorithm}
		
		\begin{lemma}
			Let $\beta, \vvec_0$ be the output of Algorithm \ref{alg:general_uncertainty}.
			Then a solution for $Z_{AR}(\hat{\U})$, where $\hat{\U}$ is constructed according to (\ref{eq:dominating_U}) with $\vvec_0$ and $\rho_i = \frac{\beta+1}{2}$ is a $\beta$ approximation for $Z_{AR}(\U)$ and $\beta\leq2\sqrt{\numu}+1$.
		\end{lemma}
		
		\begin{proof}
			In this proof let $J:=\frac{\beta-1}{2}$ be the index $j$ at termination of Algorithm \ref{alg:general_uncertainty}.
			First, we show that $\hat{\U}$ is a valid dominating set for $\U$ using Criterion (\ref{eq:validity_criterion}).
			Let $\uvec\in\U$. Then
			\begin{align*}
				\sum_{i=1}^{\numu} \frac{\left((\uvec-\vvec_0)_+\right)_i}{\rho_i}
				\overset{\text{(a)}}{=} \frac{1}{J+1} \unitvec^\intercal (\uvec-\vvec_0)_+ \overset{\text{(b)}}{\leq} 1,
			\end{align*}
			where (a) follows from our choice of $\rho_i$ and (b) from the termination criterion of Algorithm \ref{alg:general_uncertainty}.
			Next, we show $\frac{1}{\beta}\vvec_i \in \U$ for all $i\in\{1,\dots,\numu\}$ by
			\begin{align*}
				\frac{1}{\beta}\vvec_i = \frac{J+1}{\beta} \unitvec_i + \frac{1}{\beta} \sum_{j=0}^{J-1} (\uvec^j-\vvec^j)_+
				\leq \frac{1}{\beta}\left((J+1) \unitvec_i + \sum_{j=0}^{J-1} \uvec^j\right)
				\overset{\text{(a)}}{\in}\U.
			\end{align*}
			Here, the containment (a) follows from down-monotinicity, convexity and $\unitvec_i\in\U$.
			Finally, we show that Algorithm \ref{alg:general_uncertainty} terminates after at most $\frac{\beta-1}{2}=J\leq\sqrt{\numu}$ iterations, which implies $\beta\leq2\sqrt{\numu}+1$.
			Let $J':= J-1$ be the index of the last iteration that fulfills the criterion of the loop. Then
			\begin{align*}
				J^2 &= 
				\sum_{j=0}^{J'} J \overset{\text{(a)}}{\leq} 
				\sum_{j=0}^{J'} \unitvec^\intercal (\uvec^{J'} - \vvec^{J'})_+ \overset{\text{(b)}}{\leq} 
				\sum_{j=0}^{J'} \unitvec^\intercal (\uvec^{J'} - \vvec^{j})_+ \\ &\overset{\text{(c)}}{\leq} \sum_{j=0}^{J'} \max_{\uvec\in\U}\unitvec^\intercal (\uvec - \vvec^{j})_+ \overset{\text{(d)}}{=} 
				\sum_{j=0}^{J'} \unitvec^\intercal (\uvec^j - \vvec^j)_+ \\ & =
				\sum_{i=1}^{\numu} \sum_{j=0}^{J'} (\uscal^j_i - v^j_i)_+ \overset{\text{(e)}}{=}
				\sum_{i=1}^{\numu} v^J_i \overset{\text{(f)}}{=} 
				\sum_{i=1}^{\numu} \max_{j\in\left\{0,\dots,J'\right\}} \uscal^j_i \overset{\text{(g)}}{\leq} 
				\sum_{i=1}^{\numu} 1 = \numu.
			\end{align*}
			Here, inequality 
			(a) follows from the termination criterion of Algorithm \ref{alg:general_uncertainty}, 
			(b) follows as $\vvec^j$ only increases through the algorithm and thus $\forall j\leq J' \colon \vvec^{J'}\geq\vvec^j$,
			(c) follows as by $\uvec^{J'}\in\U$ the maximum over $\U$ is an upper bound,
			(d) follows from the choice of $\uvec^j$ in Algorithm \ref{alg:general_uncertainty},
			(e) follows from the construction of $\vvec^j$,
			(f) follows as by induction $v^{j+1}_i = \max(v^j_i, \uscal^j_i) = \max_{j'\leq j}\uscal^{j'}_i$,
			and (g) follows by $\unitvec^j\in\U\subseteq[0,1]^\numu$.
		\end{proof}
		
		\section{Proof of Proposition \ref{pro:stagewise}}
		\label{sec:proof_stagewise}
		\begin{proof}
			For each $t$, let $\numu_t$ be the dimension $\U_t$.
			Then $\numu = \sum_{t=1}^T \numu_t$, and an uncertainty vector $\uvec\in\U$ can be expressed as $\uvec = (\uvec^1,\dots,\uvec^T) = (\uscal^1_{1}, \dots, \uscal^1_{\numu_1}, \dots, \uscal^T_{\numu_T})$, where $\uvec^t\in\U_t$.
			For each $t$, let $\vvec^t_{0}$ be the base vertex and let $\rho^t_{1},\dots,\rho^t_{m_t}$ be the parameters inducing $\hat{\U}_t$ via Construction~(\ref{eq:dominating_U}). 
			Let $\beta = \sum_{t\in\mathcal{T}_1} \beta_t$, and construct the dominating set $\hat{\U}$ for $\U$ with base vertex 
			$
			\vvec'_0 := (\vvec'^1_{0},\dots,\vvec'^T_{0}),
			$
			where
			$$
			\vvec'^t_{0} := \begin{cases}
				\vvec^t_{0} &\text{  if } t\in\mathcal{T}_1\\
				\unitvec &\text{  if } t\in\mathcal{T}_2
			\end{cases}
			$$
			and parameters 
			$$
			\rho'^t_i := \begin{cases}
				\frac{\beta}{\beta_t} \rho^t_{i} &\text{  if } t\in\mathcal{T}_1\\
				0 & \text{  if } t\in\mathcal{T}_2
			\end{cases}
			$$
			using Construction~(\ref{eq:dominating_U}). 
			This construction corresponds to dominating each $\U_t$ with $t\in\mathcal{T}_2$ by the unit vector $\unitvec$. 
			Here $\{\unitvec\}$ is dominating $\U_t$ by Assumption~\ref{ass:uncertainty_assumption}.
			We dominate the remaining $\U_t$ for $t\in\mathcal{T}_1$ with a combined polytope. 
			Then the maximal sum of convex factors for the vertices is given by
			\begin{align*}
				\max_{\uvec\in\U}\sum_{t\in\mathcal{T}_1}\sum_{i=1}^{m_t} \frac{\left((\uvec-\vvec'_0)_+\right)^t_{i}}{\rho'^t_i} 
				&\overset{\text{(a)}}{=} \frac{1}{\beta}\max_{\uvec\in\U}\sum_{t\in\mathcal{T}_1}\beta_t\sum_{i=1}^{m_t} \frac{\left((\uvec^t-\vvec^t_0)_+\right)_{i}}{\rho^t_i}\\
				&\overset{\text{(b)}}{=} \frac{1}{\beta}\sum_{t\in\mathcal{T}_1}\beta_t\max_{\uvec^t\in\U_t}\sum_{i=1}^{m_t} \frac{\left((\uvec^t-\vvec^t_0)_+\right)_{i}}{\rho^t_i}
				\overset{\text{(c)}}{\leq}
				\frac{1}{\beta}\sum_{t\in\mathcal{T}_1}\beta_t
				\overset{\text{(d)}}{=} 1,
			\end{align*}
			where (a) follows from the definition of $\vvec'_0$ and $\rho'^t_i$, 
			(b) follows from $\U$ being stagewise independent,  
			(c) follows from Condition~(\ref{eq:validity_criterion}) and $\hat{\U}_t$ being a dominating set for $\U_t$, and 
			(d) follows from the definition of $\beta$.
			Thus $\hat{\U}$ fulfills Condition~(\ref{eq:validity_criterion}) and is a valid dominating set for $\U$.
			
			Let $\beta' := \max_{t\in\mathcal{T}_2}\beta'_t$.
			To see that $\max(\beta, \beta')$ is indeed an upper bound for the approximation factor, we find 
			\begin{align*}
				&\frac{1}{\max(\beta, \beta')} \vvec'^t_0 
				\overset{\text{(a)}}{\leq} 
				\frac{1}{\beta} (\vvec'^t_0 + \rho'^t_i \unitvec_i) 
				\overset{\text{(b)}}{=} 
				\frac{1}{\beta} (\vvec^t_0 + \frac{\beta}{\beta_t}\rho^t_i \unitvec_i) 
				\overset{\text{(c)}}{\leq} 
				\frac{1}{\beta_t} (\vvec^t_0 + \rho^t_i \unitvec_i) 
				\overset{\text{(d)}}{\in}\U_t & \forall t\in\mathcal{T}_1\\
				&\frac{1}{\max(\beta, \beta')} \vvec'^t_0 
				\overset{\text{(e)}}{\leq} 
				\frac{1}{\beta'} \unitvec 
				\overset{\text{(f)}}{\leq}
				\frac{1}{\beta'_t} \unitvec 
				\overset{\text{(g)}}{\in} \U_t
				& \forall t\in\mathcal{T}_2,
			\end{align*}
			where (a) follows from $\beta\leq\max(\beta, \beta')$ and $\rho'^t_i\geq 0$,
			(b) follows from the definitions of $\vvec'^t_0$ and $\rho'^t_i$, 
			(c) follows from $\beta_t \leq \beta$ for all $t\in\mathcal{T}_1$, 
			(d) follows as $\beta_t$ is a valid approximation factor for $\hat{\U}_t$,
			(e) follows from $\beta'\leq\max(\beta, \beta')$ and the definition of $\vvec'^t_0$,
			(f) follows from $\beta'_t\leq\beta'$ for all $t\in\mathcal{T}_2$, and
			(g) follows from the definition of $\beta'_t$. 
			Using the stagewise structure and down monotinicity of $\U$ this implies 
			\begin{align*}
				&\frac{1}{\max(\beta, \beta')} \vvec'_0 
				=\frac{1}{\max(\beta, \beta')}(\vvec'^1_0, \dots, \vvec'^T_0) \in \U,\\
				&\frac{1}{\max(\beta, \beta')} (\vvec'_0 + \rho'^t_i \unitvec^t_i) 
				=\frac{1}{\max(\beta, \beta')}(\vvec'^1_0, \dots, \vvec'^t_0 + \rho'^t_i \unitvec_i, \dots, \vvec'^T_0) \in \U.
			\end{align*}
		\end{proof}
		
		\section{Proof of Lemma \ref{lem:vertex_shift}}
		\label{sec:proof_vertex_shift}
		\begin{proof}
			Let $\umap\colon\U\to\hat{\U}$ be the domination function of $\hat{\U}$.
			Then define $\umap'\colon\U\to\hat{\U}'$ by
			\begin{equation}
				\label{eq:umap_prime}
				h'_j(\uvec) := h_j(\uvec) + s_j(1-h_j(\uvec)).
			\end{equation}
			This new dominating function $\umap'$ is nonanticipative, as $\umap$ is nonanticipative and the $j$\textsuperscript{th} component of $\umap'$ only depends on the $j$\textsuperscript{th} component of $\umap$.
			To see that (\ref{eq:umap_prime}) is also a domination function, consider the two cases $h_j(\uvec)\geq1$ and $h_j(\uvec)\leq1$.
			\begin{enumerate}
				\item For $h_j(\uvec)\geq1$ the second summand on the right hand side of (\ref{eq:umap_prime}) is negative for all $s_j\in[0,1]$.
				Thus $h'_j(\uvec)$ is minimal for $s_j=1$ and we find ${h'_j(\uvec) \geq h_j(\uvec) + 1(1-h_j(\uvec)) =1}$.
				\item On the other hand, for $h'_j(\uvec)\leq1$ the second summand is positive and $h'_j(\uvec)$ is minimal for $s_j=0$.
				Thus, we have $h'_j(\uvec) \geq h_j(\uvec)$.
			\end{enumerate}
			In any case, we find $\umap'(\uvec)\geq\uvec$ by $h'_j \geq \min(1, h_j(\uvec)) \geq \uscal_j$ where we used $\U\subseteq[0,1]^\numu$ from Assumption \ref{ass:uncertainty_assumption}.
			
			Next, we show $\umap'(\uvec)\in\hat{\U}'$.
			For this let $\umap(\uvec) = \sum_{i=0}^\numu \lambda_i \vvec_i$
			be the convex combination of $\umap(\uvec)\in\hat{\U}$.
			Then,
			$$
			h'_j(\uvec) = h_j(\uvec) + s_j(1-h_j(\uvec)) \overset{\text{(a)}}{=} \sum_{i=0}^\numu \left(\lambda_i v_{ij} + s_j(1-v_{ij})\right) = \sum_{i=0}^\numu \lambda_i v'_{ij}.
			$$
			Note that in (a) we used $\sum_{i=0}^\numu \lambda_i = 1$.
			This shows $\umap'(\uvec) = \sum_{i=0}^\numu \lambda_i \vvec'_i$ and thus $\umap'(\uvec)\in\hat{\U}'$.
		\end{proof}
		
		\section{Proof of Proposition \ref{pro:lift_eq_aff}}
		\label{sec:proof_lift_eq_aff}
		\begin{proof}
			As $Z_{LIFT}\leq Z_{AFF}$ was already shown by Georghiou et al. \cite{Georghiou2015}, we are left to show $Z_{LIFT}\geq Z_{AFF}$.
			We do so by showing that for any optimal affine solution for $Z^L_{AR}(\U^L)$, there is an affine solution for $Z_{AR}(\U)$ with the same objective value.
			First, we define the average lifting operator $\bar{L}\colon \U \to \U^L$ via
			$$
			\bar{L}_{ij}(\uvec) := (z^i_j-z^i_{j-1})\uscal_i.
			$$
			Note, that $R_i(\bar{L}(\uvec)) = \sum_{j=1}^{r_i} (z^i_j-z^i_{j-1})\uscal_i = (z^i_{r_i} - z^i_0)\uscal_i = \uscal_i$, which implies $R(\bar{L}(\uvec)) = \uvec$.
			Thus, $\bar{L}(\uvec)$ fulfills the first condition of (\ref{eq:lifteduncertainty}).
			Next, we verify the second condition by
			$$
			\frac{\bar{L}_{i,j+1}(\uvec)}{z^i_{j+1}-z^i_{j}} = \frac{(z^i_{j+1}-z^i_{j})\uscal_i}{z^i_{j+1}-z^i_{j}} = \uscal_i = \frac{(z^i_{j}-z^i_{j-1})\uscal_i}{z^i_{j}-z^i_{j-1}} = \frac{\bar{L}_{ij}(\uvec)}{z^i_j-z^i_{j-1}}.
			$$
			Finally, the third condition holds by $\bar{L}_{i1} = z^i_1 \uscal_i \leq z^i_1$ which follows from $\uscal_i\leq 1$ in Assumption \ref{ass:uncertainty_assumption}.
			Thus $\bar{L}(\uvec)\in\U^L$.
			Let $\dvec^L$ be an optimal affine solution for $Z^L_{AR}(\U^L)$.
			As $\bar{L}$ is a nonanticipative affine map, the concatenation $\dvec^L\circ\hat{L}$ is also nonanticipative and affine.
			By $\bar{L}(\U)\subseteq\U^L$ the decision $\dvec:=\dvec^L\circ\hat{L}$ is a feasible solution for $Z_{AR}(\U)$.
			Finally, we find
			$$
			Z_{AFF} \leq \max_{\uvec\in\U} \cvec^\intercal \dvec(\uvec) = \max_{\uvec^L\in\bar{L}(\U)} \cvec^\intercal \dvec^L(\uvec^L) \leq \max_{\uvec^L\in\U^L} \cvec^\intercal \dvec^L(\uvec^L) = Z_{LIFT}.
			$$
		\end{proof}
		
		\section{Proof of Proposition \ref{pro:tightenedlifting}}
		\label{sec:proof_tightenedlifting}
		\begin{proof}
			To prove $Z_{TLIFT} \leq Z_{LIFT}$ and $Z_{TLIFT} \leq Z_{SPAP}$ we show that any affine solution for $Z^L_{AR}(\U^L)$ and any picecwise affine solution with re-scaling for $Z(\hat{\U})$ can be transformed to a feasible affine solution for $Z^L_{AR}(\hat{\U}^L)$ with at most the same objective value.
			
			For $Z_{TLIFT} \leq Z_{LIFT}$, let $\dvec^L\colon \U^L \to \mathbb{R}^\numd$ be an optimal affine solution to $Z^L_{AR}(\U^L)$.
			Then $\dvec^L$ is also a feasible affine solution for $Z^L_{AR}(\hat{\U}^L)$ by $\hat{\U}^L\subseteq \U^L$.
			For the objective value, we find
			$$
			Z_{TLIFT} \leq\max_{\uvec^L\in\hat{\U}^L} \cvec^\intercal \dvec^L(\uvec^L) \overset{\text{(a)}}{\leq}
			\max_{\uvec^L\in\U^L} \cvec^\intercal \dvec^L(\uvec^L) = 
			Z_{LIFT},
			$$
			where (a) also follows from $\hat{\U}^L\subseteq \U^L$.
			
			For $Z_{TLIFT} \leq Z_{SPAP}$ let $\hat{\dvec}(\hat{\uvec})$ given by vertex solutions $\dvec_0,\dots, \dvec_\numu$ and scaling factor $\svec$ be a an optimal solution to $Z_{LP}(\hat{\uvec})$ with re-scaling as described in Section \ref{sec:rescaling}.
			We define the map $\umap^L_{\svec}\colon \hat{\U}^L\to \hat{\U}'$ from the tightened lifted uncertainty set $\hat{\U}^L$ defined in Equation (\ref{eq:tightenedlifteduncertainty}) to the re-scaled dominating uncertainty set $\hat{\U}'$ defined in Lemma \ref{lem:vertex_shift} by
			$$
			h^L_{\svec i}(\hat{\uvec}^L) := (1-s_i)(v_{0i} + \hat{\uscal}^L_{i2}) + s_i.
			$$
			This mapping is nonanticipative.
			Furthermore, $\umap^L_{\svec}(\hat{\uvec}^L)$ dominates $R(\hat{\uvec}^L)$ by
			\begin{align*}
				R(\hat{\uvec}^L)_i  = \hat{\uscal}^L_{i1} + \hat{\uscal}^L_{i2} 
				\overset{\text{(a)}}{\leq} (1-s_i)(\hat{\uscal}^L_{i1} + \hat{\uscal}^L_{i2}) + s_i
				\overset{\text{(b)}}{\leq} (1-s_i)(v_{0i} + \hat{\uscal}^L_{i2}) + s_i
				= h^L_{\svec i}(\hat{\uvec}^L).
			\end{align*}
			Here, (a) follows as $\hat{\uscal}^L_{i1} + \hat{\uscal}^L_{i2} \leq 1$ by Assumption \ref{ass:uncertainty_assumption}, and (b) follows from $\hat{\uscal}^L_{i1}\leq v_{0i}$ as $v_{0i}$ is the break-point.
			Using the definition of the re-scaled vertices $\vvec'_i$ from Lemma \ref{lem:vertex_shift}, we find
			\begin{equation}
				\label{eq:proof:tightenedlifting:quasidomination}
				\umap^L_{\svec}(\hat{\uvec}^L) = \sum_{i=1}^\numu \frac{\hat{\uscal}^L_{i2}}{\rho_i} \vvec'_i + \left(1-\sum_{i=1}^\numu \frac{\hat{\uscal}^L_{i2}}{\rho_i}\right) \vvec'_0. 
			\end{equation}
			By the tightening constraint in Definition (\ref{eq:tightenedlifteduncertainty}) this is a valid convex combination for all $\hat{\uvec}^L\in \hat{\U}^L$ and thus $\umap^L_{\svec}\left(\hat{\U}^L\right)\subseteq \hat{\U}'$ as claimed.
			Thus $\dvec^L:=\hat{\dvec}\circ \umap^L_{\svec}$ is a valid solution for $Z_{AR}^L(\hat{\U}^L)$ and
			$$
			\max_{\hat{\uvec}^L\in\hat{\U}^L}\cvec^\intercal \dvec^L(\hat{\uvec}^L) 
			= \max_{\hat{\uvec}\in \umap^L_{\svec}\left(\hat{\U}^L\right)}\cvec^\intercal \hat{\dvec}(\hat{\uvec})
			\leq \max_{\hat{\uvec}\in \hat{\U}}\cvec^\intercal \hat{\dvec}(\hat{\uvec}) = Z_{SPAP}.
			$$
			Using that $\hat{\dvec}$ is given by the vertex solutions $\dvec_0,\dots,\dvec_\numu$ together with (\ref{eq:proof:tightenedlifting:quasidomination}) we find
			$$
			\dvec^L(\hat{\uvec}^L) = \hat{\dvec}\circ \umap^L_{\svec} (\hat{\uvec}^L)
			= \sum_{i=1}^\numu \frac{\hat{\uscal}^L_{i2}}{\rho_i}\dvec_{i} + \left(1 - \sum_{i=1}^\numu \frac{\hat{\uscal}^L_{i2}}{\rho_i}\right) \dvec_0
			= \sum_{i=1}^\numu \left(\frac{\hat{\uscal}^L_{i2}}{\rho_i}\left(\dvec_{i} - \dvec_{0}\right)\right) + \dvec_0,
			$$
			which is affine in $\hat{\U}^L$ and thus
			$$
			Z_{TLIFT} \leq \max_{\hat{\uvec}^L\in\hat{\U}^L}\cvec^\intercal \dvec^L(\hat{\uvec}^L).
			$$
		\end{proof}

		\section{Sampling Uniform Uncertainty Realizations}
		\label{sec:drawing_realizations}
		In this section, we describe how to generate uniform samples from the uncertainty sets used in our numerical experiments efficiently.
		
		\textbf{Hypersphere Uncertainty:}
		To sample uniformly from the hypersphere uncertainty set $\U= \{\uvec\in\mathbb{R}_+^{\numu} \,\vert\, \sqnorm{\uvec}\leq 1\}$, we first draw i.i.d. samples $\uvec^B$ from the $\numu$-dimensional unit ball. 
		It is commonly known, that this can efficiently be done by drawing $\numu$ i.i.d. normal variables that are normalized and rescaled by the radius appropriately, see, e.g. \cite{Devroye1986}. 
		We now map the samples from the unit ball to the positive orthant by taking the absolute value in each component, i.e., we get $\uvec\in\U$ by $\uscal_i := \lvert\uscal^B_i\rvert$. 
		By rotational invariance of the unit ball, this yields a uniform sample on the hypersphere uncertainty set.
		
\begin{figure}
			\centering
			\includegraphics[width = .6\textwidth]{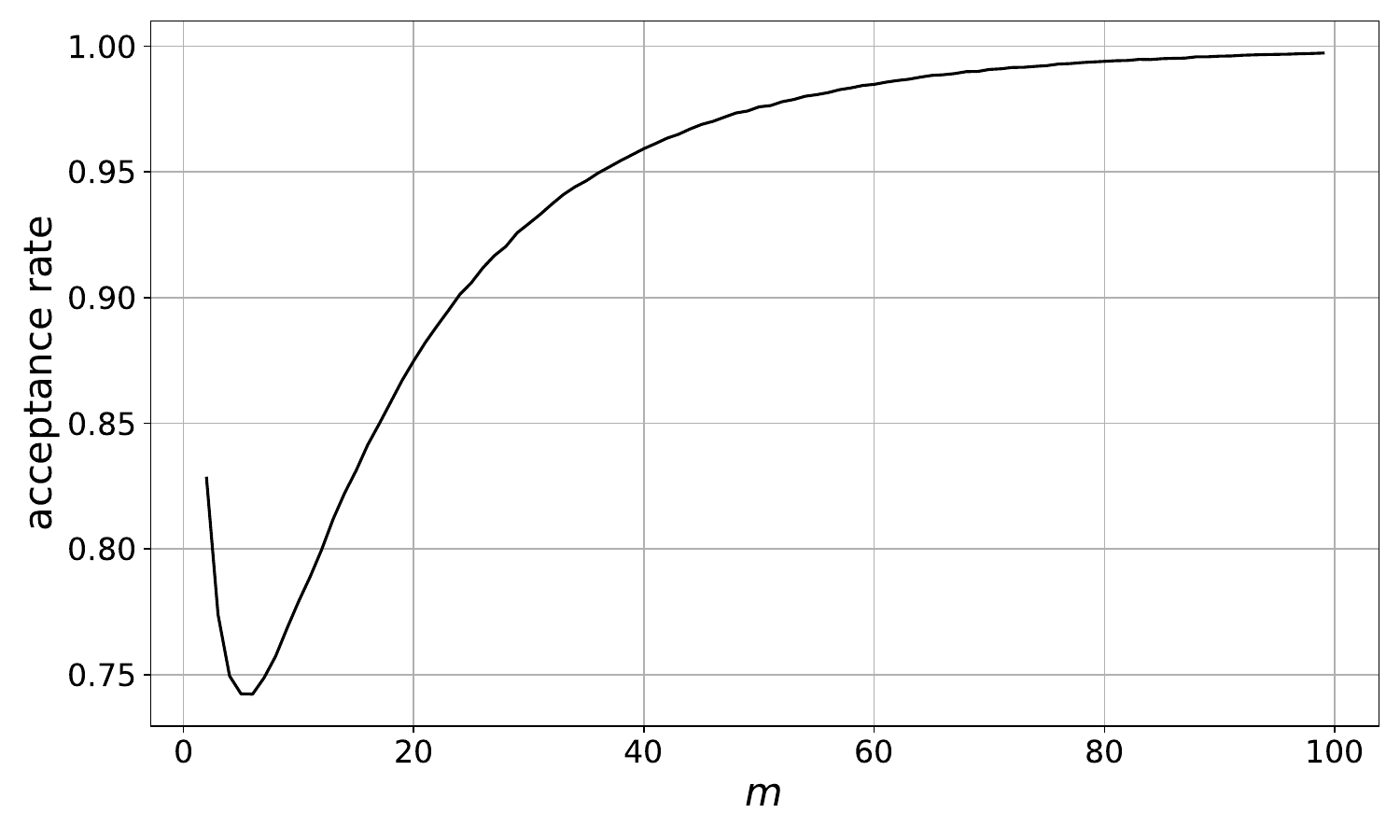}
			\caption{Acceptance rate of our proposed sampling method for budgeted uncertainty with budget $k=\sqrt{\numu}$}
			\label{fig:budeget_acceptance_rate}
		\end{figure}

		\textbf{Budgeted Uncertainty:}
		For the budgeted uncertainty set $\U=\{\uvec\in[0,1]^{\numu} \,\vert\, \norm{\uvec}_1\leq k\}$ with budged $k=\sqrt{\numu}$, we first draw i.i.d. samples $\uvec^S$ from the $\numu$-dimensional unit simplex $\convex(\nullvec, \unitvec_1, \dots, \unitvec_\numu)$.
		This can be done efficiently by taking the first $\numu$ of $\numu + 1$ i.i.d. exponentially distributed variables that were normalized with respect to the one norm, see, e.g. \cite{Devroye1986}. 
		By scaling the samples from the unit simplex with the budget $k$, and rejecting all samples that do not fulfill $\uvec := k \uvec^S \in [0,1]^\numu$, we get a uniform sample on the budgeted uncertainty set. 
		In Figure \ref{fig:budeget_acceptance_rate}, we show the acceptance rate of the samples, which shows that the proposed rejection method is efficient for our application. 

	\end{appendices}
	
	\section*{Declarations}
	
	\paragraph{Competing interests}
	The authors have no conflicts of interest to declare that are relevant to the content of this article.

\end{document}